\documentclass[preprint,11pt]{elsarticle}

 \usepackage{paralist}
  \usepackage{graphics}
  \usepackage{epsfig}

 \usepackage[all]{xy}
\usepackage{graphicx}

\usepackage{amssymb}
\usepackage{amsmath}
\usepackage{amscd}

\newtheorem{theo}{Theorem}

\newtheorem{theorem}[theo]{Theorem}
\newdefinition{definition}[theo]{Definition}
\newtheorem{lemma}[theo]{Lemma}
\newtheorem{proposition}[theo]{Proposition}
\newdefinition{remark}[theo]{Remark}
\newproof{proof}{Proof}

\makeatletter \@addtoreset{equation}{section}
\@addtoreset{theo}{section} \makeatother

\begin{document}

\begin{frontmatter}


\title{Hamilton-Jacobi Theorems for Regular Reducible Hamiltonian Systems on a Cotangent Bundle }

\author{Hong Wang\corref{cor1}}
\ead{hongwang@nankai.edu.cn}
\address{School of Mathematical Sciences and LPMC, Nankai University, Tianjin 300071, P.R.China}
\cortext[cor1]{Corresponding author. Tel.: 0086-022-23501233.
Address: School of Mathematical Sciences, Nankai University, Tianjin
300071, P.R.China.}

\markboth{Hong Wang }{Hamilton-Jacobi Theorems }

\begin{abstract} In this paper, some of formulations of Hamilton-Jacobi
equations for Hamiltonian system and regular reduced Hamiltonian
systems are given. At first, an important lemma is proved, and it is
a modification for the corresponding result of Abraham and Marsden
in \cite{abma78}, such that we can prove two types of geometric
Hamilton-Jacobi theorem for a Hamiltonian system on the cotangent
bundle of a configuration manifold, by using the symplectic
form and dynamical vector field. Then these results are
generalized to the regular reducible Hamiltonian system with
symmetry and momentum map, by using the reduced symplectic form
and the reduced dynamical vector field. The Hamilton-Jacobi theorems
are proved and two types of Hamilton-Jacobi equations, for the
regular point reduced Hamiltonian system and the regular orbit
reduced Hamiltonian system, are obtained. As an application of the
theoretical results, the regular point reducible Hamiltonian system
on a Lie group is considered, and two types of Lie-Poisson
Hamilton-Jacobi equation for the regular point reduced system are
given. In particular, the Type I and Type II of Lie-Poisson
Hamilton-Jacobi equations for the regular point reduced rigid body
and heavy top systems are shown, respectively.
\end{abstract}


\end{frontmatter}

\tableofcontents

\section{Introduction}

Symmetry is a general phenomenon in the natural world, but it is
widely used in the study of mathematics and mechanics. The reduction
theory for mechanical system with symmetry has its origin in the
classical work of Euler, Lagrange, Hamilton, Jacobi, Routh,
Liouville and Poincar\'{e} and its modern geometric formulation in
the general context of symplectic manifolds and equivariant momentum
maps is developed by Meyer, Marsden and Weinstein; see Abraham and
Marsden \cite{abma78} or Marsden and Weinstein \cite{mawe74} and
Meyer \cite{me73}. The main goal of reduction theory in mechanics is
to use conservation laws and the associated symmetries to reduce the
number of dimensions of a mechanical system required to be
described. So, such reduction theory is regarded as a useful tool
for simplifying and studying concrete mechanical systems.
Hamiltonian reduction theory is one of the most active subjects in
the study of modern analytical mechanics and applied mathematics, in
which a lot of deep and beautiful results have been obtained, see
the studies given by Abraham and Marsden \cite{abma78}, Arnold
\cite{ar89}, Marsden et al. \cite{ma92, mamiorpera07, mara99,
mawe74}, Ortega and Ratiu \cite{orra04}, Libermann and Marle \cite
{lima87}, Le\'{o}n and Rodrigues \cite{lero89} etc. on regular point
reduction and regular orbit reduction, singular point reduction and
singular orbit reduction, optimal reduction and reduction by stages
for Hamiltonian systems and so on; and there is still much to be
done in this subject.\\

At the same time, we note also that the well-known Hamilton-Jacobi
theory is an important part of classical mechanics. On the one hand,
Hamilton-Jacobi equation provides a characterization of the
generating functions of certain time-dependent canonical
transformations, such that a given Hamiltonian system in such a form
that its solutions are extremely easy to find by reduction to the
equilibrium, see Abraham and Marsden \cite{abma78}, Arnold
\cite{ar89} and Marsden and Ratiu \cite{mara99}. On the other hand,
it is possible in many cases that Hamilton-Jacobi equation provides
an immediate way to integrate the equation of motion of system, even
when the problem of Hamiltonian system itself has not been or cannot
be solved completely. In addition, the Hamilton-Jacobi equation is
also fundamental in the study of the quantum-classical relationship
in quantization, and it also plays an important role in the
development of numerical integrators that preserve the symplectic
structure and in the study of stochastic dynamical systems, see
Woodhouse \cite{wo92}, Ge and Marsden \cite{gema88}, Marsden and
West \cite{mawe01} and L\'{a}zaro-Cam\'{i} and Ortega \cite{laor09}.
For these reasons Hamilton-Jacobi theory is described as a useful
tool in the study of Hamiltonian system theory, and has been
extensively developed in past many years. We note that some
beautiful results have been obtained, see Cari$\tilde{n}$ena et al.
\cite{cagrmamamuro06} and \cite{cagrmamamuro10}, Iglesias et al.
\cite{iglema08}, for more details.\\

Now, it is a natural problem how to study the Hamilton-Jacobi theory
for a variety of reduced Hamiltonian systems by combining with
reduction theory and Hamilton-Jacobi theory of Hamiltonian systems.
This is a goal of our research. In this paper, some of formulations of Hamilton-Jacobi
equations for Hamiltonian system and regular reduced Hamiltonian systems
are given, and the main contributions are as follows:
(1) We prove a key lemma, which is an important tool
for proofs of the following theorems;
(2) We prove two types of geometric Hamilton-Jacobi
theorem for a Hamiltonian system on the cotangent bundle of a
configuration manifold, by using the symplectic form and
dynamical vector field; (3) We generalize the above results to the
regular reducible Hamiltonian system with symmetry, and obtain two types of
Hamilton-Jacobi equations for the regular point reduced Hamiltonian
system and the regular orbit reduced Hamiltonian system, see Theorem 3.3, Theorem 3.4,
Theorem 4.2 and Theorem 4.3, by using the reduced symplectic forms and the reduced
dynamical vector fields; It is worthy of note that the regular reduced symplectic spaces of the regular orbit
reduced Hamiltonian system and the regular point reduced Hamiltonian system are different,
and the symplectic forms on the reduced spaces are also different.
Thus, the assumption conditions in Theorem 4.2 and Theorem 4.3
are different from the assumption conditions in Theorem 3.3 and Theorem 3.4,
which depend on the precise analysis
of the geometric structures of the regular orbit reduced space.
(4) As an application, we give two types of
Lie-Poisson Hamilton-Jacobi equation for the regular point reduced
Hamiltonian system on a Lie group, and show the Type I and Type II of Lie-Poisson Hamilton-Jacobi
equations for the regular point reduced rigid body and heavy top systems,
respectively. In general, we
know that it is not easy to find the solutions of Hamilton's
equation. But, if we can get a solution of Hamilton-Jacobi equation
for a Hamiltonian system, by using the relationship between
Hamilton's equation and Hamilton-Jacobi equation, it is easy to give
a special solution of Hamilton's equation. Thus, it is very
important to give explicitly the various formulations of
Hamilton-Jacobi equations for Hamiltonian system and the reduced Hamiltonian systems.\\

A brief of outline of this paper is as follows. In the second
section, we first prove a key lemma, which is obtained by a careful
modification for the corresponding results of Abraham and Marsden in \cite{abma78}.
Then we prove two types of geometric version of Hamilton-Jacobi theorem of a Hamiltonian
system on the cotangent bundle of a configuration manifold, by using
the symplectic form and the dynamical vector field.
In the third section and the fourth section, we discuss the regular reducible Hamiltonian
systems with symmetry and momentum map, by combining with the Hamilton-Jacobi theory
and the regular symplectic reduction theory. The two types of Hamilton-Jacobi
equations for the regular point and the regular
orbit reduced Hamiltonian systems are obtained, respectively, by using the reduced symplectic
forms and the reduced dynamical vector fields. As the applications of
the theoretical results, in the fifth section, the regular
point reducible Hamiltonian system on a Lie group is considered, and
two types of Lie-Poisson Hamilton-Jacobi equation for
the regular point reduced system are given. In particular,
the Type I and Type II of Lie-Poisson Hamilton-Jacobi
equations for the regular point reduced rigid body and heavy top systems are shown,
respectively. These research works develop the reduction and
Hamilton-Jacobi theory of a Hamiltonian system with symmetry and make
us have much deeper understanding and recognition for the structures
of Hamiltonian systems.

\section{Geometric Hamilton-Jacobi Theorem of Hamiltonian System}

In this section, we first review briefly some basic facts about
Hamilton-Jacobi theory, and state our idea to study the problem in
this paper. Then we prove a key lemma, which is an important tool
for the proofs of geometric Hamilton-Jacobi theorems of Hamiltonian
system and the regular reducible Hamiltonian system with symmetry.
Finally, we prove two types of geometric version of Hamilton-Jacobi
theorem of a Hamiltonian system on the cotangent bundle of a
configuration manifold, by using the symplectic form and
dynamical vector field. It is worthy of note that we describe the
Hamilton-Jacobi equation by Hamiltonian vector field of the system,
it is easy to be generalized to the cases of the regular reduced
Hamiltonian systems. We shall follow the notations and conventions
introduced in Abraham and Marsden \cite{abma78}, Marsden and Ratiu
\cite{mara99}, Ortega and Ratiu \cite{orra04}, and Marsden et al.
\cite{mawazh10}. In this paper, we assume that all manifolds are
real, smooth and finite dimensional
and all actions are smooth left actions.\\

It is well-known that Hamilton-Jacobi theory from the variational
point of view is originally developed by Jacobi in 1866, which state
that the integral of Lagrangian of a system along the solution of
its Euler-Lagrange equation satisfies the Hamilton-Jacobi equation.
The classical description of this problem from the geometrical point
of view is given by Abraham and Marsden in \cite{abma78} as follows:
Let $Q$ be a smooth manifold and $TQ$ the tangent bundle, $T^* Q$
the cotangent bundle with the canonical symplectic form $\omega$£¬
and the projection $\pi_Q: T^* Q \rightarrow Q $ induces the map $
T\pi_{Q}: TT^* Q \rightarrow TQ. $

\begin{theorem}
Assume that the triple $(T^*Q,\omega,H)$ is a Hamiltonian system
with Hamiltonian vector field $X_H$, and $W: Q\rightarrow
\mathbb{R}$ is a given function. Then the following two assertions
are equivalent:

\noindent $(\mathrm{i})$ For every curve $\sigma: \mathbb{R}
\rightarrow Q $ satisfying $\dot{\sigma}(t)= T\pi_Q
(X_H(\mathbf{d}W(\sigma(t))))$, $\forall t\in \mathbb{R}$, then
$\mathbf{d}W \cdot \sigma $ is an integral curve of the Hamiltonian
vector field $X_H$.\\

\noindent $(\mathrm{ii})$ $W$ satisfies the Hamilton-Jacobi equation
$H(q^i,\frac{\partial W}{\partial q^i})=E, $ where $E$ is a
constant.
\end{theorem}

It is worthy of note that if we take that $\gamma=\mathbf{d}W$ in
the above theorem, then $\gamma$ is a closed one-form on $Q$, and
the equation $\mathbf{d}(H \cdot \mathbf{d}W)=0$ is equivalent to
the Hamilton-Jacobi equation $H(q^i,\frac{\partial W}{\partial
q^i})=E$, where $E$ is a constant. This result is used the
formulation of a geometric version of Hamilton-Jacobi theorem for
Hamiltonian system, see Cari\~{n}ena et al. \cite{cagrmamamuro06}
and Iglesias et al. \cite{iglema08}. On the other hand, this result
is developed in the context of time-dependent Hamiltonian system by
Marsden and Ratiu in \cite{mara99}. The Hamilton-Jacobi equation may
be regarded as a nonlinear partial differential equation for someone
generating function $S$, and the problem is become how to choose a
time-dependent canonical transformation $\Psi: T^*Q\times \mathbb{R}
\rightarrow T^*Q\times \mathbb{R}, $ which transforms the dynamical
vector field of a time-dependent Hamiltonian system to equilibrium,
such that the generating function $S$ of $\Psi$ satisfies the
time-dependent Hamilton-Jacobi equation, that is, the dynamical
vector field is degenerate along the solution of Hamilton-Jacobi
equation. In particular, for the time-independent Hamiltonian
system, we may look for a symplectic map as the canonical
transformation. This work offers an important idea that one can use
the dynamical vector field of a Hamiltonian system to describe
Hamilton-Jacobi equation. Moreover, we also hope to use the
dynamical vector fields of the regular reduced Hamiltonian systems
to describe the Hamilton-Jacobi equations for the regular reduced
Hamiltonian systems. These are the main works in this paper. In
order to do these, we need first to give two types of formulation of
Hamilton-Jacobi theorem for a Hamiltonian system on the cotangent
bundle of a configuration manifold. Thus, in the following we first
give an important notion and prove a key lemma, and this lemma is an
important tool for the proofs of two types of geometric
Hamilton-Jacobi theorem of the Hamiltonian system.\\

Denote by $\Omega^i(Q)$ the set of all i-forms on $Q$, $i=1,2.$
For any $\gamma \in \Omega^1(Q),\; q\in Q, $ then $\gamma(q)\in T_q^*Q, $
and we can define a map $\gamma: Q \rightarrow T^*Q, \; q \rightarrow (q, \gamma(q)).$
Hence we say often that the map $\gamma: Q
\rightarrow T^*Q$ is an one-form on $Q$. If the one-form $\gamma$ is closed,
then $\mathbf{d}\gamma(x,y)=0, \; \forall\;
x, y \in TQ$. In the following we give a weaker notion.
\begin{definition}
The one-form $\gamma$ is called to be closed with respect to $T\pi_{Q}:
TT^* Q \rightarrow TQ, $ if for any $v, w \in TT^* Q, $ we have
$\mathbf{d}\gamma(T\pi_{Q}(v),T\pi_{Q}(w))=0. $
\end{definition}

From the above definition we know that, if $\gamma$ is a closed one-form,
then it must be closed with respect to $T\pi_{Q}: TT^* Q \rightarrow
TQ. $ Conversely, if $\gamma$ is closed with respect to
$T\pi_{Q}: TT^* Q \rightarrow TQ, $ then it may not be closed. We can
prove a general result as follows.

\begin{proposition}
Assume that $\gamma: Q \rightarrow T^*Q$ is an one-form on $Q$ and
it is not closed. we define the set $N$, which is a subset of $TQ$,
such that the one-form $\gamma$ on $N$ satisfies the condition that
for any $x,y \in N, \; \mathbf{d}\gamma(x,y)\neq 0. $ Denote by
$Ker(T\pi_Q)= \{u \in TT^*Q| \; T\pi_Q(u)=0 \}, $ and $T\gamma: TQ
\rightarrow TT^* Q .$ If $T\gamma(N)\subset Ker(T\pi_Q), $ then
$\gamma$ is closed with respect to $T\pi_{Q}: TT^* Q \rightarrow TQ.
$\end{proposition}

\noindent{\bf Proof: } In fact, for any $v, w \in TT^* Q, $ if
$T\pi_{Q}(v) \notin N, $ or $T\pi_{Q}(w))\notin N, $ then by the
definition of $N$, we know that
$\mathbf{d}\gamma(T\pi_{Q}(v),T\pi_{Q}(w))=0; $ If $T\pi_{Q}(v)\in
N, $ and $T\pi_{Q}(w))\in N, $ from the condition $T\gamma(N)\subset
Ker(T\pi_Q), $ we know that $T\pi_{Q}\cdot T\gamma \cdot
T\pi_{Q}(v)= T\pi_{Q}(v)=0, $ and $T\pi_{Q}\cdot T\gamma \cdot
T\pi_{Q}(w)= T\pi_{Q}(w)=0, $ where we have used the
relation $\pi_Q\cdot \gamma\cdot \pi_Q= \pi_Q, $ and hence
$\mathbf{d}\gamma(T\pi_{Q}(v),T\pi_{Q}(w))=0. $ Thus, for any $v, w
\in TT^* Q, $ we have always that
$\mathbf{d}\gamma(T\pi_{Q}(v),T\pi_{Q}(w))=0, $ that is, $\gamma$ is
closed with respect to $T\pi_{Q}: TT^* Q \rightarrow TQ. $
\hskip 0.3cm $\blacksquare$\\

Now, we prove the following Lemma 2.4. It is worthy of note that
this lemma is obtained by a careful modification for the
corresponding result of Abraham and Marsden in \cite{abma78}.

\begin{lemma}
Assume that $\gamma: Q \rightarrow T^*Q$ is an one-form on $Q$, and
$\lambda=\gamma \cdot \pi_{Q}: T^* Q \rightarrow T^* Q .$ Then
we have that the following two assertions hold.\\
\noindent $(\mathrm{i})$ For any $x, y \in TQ, \;
\gamma^*\omega(x,y)= -\mathbf{d}\gamma (x,y),$ and for any $v, w \in
TT^* Q, \; \lambda^*\omega(v,w)=\\ -\mathbf{d}\gamma(T\pi_{Q}(v), \;
T\pi_{Q}(w)),$
since $\omega$ is the canonical symplectic form on $T^*Q$; \\
\noindent $(\mathrm{ii})$ For any $v, w \in TT^* Q, \;
\omega(T\lambda \cdot v,w)= \omega(v, w-T\lambda \cdot
w)-\mathbf{d}\gamma(T\pi_{Q}(v), \; T\pi_{Q}(w)). $
\end{lemma}

\noindent{\bf Proof:} We first prove the assertion $(\mathrm{i})$.
Since $\omega$ is the canonical symplectic form on $T^*Q$, we know
that there is an unique canonical one-form $\theta$, such that
$\omega= -\mathbf{d} \theta. $ From the Proposition 3.2.11 in
Abraham and Marsden \cite{abma78}, we have that for the one-form
$\gamma: Q \rightarrow T^*Q, \; \gamma^* \theta= \gamma. $ Then we
can obtain that
\begin{align*}
\gamma^*\omega(x,y) = \gamma^* (-\mathbf{d} \theta) (x, y) =
-\mathbf{d}(\gamma^* \theta)(x, y)= -\mathbf{d}\gamma (x, y).
\end{align*}
Note that $\lambda=\gamma \cdot \pi_{Q}: T^* Q \rightarrow T^* Q, $
and $\lambda^*= \pi_{Q}^* \cdot \gamma^*: T^*T^* Q \rightarrow
T^*T^* Q, $ then we have that
\begin{align*}
\lambda^*\omega(v,w) &= \lambda^* (-\mathbf{d} \theta) (v, w)
=-\mathbf{d}(\lambda^* \theta)(v, w)= -\mathbf{d}(\pi_{Q}^* \cdot
\gamma^* \theta)(v, w)\\ &= -\mathbf{d}(\pi_{Q}^* \cdot\gamma )(v,
w)= -\mathbf{d}\gamma(T\pi_{Q}(v), \; T\pi_{Q}(w)).
\end{align*}
It follows that the assertion $(\mathrm{i})$ holds.\\

Next, we prove the assertion $(\mathrm{ii})$. For any $v, w \in TT^*
Q,$ note that $v- T(\gamma \cdot \pi_Q)\cdot v$ is vertical, because
$$
T\pi_Q(v- T(\gamma \cdot \pi_Q)\cdot v)=T\pi_Q(v)-T(\pi_Q\cdot
\gamma\cdot \pi_Q)\cdot v= T\pi_Q(v)-T\pi_Q(v)=0,
$$
where we have used the relation $\pi_Q\cdot \gamma\cdot \pi_Q= \pi_Q. $
Thus, $\omega(v- T(\gamma \cdot \pi_Q)\cdot v,w- T(\gamma \cdot
\pi_Q)\cdot w)= 0, $ and hence,
$$\omega(T(\gamma \cdot \pi_Q)\cdot v, \; w)=
\omega(v, \; w-T(\gamma \cdot \pi_Q)\cdot w)+ \omega(T(\gamma \cdot
\pi_Q)\cdot v, \; T(\gamma \cdot \pi_Q)\cdot w). $$ However, the
second term on the right-hand side is given by
$$
\omega(T(\gamma \cdot \pi_Q)\cdot v, \; T(\gamma \cdot \pi_Q)\cdot
w)= \gamma^*\omega(T\pi_Q(v), \; T\pi_Q(w))=
-\mathbf{d}\gamma(T\pi_{Q}(v), \; T\pi_{Q}(w)),
$$
where we have used the assertion $(\mathrm{i})$. It follows that
\begin{align*}
\omega(T\lambda \cdot v,w) &=\omega(T(\gamma \cdot \pi_Q)\cdot v, \;
w)\\ &= \omega(v, \; w-T(\gamma \cdot \pi_Q)\cdot w)-\mathbf{d}\gamma(T\pi_{Q}(v), \; T\pi_{Q}(w))
\\ &= \omega(v,
w-T\lambda \cdot w)-\mathbf{d}\gamma(T\pi_{Q}(v), \; T\pi_{Q}(w)).
\end{align*}
Thus, the assertion $(\mathrm{ii})$ holds.
\hskip 0.3cm $\blacksquare$\\

Now, for a given Hamiltonian system $(T^*Q,\omega,H)$, by using
the above Lemma 2.4, we can prove the following two types of geometric
Hamilton-Jacobi theorem for the Hamiltonian system.
At first, by using the fact that the one-form $\gamma: Q
\rightarrow T^*Q $ is closed with respect to
$T\pi_Q: TT^* Q \rightarrow TQ, $ we can prove the Type I of geometric
Hamilton-Jacobi theorem for the Hamiltonian system. For convenience,
the maps involved in the following theorem and its proof are shown
in Diagram-1.

\begin{center}
\hskip 0cm \xymatrix{ & T^* Q \ar[r]^{\pi_Q}
& Q \ar[d]_{X_H^\gamma} \ar[r]^{\gamma} & T^*Q \ar[d]^{X_H} \\
  & T(T^*Q) & TQ \ar[l]_{T\gamma} & T(T^* Q)\ar[l]_{T\pi_Q}}
\end{center}
$$\mbox{Diagram-1}$$

\begin{theorem} (Type I of Hamilton-Jacobi Theorem for a Hamiltonian System)
For the Hamiltonian system $(T^*Q,\omega,H)$, assume that $\gamma: Q
\rightarrow T^*Q$ is an one-form on $Q$, and $X_H^\gamma = T\pi_{Q}\cdot X_H \cdot \gamma$,
where $X_{H}$ is the dynamical vector field of $(T^*Q,\omega,H)$.
If the one-form $\gamma: Q \rightarrow T^*Q $ is closed with respect to
$T\pi_Q: TT^* Q \rightarrow TQ, $ then $\gamma$ is a solution of the equation
$T\gamma\cdot X_H^\gamma= X_H\cdot \gamma ,$ which is called the Type I of
Hamilton-Jacobi equation for the Hamiltonian system $(T^*Q,\omega,H)$.
\end{theorem}

\noindent{\bf Proof: } If we take that $v= X_H\cdot \gamma \in TT^* Q, $ and for
any $w \in TT^* Q, \; T\pi_{Q}(w)\neq 0, $ from Lemma 2.4(ii) we have that
\begin{align*}
\omega(T\gamma \cdot X_H^\gamma, \; w) &= \omega(T(\gamma \cdot
\pi_Q)\cdot X_H\cdot \gamma, \; w)\\ &= \omega(X_H\cdot \gamma, \;
w-T(\gamma \cdot \pi_Q)\cdot
w)-\mathbf{d}\gamma(T\pi_{Q}(X_H\cdot \gamma), \; T\pi_{Q}(w))\\
& =\omega(X_H\cdot \gamma, \; w) - \omega(X_H\cdot \gamma, \;
T\lambda \cdot w)-\mathbf{d}\gamma(T\pi_{Q}(X_H\cdot \gamma), \; T\pi_{Q}(w)).
\end{align*}
Because the one-form $\gamma: Q \rightarrow T^*Q $ is closed with respect to
$T\pi_Q: TT^* Q \rightarrow TQ, $ then we have that
$$
\mathbf{d}\gamma(T\pi_{Q}(X_H\cdot \gamma), \; T\pi_{Q}(w))=0,
$$
and hence
\begin{equation}
\omega(T\gamma \cdot X_H^\gamma, \; w)- \omega(X_H\cdot \gamma, \; w)
= -\omega(X_H\cdot \gamma, \; T\lambda \cdot w).
\end{equation}
If $\gamma$ satisfies the equation $T\gamma\cdot X_H^\gamma= X_H\cdot \gamma ,$
from Lemma 2.4(i) we can obtain that
\begin{align*}
-\omega(X_H\cdot \gamma, \; T\lambda \cdot w) &
= -\omega(T\gamma \cdot X_H^\gamma, \; T\lambda \cdot w)\\
& =-\omega(T\gamma \cdot T\pi_{Q} \cdot X_{H}\cdot\gamma, \; T\lambda \cdot w)
=-\omega(T\lambda \cdot X_{H}\cdot\gamma, \; T\lambda \cdot w)\\
& = -\lambda^*\omega( X_{H}\cdot\gamma, \; w)=
\textbf{d}\gamma(T\pi_{Q}( X_{H}\cdot\gamma ), \; T\pi_{Q}(w))=0.
\end{align*}
But, because the symplectic form $\omega$ is non-degenerate, the left side of (2.1) equals zero, only when
$\gamma$ satisfies the equation $T\gamma\cdot X_H^\gamma= X_H\cdot \gamma .$ Thus,
if the one-form $\gamma: Q \rightarrow T^*Q $ is closed with respect to
$T\pi_Q: TT^* Q \rightarrow TQ, $ then $\gamma$ must be a solution of
the Type I of Hamilton-Jacobi equation
$T\gamma\cdot X_H^\gamma= X_H\cdot \gamma .$
\hskip 0.3cm $\blacksquare$\\

Next, for any symplectic map $\varepsilon: T^* Q \rightarrow T^* Q $,
we can prove the following Type II of geometric
Hamilton-Jacobi theorem for the Hamiltonian system. For convenience,
the maps involved in the following theorem and its proof are shown
in Diagram-2.

\begin{center}
\hskip 0cm \xymatrix{ & T^* Q \ar[r]^{\varepsilon}
& T^*Q \ar[d]_{X_{H\cdot \varepsilon}}
\ar[dr]^{X_H^\varepsilon} \ar[r]^{\pi_Q}
& Q \ar[r]^{\gamma} & T^*Q \ar[d]^{X_H} \\
&  & T(T^*Q) & TQ \ar[l]_{T\gamma} & T(T^* Q)\ar[l]_{T\pi_Q}}
\end{center}
$$\mbox{Diagram-2}$$

\begin{theorem} (Type II of Hamilton-Jacobi Theorem for a Hamiltonian System)
For the Hamiltonian system $(T^*Q,\omega,H)$, assume that $\gamma: Q
\rightarrow T^*Q$ is an one-form on $Q$, and
$\lambda=\gamma\cdot\pi_{Q}: T^* Q \rightarrow T^* Q $, and for any
symplectic map $\varepsilon: T^* Q \rightarrow T^* Q $, denote by
$X_H^\varepsilon = T\pi_{Q}\cdot X_H \cdot \varepsilon$,
where $X_{H}$ is the dynamical vector field of $(T^*Q,\omega,H)$.
Then $\varepsilon$ is a solution of the equation
$T\varepsilon\cdot X_{H\cdot\varepsilon}= T\lambda \cdot X_H \cdot \varepsilon,$
if and only if it is a solution of the equation $T\gamma\cdot X_H^\varepsilon= X_H\cdot
\varepsilon, $ where $ X_{H\cdot\varepsilon} \in
TT^*Q $ is the Hamiltonian vector field of the function $H\cdot\varepsilon:
T^*Q\rightarrow \mathbb{R}. $
The equation $T\gamma\cdot X_H^\varepsilon= X_H\cdot
\varepsilon ,$ is called the Type II of Hamilton-Jacobi equation
for the Hamiltonian system $(T^*Q,\omega,H)$.
\end{theorem}

\noindent{\bf Proof: } If we take that $v= X_H\cdot \varepsilon \in TT^* Q, $ and for
any $w \in TT^* Q, \; T\lambda(w)\neq 0, $ from Lemma 2.4 we have that
\begin{align*}
\omega(T\gamma \cdot X_H^\varepsilon, \; w) &= \omega(T(\gamma \cdot
\pi_Q)\cdot X_H\cdot \varepsilon, \; w)\\ &= \omega(X_H\cdot \varepsilon, \;
w-T(\gamma \cdot \pi_Q)\cdot
w)-\mathbf{d}\gamma(T\pi_{Q}(X_H\cdot \varepsilon), \; T\pi_{Q}(w))\\
& =\omega(X_H\cdot \varepsilon, \; w) - \omega(X_H\cdot \varepsilon, \;
T\lambda \cdot w)+\lambda^*\omega(X_H\cdot \varepsilon, \; w)\\
& =\omega(X_H\cdot \varepsilon, \; w) - \omega(X_H\cdot \varepsilon, \;
T\lambda \cdot w)+ \omega(T\lambda \cdot X_H\cdot \varepsilon, \; T\lambda \cdot w).
\end{align*}
Note that $\varepsilon: T^* Q
\rightarrow T^* Q $ is symplectic, and hence $ X_H\cdot \varepsilon=
T\varepsilon \cdot X_{H\cdot\varepsilon}, $ along $\varepsilon$. From the above arguments, we can obtain that
\begin{align*}
&\omega(T\gamma \cdot X_H^\varepsilon, \; w)- \omega(X_H\cdot \varepsilon, \; w)\\
& =-\omega(T\varepsilon \cdot X_{H\cdot\varepsilon}, \; T\lambda \cdot w)+ \omega(T\lambda \cdot X_H\cdot \varepsilon, \; T\lambda \cdot w)\\
& = \omega(T\lambda \cdot X_H\cdot \varepsilon -T\varepsilon \cdot X_{H\cdot\varepsilon}, \; T\lambda \cdot w).
\end{align*}
Because the symplectic form $\omega$ is non-degenerate,
it follows that $T\gamma\cdot X_H^\varepsilon= X_H\cdot
\varepsilon ,$ is equivalent to $T\varepsilon \cdot X_{H\cdot\varepsilon} = T\lambda\cdot X_H\cdot \varepsilon $.
Thus, $\varepsilon$ is a solution of the equation
$T\varepsilon\cdot X_{H\cdot\varepsilon}= T\lambda \cdot X_H \cdot\varepsilon,$ if and only if it is a solution of
the Type II of Hamilton-Jacobi equation $T\gamma\cdot X_H^\varepsilon= X_H\cdot
\varepsilon .$
\hskip 0.3cm $\blacksquare$\\

In the following we shall state the relationship between the Type I and Type II of
Hamilton-Jacobi equation and the classical Hamilton-Jacobi equation,
from the view point of generating function of a symplectic map.
At first, the following Proposition 2.7 is the Proposition 5.2.1
given by Abraham and Marsden in \cite{abma78}.

\begin{proposition}
Let $(P_i, \omega_i), \; i=1,2, $ be two symplectic manifolds,
and $\pi_i: P_1\times P_2 \rightarrow P_i $
the projection onto $P_i, \; i=1,2, $ and $\Omega= \pi_1^*\omega_1- \pi_2^*\omega_2. $
Then we have that \\
\noindent $(\mathrm{i})$ the $\Omega$ is a symplectic form on $P_1\times P_2; $\\
\noindent $(\mathrm{ii})$ a map $f: P_1 \rightarrow P_2$ is symplectic
if and only if $\mathbf{i}^*_f\Omega =0, $
where $\mathbf{i}_f: \Gamma_f \rightarrow P_1\times P_2 $ is inclusion and $\Gamma_f$
is the graph of $f$, that is, $\Gamma_f= \{(x, f(x))\in P_1\times P_2 | \; \forall x \in P_1\}.$
(In fact, $\Gamma_f$ is a Lagrangian submanifold of $P_1\times P_2$.)
\end{proposition}

Assume that $\theta_i $ is the canonical one-form of $P_i, \;
i=1,2,$ and the canonical symplectic forms $\omega_i =
-\mathbf{d}\theta_i, \; i=1,2. $ Then $\Theta= \pi_1^*\theta_1-
\pi_2^*\theta_2, $ and locally, $\Omega= -\mathbf{d}\Theta. $ Thus,
$$\mathbf{i}^*_f\Omega= -\mathbf{i}^*_f\mathbf{d}\Theta= -\mathbf{d}\mathbf{i}^*_f\Theta =0,$$
that is, $\mathbf{i}^*_f\Theta$ being close is equivalent to $f$ being symplectic.
Locally, by the Poincar\'{e} lemma, $\mathbf{i}^*_f\Theta= -\mathbf{d}S$ for
a function $S: \Gamma_f \rightarrow \mathbb{R}. $
Such a function $S$ is called a {\bf generating function}
for the symplectic map $f$. It depends on the choice of $\Theta$ and is locally defined.\\

In the following we consider that $P_1=P_2=T^*Q $ with the canonical symplectic form $\omega. $
Since the generating function $S$ is specified on the graph $\Gamma_f$,
and so can be expressed in any local coordinate system on $\Gamma_f$.
The standard choices, for the coordinates $(q,p,\tilde{q},\tilde{p})$ on $T^*Q \times T^*Q, $
are any two of the four quantities $q,\; p,\; \tilde{q}, \; \tilde{p}, $ because $\Gamma_f$
has the same dimension as $T^*Q$. In particular, we choose $(q, \tilde{q})$ as the
local coordinates on $\Gamma_f, $ and consider the generating function $S: Q\times Q \rightarrow \mathbb{R},$
then its differential is given by $\mathbf{d}S= \frac{\partial S}{\partial q}\mathbf{d}q
+ \frac{\partial S}{\partial \tilde{q}}\mathbf{d}\tilde{q}.$
On the other hand, for the canonical symplectic transformation
$f: T^*Q \rightarrow T^*Q, \; (\tilde{q}, \tilde{p})\rightarrow (q,p), $ we have that
$\mathbf{i}^*_f\Theta= \tilde{p}\mathbf{d}\tilde{q} -p\mathbf{d}q, $ and hence the condition
$\mathbf{i}^*_f\Theta= -\mathbf{d}S$ reduced to the following equations
\begin{equation}
p= \frac{\partial S}{\partial q}(q, \tilde{q}), \;\;\;\;\;\;
\tilde{p}= -\frac{\partial S}{\partial \tilde{q}}(q, \tilde{q}).
\end{equation}
Moveover, we consider the Hamiltonian system $(T^*Q, \omega, H),$ the flow $f$ of
Hamiltonian vector field $X_H$ is a symplectic map $f: T^*Q \rightarrow T^*Q. $
From the generating function theory we know that it must have a
generating function $S(q,\tilde{q}).$ Let $(q^i,p_i)=(q^1,\cdots,q^n, p_1,\cdots,p_n)$
denote canonical coordinates with respect to $\omega$ on $T^*Q$, then the Hamilton's equations
in canonical coordinates are
\begin{equation}
\frac{dq^i}{dt}= \frac{\partial H}{\partial p_i}, \;\;\;\;\;\;
\frac{dp_i}{dt}= -\frac{\partial H}{\partial q^i}, \;\;\;\;\;\; i=1, \cdots, n.
\end{equation}
From the conservation of energy we have that
$$
\frac{dH(q,p)}{dt}=\sum^n_{i=1}(\frac{\partial H}{\partial q^i}\cdot \frac{dq^i}{dt}
+ \frac{\partial H}{\partial p_i}\cdot \frac{dp_i}{dt})
=\sum^n_{i=1}(-\frac{dp_i}{dt}\cdot \frac{dq^i}{dt}+\frac{dq^i}{dt}\cdot \frac{dp_i}{dt})=0.
$$
Thus, $H(q,p)=E$ is a constant in $t$. Moreover, consider that $p=
\frac{\partial S}{\partial q}(q, \tilde{q}), $ for the generating
function $S(q,\tilde{q}),$ then we have the classical
Hamilton-Jacobi equation $H(q,\frac{\partial W}{\partial q})=E, $
which is given by Theorem 2.1, where $W=S(q,\tilde{q})$ is function
of $q$ with the parameters $\tilde{q}$, see Abraham and Marsden
\cite{abma78}, Arnold \cite{ar89} and Marsden and Ratiu
\cite{mara99}.\\

For the generating function $S$ of a symplectic map $f: T^*Q \rightarrow T^*Q, $
assume that $\gamma= (\pi_Q)_*(\mathbf{d}S)= \frac{\partial S}{\partial q}\mathbf{d}q$
is an one-form on $Q$,
where $\pi_Q: T^*Q \rightarrow Q$ and $(\pi_Q)_*: T^*T^*Q \rightarrow T^*Q$.
Because
$\mathbf{d}\gamma= \mathbf{d}(\pi_Q)_*(\mathbf{d}S)=(\pi_Q)_*(\mathbf{d}^2S)=0, $
then the one-form $\gamma: Q \rightarrow T^*Q $ is closed with respect to
$T\pi_Q: TT^* Q \rightarrow TQ, $ and from Lemma 2.4(i), hence we have that
for any $v,\; w \in TT^*Q, $
$$\lambda^*\omega(v,w)= -\mathbf{d}\gamma(T\pi_Q(v), T\pi_Q(w))=0. $$
Moreover, from Lemma 2.4(ii), we can obtain that
\begin{equation}
\omega(T\lambda \cdot v, \; w)-\omega(v,\; w)= -\omega(v, \; T\lambda \cdot w).
\end{equation}
If $v$ satisfies the equation $T\lambda\cdot v= v, $ then for any
$w \in TT^*Q, \; T\lambda(w)\neq 0, $ we have that $-\omega(v, \; T\lambda \cdot w)
= -\omega(T\lambda \cdot v, \; T\lambda \cdot w)= -\lambda^*\omega(v,w)=0. $
But, because the symplectic form $\omega$ is non-degenerate,
the left side of (2.4) equals zero, only when
$v$ satisfies the equation $T\lambda\cdot v= v. $ Thus,
for any $v,\; w \in TT^*Q, \; T\lambda(w)\neq 0 $,
we must have that $T\lambda\cdot v= v. $\\

For Hamiltonian system $(T^*Q, \omega, H)$, if we take that
$v= X_H\cdot \gamma \in TT^* Q, $ and for any $w \in TT^* Q, \; T\lambda(w)\neq 0, $
from the equation $T\lambda \cdot v= v, $ we have that
$$T\gamma \cdot X_H^\gamma
= T\gamma \cdot T\pi_Q \cdot X_H \cdot\gamma= T\lambda\cdot X_H \cdot\gamma= X_H \cdot\gamma, $$
that is, $T\gamma \cdot X_H^\gamma= X_H \cdot\gamma, $ this is the Type I of
Hamilton-Jacobi equation for the Hamiltonian system $(T^*Q, \omega, H). $
Moreover, for any symplectic map $\varepsilon: T^*Q \rightarrow T^*Q, $ we take that
$v= X_H\cdot \varepsilon \in TT^* Q, $ and for any $w \in TT^* Q, \; T\lambda(w)\neq 0, $
from the equation $T\lambda \cdot v= v, $ we have that
$T\lambda \cdot X_H \cdot \varepsilon= X_H \cdot\varepsilon. $
Since $\varepsilon: T^*Q \rightarrow T^*Q $ is symplectic, we have that
$ X_H\cdot \varepsilon= T\varepsilon \cdot X_{H\cdot\varepsilon}, $ along $\varepsilon$,
and hence $T\varepsilon \cdot X_{H\cdot\varepsilon}= T\lambda \cdot X_H \cdot \varepsilon. $
On the other hand, note that $\lambda:= \gamma \cdot \pi_Q, $ we can obtain that
$$X_H \cdot\varepsilon= T\lambda \cdot X_H \cdot \varepsilon
= T\gamma \cdot T\pi_Q \cdot X_H \cdot\varepsilon= T\gamma \cdot X_H^\varepsilon .$$
Thus, $\varepsilon$ is a solution of the equation
$T\varepsilon\cdot X_{H\cdot\varepsilon}= T\lambda \cdot X_H \cdot\varepsilon,$ if and only if it is a solution of
the Type II of Hamilton-Jacobi equation $T\gamma\cdot X_H^\varepsilon= X_H\cdot
\varepsilon ,$ for the Hamiltonian system $(T^*Q, \omega, H). $\\

To sum up the above discussion, if the one-form $\gamma= (\pi_Q)_*(\mathbf{d}S)= \frac{\partial S}{\partial q}\mathbf{d}q$ is
given by a generating function $S$ of a symplectic map, then the
classical Hamilton-Jacobi equation $H(q, \gamma(q))=E,$(constant in $t$),
or equivalently, $\mathbf{d}(H\cdot \gamma)=0, $ as well as the Type I of
Hamilton-Jacobi equation $T\gamma \cdot X_H^\gamma= X_H \cdot\gamma, $ and the Type II of
Hamilton-Jacobi theorem, all of them hold.

\begin{remark}
It is worthy of note that, we can obtain the Type I and Type II of
Hamilton-Jacobi equation from Theorem $2.5$ and Theorem $2.6$, even
if the one-form $\gamma: Q \rightarrow T^*Q$ may not be given by a
generating function of a symplectic map. Thus, the formulations of
Type I and Type II of Hamilton-Jacobi equation have more extensive
sense. On the other hand, if $\gamma$ is a solution of the classical
Hamilton-Jacobi equation, that is, $\mathbf{d}(H\cdot \gamma)=0, $
or equivalently, $X_H\cdot \gamma=0, $ which shows that the
dynamical vector field of the Hamiltonian system $(T^*Q,\omega,H)$
is degenerate along $\gamma$, in this case, $X_H^\gamma= T\pi_Q\cdot
X_H\cdot \gamma=0,$ and hence the Type I of Hamilton-Jacobi
equation, $X_{H}\cdot \gamma= T\gamma\cdot X_H^\gamma, $ holds
trivially. In addition, for a symplectic map $\varepsilon: T^* Q
\rightarrow T^* Q $, if $X_H\cdot \varepsilon=0, $ then from the
Type II of Hamilton-Jacobi equation, we have that $X_{H}\cdot
\varepsilon= T\gamma\cdot X_H^\varepsilon=0. $ Moreover, from the
equation $T\varepsilon\cdot X_{H\cdot\varepsilon}= X_H
\cdot\varepsilon,$ we know that $X_{H}\cdot \varepsilon=0 $ is
equivalent to $X_{H\cdot\varepsilon}=0.$
\end{remark}

In the following we shall generalize the above Type I and Type II of
Hamilton-Jacobi theorem to the regular point and the regular orbit
reducible Hamiltonian systems with symmetries, and give a variety of
Hamilton-Jacobi theorems for the regular reduced Hamiltonian
systems.

\section{Hamilton-Jacobi Theorem of Regular Point Reduced Hamiltonian
System}

In this section, we first give the regular point reducible
Hamiltonian system with symmetry. Then we prove the Type I and Type II of
Hamilton-Jacobi theorems for the regular point reduced Hamiltonian system,
by using Lemma 2.4, the regular point reduced symplectic form and the reduced
dynamical vector field.\\

At first, we consider the regular point reducible Hamiltonian system.
Let $Q$ be a smooth manifold and $T^\ast Q$ its cotangent bundle
with the symplectic form $\omega$. Let $\Phi: G\times Q\rightarrow
Q$ be a smooth left action of a Lie group $G$ on $Q$, which is free
and proper. Then the cotangent lifted left action $\Phi^{T^\ast}:
G\times T^\ast Q\rightarrow T^\ast Q$ is symplectic, free and
proper. Assume that the action admits an
$\operatorname{Ad}^\ast$-equivariant momentum map $\mathbf{J}:T^\ast
Q\rightarrow \mathfrak{g}^\ast$, where $\mathfrak{g}$ is the Lie
algebra of $G$ and $\mathfrak{g}^\ast$ is the dual of
$\mathfrak{g}$. Let $\mu\in\mathfrak{g}^\ast$ be a regular value of
$\mathbf{J}$ and denote by $G_\mu$ the isotropy subgroup of the
coadjoint $G$-action at the point $\mu\in\mathfrak{g}^\ast$, which
is defined by $G_\mu=\{g\in G|\operatorname{Ad}_g^\ast \mu=\mu \}$.
Since $G_\mu (\subset G)$ acts freely and properly on $Q$ and on
$T^\ast Q$, then $Q_\mu=Q/G_\mu$ is a smooth manifold and that the
canonical projection $\rho_\mu:Q\rightarrow Q_\mu$ is a surjective
submersion. It follows that $G_\mu$ acts also freely and properly on
$\mathbf{J}^{-1}(\mu)$, so that the space $(T^\ast
Q)_\mu=\mathbf{J}^{-1}(\mu)/G_\mu$ is a symplectic manifold with the
symplectic form $\omega_\mu$ uniquely characterized by the relation
\begin{equation}\pi_\mu^\ast \omega_\mu=i_\mu^\ast
\omega. \label{3.1}\end{equation} The map
$i_\mu:\mathbf{J}^{-1}(\mu)\rightarrow T^\ast Q$ is the inclusion
and $\pi_\mu:\mathbf{J}^{-1}(\mu)\rightarrow (T^\ast Q)_\mu$ is the
projection. The pair $((T^\ast Q)_\mu,\omega_\mu)$ is called
Marsden-Weinstein reduced space of $(T^\ast Q,\omega)$ at $\mu$.\\

\begin{remark}
If $(T^\ast Q, \omega)$ is a connected symplectic manifold, and
$\mathbf{J}:T^\ast Q\rightarrow \mathfrak{g}^\ast$ is a
non-equivariant momentum map with a non-equivariance group
one-cocycle $\sigma: G\rightarrow \mathfrak{g}^\ast$, which is
defined by $\sigma(g):=\mathbf{J}(g\cdot
z)-\operatorname{Ad}^\ast_{g^{-1}}\mathbf{J}(z)$, where $g\in G$ and
$z\in T^\ast Q$. Then we know that $\sigma$ produces a new affine
action $\Theta: G\times \mathfrak{g}^\ast \rightarrow
\mathfrak{g}^\ast $ defined by
$\Theta(g,\mu):=\operatorname{Ad}^\ast_{g^{-1}}\mu + \sigma(g)$,
where $\mu \in \mathfrak{g}^\ast$, with respect to which the given
momentum map $\mathbf{J}$ is equivariant. Assume that $G$ acts
freely and properly on $T^\ast Q$, and $\tilde{G}_\mu$ denotes the
isotropy subgroup of $\mu \in \mathfrak{g}^\ast$ relative to this
affine action $\Theta$ and $\mu$ is a regular value of $\mathbf{J}$.
Then the quotient space $(T^\ast
Q)_\mu=\mathbf{J}^{-1}(\mu)/\tilde{G}_\mu$ is also a symplectic
manifold with the symplectic form $\omega_\mu$ uniquely characterized by
$(3.1)$, see Ortega and Ratiu \cite{orra04}.
\end{remark}

Let $H: T^\ast Q\rightarrow \mathbb{R}$ be a $G$-invariant
Hamiltonian, the flow $F_t$ of the Hamiltonian vector field $X_H$
leaves the connected components of $\mathbf{J}^{-1}(\mu)$ invariant
and commutes with the $G$-action, so it induces a flow $f_t^\mu$ on
$(T^\ast Q)_\mu$, defined by $f_t^\mu\cdot \pi_\mu=\pi_\mu \cdot
F_t\cdot i_\mu$, and the vector field $X_{h_\mu}$ generated by the
flow $f_t^\mu$ on $((T^\ast Q)_\mu,\omega_\mu)$ is Hamiltonian with
the associated regular point reduced Hamiltonian function
$h_\mu:(T^\ast Q)_\mu\rightarrow \mathbb{R}$ defined by
$h_\mu\cdot\pi_\mu=H\cdot i_\mu$, and the Hamiltonian vector fields
$X_H$ and $X_{h_\mu}$ are $\pi_\mu$-related. Thus, we can define a
regular point reducible Hamiltonian system as follows.

\begin{definition}
(Regular Point Reducible Hamiltonian System) A 4-tuple $(T^\ast Q,
G,\omega,H )$, where the Hamiltonian $H: T^\ast Q\rightarrow
\mathbb{R}$ is $G$-invariant, is called a regular point reducible
Hamiltonian system, if there exists a point
$\mu\in\mathfrak{g}^\ast$, which is a regular value of the momentum
map $\mathbf{J}$, such that the regular point reduced system, that
is, the 3-tuple $((T^\ast Q)_\mu, \omega_\mu,h_\mu )$, where
$(T^\ast Q)_\mu=\mathbf{J}^{-1}(\mu)/G_\mu$, $\pi_\mu^\ast
\omega_\mu=i_\mu^\ast\omega$, $h_\mu\cdot \pi_\mu=H\cdot i_\mu$, is
a Hamiltonian system, which is also called Marsden-Weinstein reduced
Hamiltonian system. Here $((T^\ast Q)_\mu,\omega_\mu)$ is
Marsden-Weinstein reduced space, the function $h_\mu:(T^\ast
Q)_\mu\rightarrow \mathbb{R}$ is called the reduced Hamiltonian.
\end{definition}

For the regular point reducible Hamiltonian system
$(T^*Q,G,\omega,H)$, by using Lemma 2.4, the regular point reduced symplectic form and the reduced
dynamical vector field, we can prove the following two types of Hamilton-Jacobi
theorem for the regular point reduced Hamiltonian system $((T^\ast
Q)_\mu, \omega_\mu,h_\mu )$.
At first, by using the fact that the one-form $\gamma: Q
\rightarrow T^*Q $ is closed with respect to
$T\pi_Q: TT^* Q \rightarrow TQ, $ and $\textmd{Im}(\gamma)\subset
\mathbf{J}^{-1}(\mu), $ and it is $G_\mu$-invariant, we can prove the following Type I of
Hamilton-Jacobi theorem for the regular point reduced Hamiltonian system $((T^\ast
Q)_\mu, \omega_\mu,h_\mu )$.
For convenience, the maps involved in
the following theorem and its proof are shown in Diagram-3.

\begin{center}
\hskip 0cm \xymatrix{ \mathbf{J}^{-1}(\mu) \ar[r]^{i_\mu} & T^* Q \ar[r]^{\pi_Q}
& Q \ar[d]_{X_H^\gamma} \ar[r]^{\gamma}
& T^*Q \ar[d]_{X_H} \ar[r]^{\pi_\mu}
& (T^* Q)_\mu \ar[d]_{X_{h_\mu}} \\
& T(T^*Q)  & TQ \ar[l]^{T\gamma}
& T(T^*Q) \ar[l]^{T\pi_Q} \ar[r]_{T\pi_\mu} & T(T^* Q)_\mu }
\end{center}
$$\mbox{Diagram-3}$$

\begin{theorem} (Type I of Hamilton-Jacobi Theorem for a Regular Point Reduced Hamiltonian System)
For the regular point reducible Hamiltonian system
$(T^*Q,G,\omega,H)$, assume that $\gamma: Q \rightarrow T^*Q$ is an one-form
on $Q$, and $X_H^\gamma = T\pi_{Q}\cdot X_H \cdot \gamma$,
where $X_{H}$ is the dynamical vector field of $(T^*Q,G,\omega,H)$. Moreover,
assume that $\mu \in \mathfrak{g}^\ast $ is a regular value of the momentum
map $\mathbf{J}$, and $\textmd{Im}(\gamma)\subset
\mathbf{J}^{-1}(\mu), $ and it is $G_\mu$-invariant, and
$\bar{\gamma}=\pi_\mu(\gamma): Q \rightarrow (T^* Q)_\mu. $
If the one-form $\gamma: Q \rightarrow T^*Q $ is closed with respect to
$T\pi_Q: TT^* Q \rightarrow TQ, $
then $\bar{\gamma}$ is a solution of the equation
$T\bar{\gamma}\cdot X_H^\gamma= X_{h_\mu}\cdot \bar{\gamma}, $
which is called the Type I of Hamilton-Jacobi equation for the Marsden-Weinstein
reduced Hamiltonian system $((T^\ast Q)_\mu, \omega_\mu,h_\mu)$.
\end{theorem}

\noindent{\bf Proof: } At first, from Theorem 2.5, we know that
$\gamma$ is a solution of the Type I of Hamilton-Jacobi equation
$T\gamma\cdot X_H^\gamma= X_H\cdot \gamma. $ Next, we note that
$\textmd{Im}(\gamma)\subset \mathbf{J}^{-1}(\mu), $ and it
is $G_\mu$-invariant, in this case $\pi_\mu^*\omega_\mu=
i_\mu^*\omega= \omega, $ along $\textmd{Im}(\gamma)$.
By using the reduced symplectic form $\omega_\mu$, if we take
that $v= X_H\cdot \gamma \in TT^* Q,$ and for any $w \in TT^* Q, \; T\pi_{Q}(w)\neq 0,$
and $T\pi_{\mu} (w) \neq 0, $ from Lemma 2.4(ii) we have that
\begin{align*}
& \omega_\mu(T\bar{\gamma} \cdot X_H^\gamma, \; T\pi_\mu \cdot w) =
\omega_\mu(T(\pi_\mu \cdot \gamma) \cdot X_H^\gamma, \; T\pi_\mu
\cdot w )\\
& = \pi_\mu^*\omega_\mu(T\gamma \cdot X_H^\gamma, \; w)
= \omega(T(\gamma \cdot \pi_Q)\cdot X_H\cdot \gamma, \; w)\\
& = \omega(X_H\cdot \gamma, \; w-T(\gamma \cdot \pi_Q)\cdot w)
-\mathbf{d}\gamma(T\pi_{Q}(X_H\cdot \gamma), \; T\pi_{Q}(w))\\
& =\omega(X_H\cdot \gamma, \; w) - \omega(X_H\cdot \gamma, \;
T(\gamma \cdot \pi_Q)\cdot w)-\mathbf{d}\gamma(T\pi_{Q}(X_H\cdot \gamma), \; T\pi_{Q}(w))\\
& =\pi_\mu^*\omega_\mu(X_H\cdot
\gamma, \; w) - \pi_\mu^*\omega_\mu(X_H\cdot \gamma, \; T(\gamma \cdot \pi_Q)\cdot w)
-\mathbf{d}\gamma(T\pi_{Q}(X_H\cdot \gamma), \; T\pi_{Q}(w))\\
& = \omega_\mu(T\pi_\mu(X_H\cdot \gamma), \;
T\pi_\mu \cdot w) - \omega_\mu(T\pi_\mu\cdot(X_H\cdot \gamma), \;
T(\pi_\mu \cdot\gamma \cdot \pi_Q) \cdot w) \\
& \;\;\;\; -\mathbf{d}\gamma(T\pi_{Q}(X_H\cdot \gamma), \; T\pi_{Q}(w))\\
& = \omega_\mu(T\pi_\mu(X_H)\cdot \pi_\mu(\gamma), \; T\pi_\mu \cdot w)
- \omega_\mu(T\pi_\mu(X_H)\cdot \pi_\mu(\gamma), \; T\bar{\gamma}\cdot T\pi_{Q}(w))\\
& \;\;\;\; -\mathbf{d}\gamma(T\pi_{Q}(X_H\cdot \gamma), \; T\pi_{Q}(w))\\
& = \omega_\mu(X_{h_\mu} \cdot
\bar{\gamma}, \; T\pi_\mu \cdot w)- \omega_\mu(X_{h_\mu} \cdot
\bar{\gamma}, \; T\bar{\gamma} \cdot T\pi_{Q}(w))-\mathbf{d}\gamma(T\pi_{Q}(X_H\cdot \gamma), \; T\pi_{Q}(w)),
\end{align*}
where we have used that $T\pi_\mu(X_H)= X_{h_\mu}. $
Since the one-form $\gamma: Q \rightarrow T^*Q $ is closed with respect to
$T\pi_Q: TT^* Q \rightarrow TQ, $ then we have that $
\mathbf{d}\gamma(T\pi_{Q}(X_H\cdot \gamma), \; T\pi_{Q}(w))=0, $
and hence
\begin{equation}
 \omega_\mu(T\bar{\gamma} \cdot X_H^\gamma, \; T\pi_\mu \cdot w)-
\omega_\mu(X_{h_\mu} \cdot \bar{\gamma}, \; T\pi_\mu \cdot w)
 = - \omega_\mu(X_{h_\mu} \cdot
\bar{\gamma}, \; T\bar{\gamma} \cdot T\pi_{Q}(w)).
\end{equation}
If $\bar{\gamma}$ satisfies the equation $T\bar{\gamma}\cdot X_H^\gamma= X_{h_\mu}\cdot \bar{\gamma}, $
from Lemma 2.4(i) we can obtain that
\begin{align*}
- \omega_\mu(X_{h_\mu} \cdot
\bar{\gamma}, \; T\bar{\gamma} \cdot T\pi_{Q}(w))
& = -\omega_\mu (T\bar{\gamma} \cdot X_H^\gamma, \; T\bar{\gamma} \cdot T\pi_{Q}(w))\\
& = -\bar{\gamma}^*\omega_\mu (T\pi_{Q} \cdot X_{H}\cdot\gamma, \; T\pi_Q(w))\\
& = -\gamma^* \cdot \pi_\mu^*\omega_\mu (T\pi_{Q} \cdot X_{H}\cdot\gamma, \; T\pi_{Q}(w))\\
& = -\gamma^*\omega( T\pi_{Q}(X_{H}\cdot\gamma), \; T\pi_{Q}(w))\\
& = \textbf{d}\gamma(T\pi_{Q}( X_{H}\cdot\gamma ), \; T\pi_{Q}(w))=0.
\end{align*}
Because the reduced symplectic form $\omega_\mu$ is non-degenerate, the left side of (3.2) equals zero, only when
$\bar{\gamma}$ satisfies the equation $T\bar{\gamma}\cdot X_H^\gamma= X_{h_\mu}\cdot \bar{\gamma}.$ Thus,
if the one-form $\gamma: Q \rightarrow T^*Q $ is closed with respect to
$T\pi_Q: TT^* Q \rightarrow TQ, $ then $\bar{\gamma}$ must be a solution of the Type I of Hamilton-Jacobi equation
$T\bar{\gamma}\cdot X_H^\gamma= X_{h_\mu}\cdot \bar{\gamma}.$
\hskip 0.3cm $\blacksquare$\\

Next, for any $G_\mu$-invariant symplectic map $\varepsilon: T^* Q \rightarrow T^* Q $,
we can prove the following Type II of
Hamilton-Jacobi theorem for the regular point reduced Hamiltonian system $((T^\ast
Q)_\mu, \omega_\mu,h_\mu )$.
For convenience, the maps involved in
the following theorem and its proof are shown in Diagram-4.

\begin{center}
\hskip 0cm \xymatrix{ \mathbf{J}^{-1}(\mu) \ar[r]^{i_\mu} & T^* Q
\ar[d]_{X_{H\cdot \varepsilon}} \ar[dr]^{X_H^\varepsilon} \ar[r]^{\pi_Q}
& Q \ar[r]^{\gamma} & T^*Q \ar[d]_{X_H} \ar[dr]^{X_{h_\mu} \cdot\bar{\varepsilon}} \ar[r]^{\pi_\mu} & (T^* Q)_\mu \ar[d]^{X_{h_\mu}} \\
  & T(T^*Q)  & TQ \ar[l]^{T\gamma} & T(T^*Q) \ar[l]^{T\pi_Q} \ar[r]_{T\pi_\mu} & T(T^* Q)_\mu }
\end{center}
$$\mbox{Diagram-4}$$

\begin{theorem} (Type II of Hamilton-Jacobi Theorem for a Regular Point Reduced Hamiltonian System)
For the regular point reducible Hamiltonian system
$(T^*Q,G,\omega,H)$, assume that $\gamma: Q \rightarrow T^*Q$ is an one-form
on $Q$, and $\lambda=\gamma \cdot \pi_{Q}: T^* Q \rightarrow T^* Q, $
and for any symplectic map $\varepsilon:T^* Q \rightarrow T^* Q, $
denote by $X_H^\varepsilon = T\pi_{Q}\cdot X_H \cdot \varepsilon$,
where $X_{H}$ is the dynamical vector field of $(T^*Q,G,\omega,H)$. Moreover,
assume that $\mu \in \mathfrak{g}^\ast $ is a regular value of the momentum
map $\mathbf{J}$, and $\textmd{Im}(\gamma)\subset
\mathbf{J}^{-1}(\mu), $ and it is $G_\mu$-invariant,
and $\varepsilon$ is $G_\mu$-invariant and $\varepsilon(\mathbf{J}^{-1}(\mu)) \subset \mathbf{J}^{-1}(\mu). $
Denote by $\bar{\gamma}=\pi_\mu(\gamma): Q \rightarrow (T^* Q)_\mu $,
$\bar{\lambda}=\pi_\mu(\lambda): \mathbf{J}^{-1}(\mu) (\subset T^*Q) \rightarrow (T^* Q)_\mu $,
and $\bar{\varepsilon}=\pi_\mu(\varepsilon): \mathbf{J}^{-1}(\mu) (\subset T^*Q) \rightarrow (T^* Q)_\mu $.
Then $\varepsilon$ and $\bar{\varepsilon}$ satisfy
the equation $T\bar{\varepsilon}\cdot(X_{h_\mu \cdot \bar{\varepsilon}})= T\bar{\lambda}\cdot X_H \cdot \varepsilon, $
if and only if they satisfy the equation
$T\bar{\gamma}\cdot X_H^\varepsilon= X_{h_\mu}\cdot \bar{\varepsilon}, $
where $X_{h_\mu \cdot \bar{\varepsilon}} \in TT^* Q$ is the Hamiltonian
vector field of the function $h_\mu \cdot \bar{\varepsilon}: T^* Q
\rightarrow \mathbb{R}. $
The equation $T\bar{\gamma}\cdot X_H^\varepsilon=
X_{h_\mu}\cdot \bar{\varepsilon}$ is called the Type II of Hamilton-Jacobi equation for
the Marsden-Weinstein
reduced Hamiltonian system $((T^\ast Q)_\mu, \omega_\mu,h_\mu)$.
\end{theorem}

\noindent{\bf Proof: } At first, we note that
$\textmd{Im}(\gamma)\subset \mathbf{J}^{-1}(\mu), $ and it
is $G_\mu$-invariant, in this case, $\pi_\mu^*\omega_\mu=
i_\mu^*\omega= \omega, $ along $\textmd{Im}(\gamma)$.
By using the reduced symplectic form $\omega_\mu$, if we take
that $v= X_H\cdot \varepsilon \in TT^* Q,$ and for any $w \in TT^* Q, \; T\bar{\lambda}(w)\neq 0,$
and $T\pi_{\mu} (w) \neq 0, $ from Lemma 2.4 we have that
\begin{align*}
& \omega_\mu(T\bar{\gamma} \cdot X_H^\varepsilon, \; T\pi_\mu \cdot w) =
\omega_\mu(T(\pi_\mu \cdot \gamma) \cdot X_H^\varepsilon, \; T\pi_\mu
\cdot w )\\
& = \pi_\mu^*\omega_\mu(T\gamma \cdot X_H^\varepsilon, \; w)
= \omega(T(\gamma \cdot \pi_Q)\cdot X_H\cdot \varepsilon, \; w)\\
& = \omega(X_H\cdot \varepsilon, \; w-T(\gamma \cdot \pi_Q)\cdot w)
-\mathbf{d}\gamma(T\pi_{Q}(X_H\cdot \varepsilon), \; T\pi_{Q}(w))\\
& =\omega(X_H\cdot \varepsilon, \; w) - \omega(X_H\cdot \varepsilon, \;
T\lambda\cdot w)+ \lambda^*\omega(X_H\cdot \varepsilon, \; w)\\
& =\pi_\mu^*\omega_\mu(X_H\cdot
\varepsilon, \; w) - \pi_\mu^*\omega_\mu(X_H\cdot \varepsilon, \; T\lambda\cdot w)+ \lambda^*\cdot \pi_\mu^*\omega_\mu(X_H\cdot \varepsilon, \; w)\\
& = \omega_\mu(T\pi_\mu(X_H\cdot \varepsilon), \;
T\pi_\mu \cdot w) - \omega_\mu(T\pi_\mu\cdot(X_H\cdot \varepsilon), \;
T(\pi_\mu \cdot\lambda) \cdot w)+ (\pi_\mu\cdot\lambda)^*\cdot \omega_\mu(X_H\cdot \varepsilon, \; w)\\
& = \omega_\mu(T\pi_\mu(X_H)\cdot \pi_\mu(\varepsilon), \; T\pi_\mu \cdot w)
- \omega_\mu(T\pi_\mu(X_H)\cdot \pi_\mu(\varepsilon), \; T\bar{\lambda}\cdot w)
+ \omega_\mu(T\bar{\lambda}\cdot X_H\cdot \varepsilon, \; T\bar{\lambda}\cdot w)\\
& = \omega_\mu(X_{h_\mu} \cdot
\bar{\varepsilon}, \; T\pi_\mu \cdot w)- \omega_\mu(X_{h_\mu} \cdot
\bar{\varepsilon}, \; T\bar{\lambda} \cdot w)+ \omega_\mu(T\bar{\lambda}\cdot X_H\cdot \varepsilon, \; T\bar{\lambda}\cdot w),
\end{align*}
where we have used that $T\pi_\mu(X_H)= X_{h_\mu}. $
Note that $\varepsilon: T^* Q \rightarrow T^* Q $ is
symplectic, and $\pi_\mu^*\omega_\mu=
i_\mu^*\omega= \omega, $ along $\textmd{Im}(\gamma)$, and hence
$\bar{\varepsilon}= \pi_\mu(\varepsilon): \mathbf{J}^{-1}(\mu)(\subset T^* Q) \rightarrow (T^*
Q)_\mu $ is also symplectic along
$\textmd{Im}(\gamma)$, and $X_{h_\mu}\cdot \bar{\varepsilon}=
T\bar{\varepsilon} \cdot X_{h_\mu \cdot \bar{\varepsilon}}, $ along $\textmd{Im}(\gamma)\cap\textmd{Im}(\varepsilon)$.
From the above arguments, we can obtain that
\begin{align*}
& \omega_\mu(T\bar{\gamma} \cdot X_H^\varepsilon, \; T\pi_\mu \cdot w)-
\omega_\mu(X_{h_\mu} \cdot \bar{\varepsilon}, \; T\pi_\mu \cdot w) \\
& = \omega_\mu(T\bar{\lambda}\cdot X_H\cdot \varepsilon, \; T\bar{\lambda}\cdot w)
-\omega_\mu(T\bar{\varepsilon} \cdot X_{h_\mu \cdot \bar{\varepsilon}}, \; T\bar{\lambda} \cdot w)\\
& = \omega_\mu(T\bar{\lambda}\cdot X_H\cdot \varepsilon- T\bar{\varepsilon} \cdot X_{h_\mu \cdot \bar{\varepsilon}}, \; T\bar{\lambda}\cdot w).
\end{align*}
Because the reduced symplectic form $\omega_\mu$ is non-degenerate,
it follows that
$T\bar{\gamma}\cdot X_H^\varepsilon=
X_{h_\mu}\cdot \bar{\varepsilon}, $ is equivalent to
$T\bar{\lambda}\cdot X_H\cdot \varepsilon= T\bar{\varepsilon} \cdot X_{h_\mu \cdot \bar{\varepsilon}}. $
Thus, we know that the $\varepsilon$ and $\bar{\varepsilon}$ satisfy
the equation $T\bar{\varepsilon}\cdot(X_{h_\mu \cdot \bar{\varepsilon}})= T\bar{\lambda}\cdot X_H \cdot \varepsilon, $
if and only if they satisfy the Type II of Hamilton-Jacobi equation
$T\bar{\gamma}\cdot X_H^\varepsilon= X_{h_\mu}\cdot \bar{\varepsilon}. $
\hskip 0.3cm $\blacksquare$\\

Moreover, for the regular point reducible Hamiltonian system
$(T^*Q,G,\omega,H)$, we know that the dynamical vector fields
$X_{H}$ and $X_{h_\mu}$ are $\pi_\mu$-related, that is,
$X_{h_\mu}\cdot \pi_\mu=T\pi_\mu\cdot X_{H}\cdot i_\mu.$ Then we can
prove the following Theorem 3.5, which states the relationship
between the solution of Type II of Hamilton-Jacobi equation and the
Marsden-Weinstein reduction.

\begin{theorem}
For the regular point reducible Hamiltonian system
$(T^*Q,G,\omega,H)$, assume that $\gamma: Q \rightarrow T^*Q$ is an one-form on $Q$,
and $\lambda=\gamma \cdot \pi_{Q}: T^* Q
\rightarrow T^* Q $, and $\varepsilon: T^* Q \rightarrow T^* Q $ is a symplectic map.
Moreover, assume that $\mu \in
\mathfrak{g}^\ast $ is a regular value of the momentum
map $\mathbf{J}$, and $\textmd{Im}(\gamma)\subset
\mathbf{J}^{-1}(\mu), $ and it is $G_\mu$-invariant, and $\varepsilon$ is $G_\mu$-invariant
and $\varepsilon(\mathbf{J}^{-1}(\mu)) \subset \mathbf{J}^{-1}(\mu). $
Denote by $\bar{\gamma}=\pi_\mu(\gamma): Q \rightarrow (T^* Q)_\mu $,
$\bar{\lambda}=\pi_\mu(\lambda): \mathbf{J}^{-1}(\mu) (\subset T^*Q) \rightarrow (T^* Q)_\mu $,
and $\bar{\varepsilon}=\pi_\mu(\varepsilon): \mathbf{J}^{-1}(\mu) (\subset T^*Q) \rightarrow (T^* Q)_\mu $.
Then $\varepsilon$ is a solution of the Type II of Hamilton-Jacobi equation
$T\gamma\cdot X_H^\varepsilon= X_H\cdot \varepsilon, $ for the
regular point reducible Hamiltonian system $(T^*Q,G,\omega,H), $ if and only if
$\varepsilon$ and $\bar{\varepsilon} $ satisfy the Type II of Hamilton-Jacobi equation
$T\bar{\gamma}\cdot X_H^\varepsilon= X_{h_\mu}\cdot \bar{\varepsilon}, $ for the
Marsden-Weinstein reduced Hamiltonian system $((T^\ast Q)_\mu,
\omega_\mu,h_\mu)$.
\end{theorem}

\noindent{\bf Proof: }
Note that $\textmd{Im}(\gamma)\subset \mathbf{J}^{-1}(\mu), $ and it
is $G_\mu$-invariant, in this case, $\pi_\mu^*\omega_\mu=
i_\mu^*\omega= \omega, $ along $\textmd{Im}(\gamma)$.
Since the dynamical vector fields
$X_{H}$ and $X_{h_\mu}$ are $\pi_\mu$-related, that is,
$X_{h_\mu}\cdot \pi_\mu= T\pi_\mu\cdot X_{H}\cdot i_\mu, $ and
by using the reduced symplectic form $\omega_\mu$, we have that
\begin{align*}
& \omega_\mu(T\bar{\gamma} \cdot X_H^\varepsilon
- X_{h_\mu} \cdot \bar{\varepsilon}, \; T\pi_\mu \cdot w) \\
& = \omega_\mu(T\bar{\gamma} \cdot X_H^\varepsilon, \; T\pi_\mu \cdot w)-
\omega_\mu(X_{h_\mu} \cdot \bar{\varepsilon}, \; T\pi_\mu \cdot w) \\
& = \omega_\mu(T\pi_\mu \cdot T\gamma \cdot X_H^\varepsilon, \; T\pi_\mu \cdot w)-
\omega_\mu(X_{h_\mu} \cdot \pi_\mu \cdot \varepsilon, \; T\pi_\mu \cdot w) \\
& = \pi_\mu^*\omega_\mu(T\gamma \cdot X_H^\varepsilon, \; w)
-\omega_\mu(T\pi_\mu\cdot X_{H}\cdot \varepsilon, \; T\pi_\mu \cdot w) \\
& = \pi_\mu^*\omega_\mu(T\gamma \cdot X_H^\varepsilon, \; w)
-\pi_\mu^*\omega_\mu(X_{H}\cdot \varepsilon, \; w)\\
& = \omega(T\gamma \cdot X_H^\varepsilon- X_{H}\cdot \varepsilon, \; w).
\end{align*}
Because the symplectic form $\omega$ and the reduced symplectic form $\omega_\mu$ are non-degenerate,
it follows that the equation
$T\bar{\gamma}\cdot X_H^\varepsilon= X_{h_\mu}\cdot \bar{\varepsilon}, $
is equivalent to the equation $T\gamma\cdot X_H^\varepsilon= X_H\cdot \varepsilon$. Thus,
$\varepsilon$ is a solution of the Type II of Hamilton-Jacobi equation
$T\gamma\cdot X_H^\varepsilon= X_H\cdot \varepsilon, $ for the
regular point reducible Hamiltonian system $(T^*Q,G,\omega,H), $ if and only if
$\varepsilon$ and $\bar{\varepsilon} $ satisfy the Type II of Hamilton-Jacobi equation
$T\bar{\gamma}\cdot X_H^\varepsilon= X_{h_\mu}\cdot \bar{\varepsilon}, $ for the
Marsden-Weinstein reduced Hamiltonian system $((T^\ast Q)_\mu,
\omega_\mu,h_\mu)$.  \hskip 0.3cm
$\blacksquare$\\

\begin{remark}
If $(T^\ast Q, \omega)$ is a connected symplectic manifold, and
$\mathbf{J}:T^\ast Q\rightarrow \mathfrak{g}^\ast$ is a
non-equivariant momentum map with a non-equivariance group
one-cocycle $\sigma: G\rightarrow \mathfrak{g}^\ast$, in this case,
for the given regular point reducible Hamiltonian system
$(T^*Q,G,\omega,H)$, we can also prove the Type I and Type II of
Hamilton-Jacobi theorem for the regular point reduced Hamiltonian
system $((T^\ast Q)_\mu, \omega_\mu,h_\mu )$ by using the above
similar ways, where the reduced space $((T^\ast Q)_\mu, \omega_\mu
)$ is determined by the affine action given in Remark 3.1.
\end{remark}

\begin{remark}
It is worthy of note that, the one-form $\gamma: Q \rightarrow T^*Q$
may not be given by a generating function of a symplectic map, and
hence the formulations of Type I and Type II of Hamilton-Jacobi
equation for a regular point reducible Hamiltonian system, given by
Theorem $3.3$ and Theorem $3.4$, have more extensive sense. On the
other hand, if $\gamma$ is a solution of the classical
Hamilton-Jacobi equation, that is, $X_H\cdot \gamma=0, $ then
$X_H^\gamma= T\pi_Q\cdot X_H\cdot \gamma=0,$ and hence from the Type
I of Hamilton-Jacobi equation, we have that $X_{h_\mu}\cdot
\bar{\gamma}= T\bar{\gamma}\cdot X_H^\gamma=0. $ Since the classical
Hamilton-Jacobi equation $X_H\cdot \gamma=0, $ shows that the
dynamical vector field of the regular point reducible Hamiltonian
system $(T^*Q,G,\omega,H)$ is degenerate along $\gamma$, then the
equation $X_{h_\mu}\cdot \bar{\gamma}=0,$ shows that the dynamical
vector field of the Marsden-Weinstein reduced Hamiltonian system
$((T^\ast Q)_\mu, \omega_\mu,h_\mu)$ is degenerate along
$\bar{\gamma}$. The equation $X_{h_\mu}\cdot \bar{\gamma}=0$ is
called the classical Hamilton-Jacobi equation for the
Marsden-Weinstein reduced Hamiltonian system $((T^\ast Q)_\mu,
\omega_\mu,h_\mu). $ In addition, for a symplectic map $\varepsilon:
T^* Q \rightarrow T^* Q $, if $X_H\cdot \varepsilon=0, $ then from
the Type II of Hamilton-Jacobi equation, we have that
$X_{h_\mu}\cdot \bar{\varepsilon}= T\bar{\gamma}\cdot
X_H^\varepsilon=0. $ Moreover, from the equation
$T\bar{\varepsilon}\cdot(X_{h_\mu \cdot \bar{\varepsilon}})=
T\bar{\lambda}\cdot X_H \cdot \varepsilon, $ we know that
$X_{h_\mu}\cdot \bar{\varepsilon}=0 $ is equivalent to $X_{h_\mu
\cdot \bar{\varepsilon}}=0.$
\end{remark}

\section{Hamilton-Jacobi Theorem of Regular Orbit Reduced Hamiltonian System }

The orbit reduction is an alternative approach to symplectic reduction
given by Marle \cite{ma76}
and Kazhdan, Kostant and Sternberg \cite{kakost78}, which is different from the
Marsden-Weinstein reduction.
In this section, we first give the regular orbit reducible
Hamiltonian system with symmetry. Then we prove the Type I and Type II of
Hamilton-Jacobi theorems for the regular orbit reduced Hamiltonian system,
by using Lemma 2.4, the orbit reduced symplectic form and the reduced
dynamical vector field.\\

At first, we consider the regular orbit reducible Hamiltonian system.
For the cotangent lifted
left action $\Phi^{T^\ast}:G\times T^\ast Q\rightarrow T^\ast Q$, which is
symplectic, free and proper, assume that the action admits an
$\operatorname{Ad}^\ast$-equivariant momentum map $\mathbf{J}:T^\ast
Q\rightarrow \mathfrak{g}^\ast$. Let $\mu\in \mathfrak{g}^\ast$ be a
regular value of the momentum map $\mathbf{J}$ and
$\mathcal{O}_\mu=G\cdot \mu\subset \mathfrak{g}^\ast$ be the
$G$-orbit of the coadjoint $G$-action through the point $\mu$. Since
$G$ acts freely, properly and symplectically on $T^\ast Q$, then the
quotient space $(T^\ast Q)_{\mathcal{O}_\mu}=
\mathbf{J}^{-1}(\mathcal{O}_\mu)/G$ is a regular quotient symplectic
manifold with the symplectic form $\omega_{\mathcal{O}_\mu}$
uniquely characterized by the relation
\begin{equation}i_{\mathcal{O}_\mu}^\ast \omega=\pi_{\mathcal{O}_{\mu}}^\ast
\omega_{\mathcal{O}
_\mu}+\mathbf{J}_{\mathcal{O}_\mu}^\ast\omega_{\mathcal{O}_\mu}^+,
\label{4.1}\end{equation} where $\mathbf{J}_{\mathcal{O}_\mu}$ is
the restriction of the momentum map $\mathbf{J}$ to
$\mathbf{J}^{-1}(\mathcal{O}_\mu)$, that is,
$\mathbf{J}_{\mathcal{O}_\mu}=\mathbf{J}\cdot i_{\mathcal{O}_\mu}$
and $\omega_{\mathcal{O}_\mu}^+$ is the $+$-symplectic structure on
the orbit $\mathcal{O}_\mu$ given by
\begin{equation}\omega_{\mathcal{O}_\mu}^
+(\nu)(\xi_{\mathfrak{g}^\ast}(\nu),\eta_{\mathfrak{g}^\ast}(\nu))
=<\nu,[\xi,\eta]>,\;\; \forall\;\nu\in\mathcal{O}_\mu, \;
\xi,\eta\in \mathfrak{g}. \label{3.3}\end{equation} The maps
$i_{\mathcal{O}_\mu}:\mathbf{J}^{-1}(\mathcal{O}_\mu)\rightarrow
T^\ast Q$ and
$\pi_{\mathcal{O}_\mu}:\mathbf{J}^{-1}(\mathcal{O}_\mu)\rightarrow
(T^\ast Q)_{\mathcal{O}_\mu}$ are natural injection and the
projection, respectively. The pair $((T^\ast
Q)_{\mathcal{O}_\mu},\omega_{\mathcal{O}_\mu})$ is called the regular
orbit reduced symplectic space of $(T^\ast Q,\omega)$ at $\mu$.\\

Let $H:T^\ast Q\rightarrow \mathbb{R}$ be a $G$-invariant
Hamiltonian, the flow $F_t$ of the Hamiltonian vector field $X_H$
leaves the connected components of
$\mathbf{J}^{-1}(\mathcal{O}_\mu)$ invariant and commutes with the
$G$-action, so it induces a flow $f_t^{\mathcal{O}_\mu}$ on $(T^\ast
Q)_{\mathcal{O}_\mu}$, defined by $f_t^{\mathcal{O}_\mu}\cdot
\pi_{\mathcal{O}_\mu}=\pi_{\mathcal{O}_\mu} \cdot F_t\cdot
i_{\mathcal{O}_\mu}$, and the vector field $X_{h_{\mathcal{O}_\mu}}$
generated by the flow $f_t^{\mathcal{O}_\mu}$ on $((T^\ast
Q)_{\mathcal{O}_\mu},\omega_{\mathcal{O}_\mu})$ is Hamiltonian with
the associated regular orbit reduced Hamiltonian function
$h_{\mathcal{O}_\mu}:(T^\ast Q)_{\mathcal{O}_\mu}\rightarrow
\mathbb{R}$ defined by $h_{\mathcal{O}_\mu}\cdot
\pi_{\mathcal{O}_\mu}= H\cdot i_{\mathcal{O}_\mu}$, and the
Hamiltonian vector fields $X_H$ and $X_{h_{\mathcal{O}_\mu}}$ are
$\pi_{\mathcal{O}_\mu}$-related. Thus, we can define a regular orbit
reducible Hamiltonian system as follows.

\begin{definition}
(Regular Orbit Reducible Hamiltonian System) A 4-tuple $(T^\ast Q,
G,\omega,H)$, where the Hamiltonian $H: T^\ast Q\rightarrow
\mathbb{R} $ is $G$-invariant, is called a regular orbit reducible
Hamiltonian system, if there exists an orbit $\mathcal{O}_\mu, \;
\mu\in\mathfrak{g}^\ast$, where $\mu$ is a regular value of the
momentum map $\mathbf{J}$, such that the regular orbit reduced
system, that is, the 3-tuple $((T^\ast
Q)_{\mathcal{O}_\mu},\omega_{\mathcal{O}_\mu},h_{\mathcal{O}_\mu})$,
where $(T^\ast
Q)_{\mathcal{O}_\mu}=\mathbf{J}^{-1}(\mathcal{O}_\mu)/G$,
$\pi_{\mathcal{O}_\mu}^\ast \omega_{\mathcal{O}_\mu}
=i_{\mathcal{O}_\mu}^\ast\omega-\mathbf{J}_{\mathcal{O}_\mu}^\ast\omega_{\mathcal{O}_\mu}^+$,
$h_{\mathcal{O}_\mu}\cdot \pi_{\mathcal{O}_\mu} =H\cdot
i_{\mathcal{O}_\mu}$, is a Hamiltonian system. Here $((T^\ast
Q)_{\mathcal{O}_\mu},\omega_{\mathcal{O}_\mu})$ is the regular orbit
reduced symplectic space, and the function $h_{\mathcal{O}_\mu}:(T^\ast
Q)_{\mathcal{O}_\mu}\rightarrow \mathbb{R}$ is called the regular
orbit reduced Hamiltonian.
\end{definition}

It is worthy of note that the regular reduced symplectic spaces $((T^\ast
Q)_{\mathcal{O}_\mu},\omega_{\mathcal{O}_\mu})$ and $((T^\ast Q)_\mu,\omega_\mu), $ of the regular orbit
reduced Hamiltonian system and the regular point reduced Hamiltonian system, are different,
and the symplectic forms on the reduced spaces, given by (4.1) for the regular orbit
reduced Hamiltonian system and given by (3.1) for the regular point
reduced Hamiltonian system, are also different.
Since the regular orbit reduced symplectic space
$(T^\ast Q)_{\mathcal{O}_\mu}= \mathbf{J}^{-1}(\mathcal{O}_\mu)/G
\cong \mathbf{J}^{-1}(\mu)/G \times \mathcal{O}_\mu, $ if we give a stronger
assumption condition, that is, for the one-form $\gamma: Q \rightarrow T^*Q$ on $Q,$
assume that $\textmd{Im}(\gamma)\subset
\mathbf{J}^{-1}(\mu), $ and it is $G$-invariant, then in this case for any $V\in TQ, $ and $w\in TT^*Q, $
we have that
$\mathbf{J}_{\mathcal{O}_\mu}^\ast\omega_{\mathcal{O}_\mu}^+(T\gamma
\cdot V, \; w)=0, $ and hence from (4.1), $i_{\mathcal{O}_\mu}^\ast
\omega=\pi_{\mathcal{O}_{\mu}}^\ast \omega_{\mathcal{O}
_\mu}+\mathbf{J}_{\mathcal{O}_\mu}^\ast\omega_{\mathcal{O}_\mu}^+, $
we have that $\pi_{\mathcal{O}_\mu}^*\omega_{\mathcal{O}_\mu}=
i_{\mathcal{O}_\mu}^*\omega= \omega, $ along $\textmd{Im}(\gamma)$.
Thus, we can use the Lemma 2.4 for the regular orbit reduced symplectic form $\omega_{\mathcal{O}_\mu}$.
In particular, note that it is easy to give the wrong results
without the precise analysis for the regular orbit reduction case. For example, by analogizing for the regular point reduction case,
we assume that $\textmd{Im}(\gamma)\subset
\mathbf{J}^{-1}(\mathcal{O}_\mu), $ and it is $G$-invariant, then we can not guarantee
that $\mathbf{J}_{\mathcal{O}_\mu}^\ast\omega_{\mathcal{O}_\mu}^+=0, $ along $\textmd{Im}(\gamma), $
such that the Lemma 2.4 can be used in the following Theorem 4.2 and Theorem 4.3.\\

In the following for the regular orbit reducible Hamiltonian system
$(T^*Q,G,\omega,H)$, by using Lemma 2.4, the regular orbit reduced symplectic form and the reduced
dynamical vector field, we can prove the following two types of Hamilton-Jacobi
theorem for the regular orbit reduced Hamiltonian system.
At first, by using Lemma 2.4 and the orbit reduced symplectic form $\omega_{\mathcal{O}_\mu}$,
and the fact that the one-form $\gamma: Q
\rightarrow T^*Q $ is closed with respect to
$T\pi_Q: TT^* Q \rightarrow TQ, $ and $\textmd{Im}(\gamma)\subset
\mathbf{J}^{-1}(\mu), $ and it is $G$-invariant, we can also prove the following Type I of Hamilton-Jacobi
theorem for the regular orbit reduced Hamiltonian system $((T^\ast
Q)_{\mathcal{O}_\mu}, \omega_{\mathcal{O}_\mu},h_{\mathcal{O}_\mu}
)$.

\begin{theorem} (Type I of Hamilton-Jacobi Theorem for a Regular
Orbit Reduced Hamiltonian System)
For the regular orbit reducible Hamiltonian system
$(T^*Q,G,\omega,H)$, assume that $\gamma: Q \rightarrow T^*Q$ is an
one-form on $Q$, and $X_H^\gamma = T\pi_{Q}\cdot X_H \cdot \gamma$,
where $X_{H}$ is the dynamical vector field of $(T^*Q,G,\omega,H)$.
Moreover, assume that $\mu \in \mathfrak{g}^\ast $ is a regular
value of the momentum map $\mathbf{J}$, and $\mathcal{O}_\mu, \;
(\mu \in \mathfrak{g}^\ast), $ is the regular reducible orbit of the
Hamiltonian system, and $\textmd{Im}(\gamma)\subset
\mathbf{J}^{-1}(\mu), $ and it is $G$-invariant,
$\bar{\gamma}=\pi_{\mathcal{O}_\mu}(\gamma): Q \rightarrow (T^*
Q)_{\mathcal{O}_\mu}. $ If the one-form $\gamma: Q \rightarrow T^*Q
$ is closed with respect to $T\pi_Q: TT^* Q \rightarrow TQ, $ then
$\bar{\gamma}$ is a solution of the equation $T\bar{\gamma}\cdot
X_H^\gamma= X_{h_{\mathcal{O}_\mu}}\cdot \bar{\gamma}, $ which is
called the Type I of Hamilton-Jacobi equation for the regular orbit
reduced Hamiltonian system $((T^\ast Q)_{\mathcal{O}_\mu},
\omega_{\mathcal{O}_\mu},h_{\mathcal{O}_\mu})$. Here the maps involved in
the theorem are shown in the following Diagram-5.
\begin{center}
\hskip 0cm \xymatrix{ \mathbf{J}^{-1}(\mathcal{O}_\mu)
\ar[r]^{i_{\mathcal{O}_\mu}} & T^* Q \ar[r]^{\pi_Q}
& Q \ar[d]_{X_H^\gamma} \ar[r]^{\gamma} & T^*Q \ar[d]_{X_H} \ar[r]^{\pi_{\mathcal{O}_\mu}}
& (T^* Q)_{\mathcal{O}_\mu} \ar[d]_{X_{h_{\mathcal{O}_\mu}}} \\
  & T(T^*Q)  & TQ \ar[l]^{T\gamma} & T(T^*Q) \ar[l]^{T\pi_Q} \ar[r]_{T\pi_{\mathcal{O}_\mu}}
  & T(T^* Q)_{\mathcal{O}_\mu} }
\end{center}
$$\mbox{Diagram-5}$$
\end{theorem}

\noindent{\bf Proof: } At first, from Theorem 2.5, we know that
$\gamma$ is a solution of the Type I of Hamilton-Jacobi equation
$T\gamma\cdot X_H^\gamma= X_H\cdot \gamma .$
Next, we note that the regular orbit reduced symplectic space
$(T^\ast Q)_{\mathcal{O}_\mu}= \mathbf{J}^{-1}(\mathcal{O}_\mu)/G
\cong \mathbf{J}^{-1}(\mu)/G \times \mathcal{O}_\mu, $ with the
symplectic form $\omega_{\mathcal{O}_\mu}$ uniquely characterized by
the relation $i_{\mathcal{O}_\mu}^\ast
\omega=\pi_{\mathcal{O}_{\mu}}^\ast \omega_{\mathcal{O}
_\mu}+\mathbf{J}_{\mathcal{O}_\mu}^\ast\omega_{\mathcal{O}_\mu}^+. $
Since $\textmd{Im}(\gamma)\subset \mathbf{J}^{-1}(\mu), $ and it is
$G$-invariant, in this case for any $V\in TQ, $ and $w\in TT^*Q, $
we have that
$\mathbf{J}_{\mathcal{O}_\mu}^\ast\omega_{\mathcal{O}_\mu}^+(T\gamma
\cdot V, \; w)=0, $ and hence
$\pi_{\mathcal{O}_\mu}^*\omega_{\mathcal{O}_\mu}=
i_{\mathcal{O}_\mu}^*\omega= \omega, $ along $\textmd{Im}(\gamma)$.
By using the reduced symplectic form $\omega_{\mathcal{O}_\mu}$, and if
we take that $v= X_H\cdot \gamma \in TT^* Q,$ and
for any $w \in TT^* Q, \; T\pi_{Q}(w)\neq 0,$ and $T\pi_{\mathcal{O}_\mu} (w) \neq 0,
$ from Lemma 2.4(ii) we have that
\begin{align*}
& \omega_{\mathcal{O}_\mu}(T\bar{\gamma} \cdot X_H^\gamma, \;
T\pi_{\mathcal{O}_\mu} \cdot w) =
\omega_{\mathcal{O}_\mu}(T(\pi_{\mathcal{O}_\mu} \cdot \gamma) \cdot
X_H^\gamma, \; T\pi_{\mathcal{O}_\mu} \cdot w )\\
& = \pi_{\mathcal{O}_\mu}^*\omega_{\mathcal{O}_\mu}(T\gamma \cdot
X_H^\gamma, \; w) = \omega(T(\gamma \cdot \pi_Q)\cdot X_H\cdot
\gamma, \; w)\\
& = \omega(X_H\cdot \gamma, \; w-T(\gamma \cdot
\pi_Q)\cdot w)-\mathbf{d}\gamma(T\pi_{Q}(X_H\cdot \gamma), \; T\pi_{Q}(w))\\
& = \omega_{\mathcal{O}_\mu}(X_{h_{\mathcal{O}_\mu}} \cdot
\bar{\gamma}, \; T\pi_{\mathcal{O}_\mu} \cdot w)-
\omega_{\mathcal{O}_\mu}(X_{h_{\mathcal{O}_\mu}} \cdot
\bar{\gamma}, \; T\bar{\gamma} \cdot T\pi_{Q}(w))-\mathbf{d}\gamma(T\pi_{Q}(X_H\cdot \gamma), \; T\pi_{Q}(w)),
\end{align*}
in which we have used that $ T\pi_{\mathcal{O}_\mu} \cdot X_H= X_{
h_{\mathcal{O}_\mu}}. $
Since the one-form $\gamma: Q \rightarrow T^*Q $ is closed with respect to
$T\pi_Q: TT^* Q \rightarrow TQ, $ then we have that $
\mathbf{d}\gamma(T\pi_{Q}(X_H\cdot \gamma), \; T\pi_{Q}(w))=0, $
and hence
\begin{equation}
 \omega_{\mathcal{O}_\mu}(T\bar{\gamma} \cdot X_H^\gamma, \; T\pi_{\mathcal{O}_\mu} \cdot w)-
\omega_{\mathcal{O}_\mu}(X_{h_{\mathcal{O}_\mu}} \cdot \bar{\gamma}, \; T\pi_{\mathcal{O}_\mu} \cdot w)
 = - \omega_{\mathcal{O}_\mu}(X_{h_{\mathcal{O}_\mu}} \cdot
\bar{\gamma}, \; T\bar{\gamma} \cdot T\pi_{Q}(w)).
\end{equation}
If $\bar{\gamma}$ satisfies the equation $T\bar{\gamma}\cdot X_H^\gamma= X_{h_{\mathcal{O}_\mu}}\cdot \bar{\gamma}, $
from Lemma 2.4(i) we can obtain that
\begin{align*}
- \omega_{\mathcal{O}_\mu}(X_{h_{\mathcal{O}_\mu}} \cdot
\bar{\gamma}, \; T\bar{\gamma} \cdot T\pi_{Q}(w))
& = -\omega_{\mathcal{O}_\mu} (T\bar{\gamma} \cdot X_H^\gamma, \; T\bar{\gamma} \cdot T\pi_{Q}(w))\\
& = -\bar{\gamma}^*\omega_{\mathcal{O}_\mu} (T\pi_{Q} \cdot X_{H}\cdot\gamma, \; T\pi_Q(w))\\
& = -\gamma^* \cdot \pi_{\mathcal{O}_\mu}^*\omega_{\mathcal{O}_\mu}(T\pi_{Q} \cdot X_{H}\cdot\gamma, \; T\pi_{Q}(w))\\
& = -\gamma^*\omega( T\pi_{Q}(X_{H}\cdot\gamma), \; T\pi_{Q}(w))\\
& = \textbf{d}\gamma(T\pi_{Q}( X_{H}\cdot\gamma ), \; T\pi_{Q}(w))=0.
\end{align*}
Because the reduced symplectic form $\omega_{\mathcal{O}_\mu}$
is non-degenerate, the left side of (4.3) equals zero, only when
$\bar{\gamma}$ satisfies the equation
$T\bar{\gamma}\cdot X_H^\gamma= X_{h_{\mathcal{O}_\mu}}\cdot \bar{\gamma}.$
Thus, if the one-form $\gamma: Q \rightarrow T^*Q $ is closed with respect to
$T\pi_Q: TT^* Q \rightarrow TQ, $ then $\bar{\gamma}$ must be
a solution of the Type I of Hamilton-Jacobi equation
$T\bar{\gamma}\cdot X_H^\gamma= X_{h_{\mathcal{O}_\mu}}\cdot \bar{\gamma}.$
\hskip 0.3cm $\blacksquare$\\

Next, for any $G$-invariant symplectic map $\varepsilon: T^* Q \rightarrow T^* Q $,
we can also prove the following Type II of
Hamilton-Jacobi theorem for the regular orbit reduced Hamiltonian system $((T^\ast
Q)_{\mathcal{O}_\mu}, \omega_{\mathcal{O}_\mu},h_{\mathcal{O}_\mu}
)$.

\begin{theorem} (Type II of Hamilton-Jacobi Theorem for a Regular
Orbit Reduced Hamiltonian System)
For the regular orbit reducible Hamiltonian system
$(T^*Q,G,\omega,H)$, assume that $\gamma: Q \rightarrow T^*Q$ is an one-form on $Q$,
and $\lambda=\gamma \cdot \pi_{Q}: T^* Q
\rightarrow T^* Q $, and for any $G$-invariant symplectic map $\varepsilon: T^* Q \rightarrow T^* Q $,
denote by $X_H^\varepsilon = T\pi_{Q}\cdot
X_H \cdot \varepsilon$, where $X_{H}$ is the dynamical vector field of $(T^*Q,G,\omega,H)$.
Moreover, assume that $\mu \in \mathfrak{g}^\ast $ is a regular value of the momentum
map $\mathbf{J}$, and $\mathcal{O}_\mu, \; (\mu \in \mathfrak{g}^\ast), $ is the
regular reducible orbit of the Hamiltonian system, and
$\textmd{Im}(\gamma)\subset \mathbf{J}^{-1}(\mu), $ and it is
$G$-invariant, and $\varepsilon(\mathbf{J}^{-1}(\mathcal{O}_\mu)) \subset \mathbf{J}^{-1}(\mathcal{O}_\mu). $
Denote by $\bar{\gamma}=\pi_{\mathcal{O}_\mu}(\gamma): Q
\rightarrow (T^* Q)_{\mathcal{O}_\mu} $,
$\bar{\lambda}=\pi_{\mathcal{O}_\mu}(\lambda): \mathbf{J}^{-1}(\mathcal{O}_\mu) (\subset T^*Q) \rightarrow (T^*
Q)_{\mathcal{O}_\mu} $, and $\bar{\varepsilon}=\pi_{\mathcal{O}_\mu}(\varepsilon):
\mathbf{J}^{-1}(\mathcal{O}_\mu) (\subset T^*Q) \rightarrow (T^*
Q)_{\mathcal{O}_\mu} $. Then $\varepsilon$ and $\bar{\varepsilon}$ satisfy the equation
$T\bar{\varepsilon}\cdot X_{h_{\mathcal{O}_\mu}\cdot\bar{\varepsilon}}= T\bar{\lambda} \cdot X_H\cdot\varepsilon, $
if and only if they satisfy the equation $T\bar{\gamma}\cdot X_H^\varepsilon=
X_{h_{\mathcal{O}_\mu}}\cdot \bar{\varepsilon}$, where
$ X_{h_{\mathcal{O}_\mu} \cdot
\bar{\varepsilon}} \in TT^*Q $ is the Hamiltonian vector field of the
function $h_{\mathcal{O}_\mu} \cdot \bar{\varepsilon}: T^*Q\rightarrow
\mathbb{R}. $
The equation $T\bar{\gamma}\cdot X_H^\varepsilon=
X_{h_{\mathcal{O}_\mu}}\cdot \bar{\varepsilon}$ is called the Type II of Hamilton-Jacobi
equation for the regular orbit reduced Hamiltonian system $((T^\ast
Q)_{\mathcal{O}_\mu}, \omega_{\mathcal{O}_\mu},h_{\mathcal{O}_\mu})$.
Here the maps involved in the theorem are shown in the following Diagram-6.

\begin{center}
\hskip 0cm \xymatrix{ \mathbf{J}^{-1}(\mathcal{O}_\mu)
\ar[r]^{i_{\mathcal{O}_\mu}} & T^* Q \ar[d]_{X_{H\cdot \varepsilon}}
\ar[dr]^{X_H^\varepsilon} \ar[r]^{\pi_Q}
& Q \ar[r]^{\gamma} & T^*Q \ar[d]_{X_H} \ar[dr]^{X_{h_{\mathcal{O}_\mu}} \cdot\bar{\varepsilon}} \ar[r]^{\pi_{\mathcal{O}_\mu}}
& (T^* Q)_{\mathcal{O}_\mu} \ar[d]^{X_{h_{\mathcal{O}_\mu}}} \\
  & T(T^*Q)  & TQ \ar[l]^{T\gamma} & T(T^*Q) \ar[l]^{T\pi_Q} \ar[r]_{T\pi_{\mathcal{O}_\mu}} & T(T^* Q)_{\mathcal{O}_\mu} }
\end{center}
$$\mbox{Diagram-6}$$
\end{theorem}

\noindent{\bf Proof: } Note that the regular orbit reduced symplectic space
$(T^\ast Q)_{\mathcal{O}_\mu}= \mathbf{J}^{-1}(\mathcal{O}_\mu)/G
\cong (\mathbf{J}^{-1}(\mu)/G) \times \mathcal{O}_\mu, $ with the reduced
symplectic form $\omega_{\mathcal{O}_\mu}$ uniquely characterized by
the relation $i_{\mathcal{O}_\mu}^\ast
\omega=\pi_{\mathcal{O}_{\mu}}^\ast \omega_{\mathcal{O}
_\mu}+\mathbf{J}_{\mathcal{O}_\mu}^\ast\omega_{\mathcal{O}_\mu}^+. $
Since $\textmd{Im}(\gamma)\subset \mathbf{J}^{-1}(\mu), $ and it is
$G$-invariant, in this case for any $V\in TQ, $ and $w\in TT^*Q, $
we have that
$\mathbf{J}_{\mathcal{O}_\mu}^\ast\omega_{\mathcal{O}_\mu}^+(T\gamma
\cdot V, \; w)=0, $ and hence
$\pi_{\mathcal{O}_\mu}^*\omega_{\mathcal{O}_\mu}=
i_{\mathcal{O}_\mu}^*\omega= \omega, $ along $\textmd{Im}(\gamma)$.
By using the reduced symplectic form $\omega_{\mathcal{O}_\mu}$, and if
we take that $v= X_H\cdot \varepsilon \in TT^* Q, $ and
for any $w \in TT^* Q, \; T\bar{\lambda}(w)\neq 0, $ and $T\pi_{\mathcal{O}_\mu} (w) \neq 0,
$ from Lemma 2.4 we have that
\begin{align*}
& \omega_{\mathcal{O}_\mu}(T\bar{\gamma} \cdot X_H^\varepsilon, \;
T\pi_{\mathcal{O}_\mu} \cdot w) =
\omega_{\mathcal{O}_\mu}(T(\pi_{\mathcal{O}_\mu} \cdot \gamma) \cdot
X_H^\varepsilon, \; T\pi_{\mathcal{O}_\mu} \cdot w )\\
& = \pi_{\mathcal{O}_\mu}^*\omega_{\mathcal{O}_\mu}(T\gamma \cdot
X_H^\varepsilon, \; w) = \omega(T(\gamma \cdot \pi_Q)\cdot X_H\cdot
\varepsilon, \; w)\\
& = \omega(X_H\cdot \varepsilon, \; w-T(\gamma \cdot
\pi_Q)\cdot w)-\mathbf{d}\gamma(T\pi_{Q}(X_H\cdot \varepsilon), \; T\pi_{Q}(w))\\
& =\omega(X_H\cdot \varepsilon, \; w) - \omega(X_H\cdot \varepsilon, \;
T\lambda\cdot w)+ \lambda^*\omega(X_H\cdot \varepsilon, \; w)\\
& = \omega_{\mathcal{O}_\mu}(X_{h_{\mathcal{O}_\mu}} \cdot
\bar{\varepsilon}, \; T\pi_{\mathcal{O}_\mu} \cdot w)-
\omega_{\mathcal{O}_\mu}(X_{h_{\mathcal{O}_\mu}} \cdot
\bar{\varepsilon}, \; T\bar{\lambda} \cdot w)
+ \omega_{\mathcal{O}_\mu}(T\bar{\lambda} \cdot X_H\cdot \varepsilon, \; T\bar{\lambda} \cdot w),
\end{align*}
in which we have used that $ T\pi_{\mathcal{O}_\mu} \cdot X_H= X_{
h_{\mathcal{O}_\mu}}. $
Note that $\varepsilon: T^* Q
\rightarrow T^* Q $ is symplectic, and $\pi_{\mathcal{O}_\mu}^*\omega_{\mathcal{O}_\mu}=
i_{\mathcal{O}_\mu}^*\omega= \omega, $ along $\textmd{Im}(\gamma)$,
and hence $\bar{\varepsilon}=
\pi_{\mathcal{O}_\mu}\cdot \varepsilon: \mathbf{J}^{-1}(\mathcal{O}_\mu) (\subset T^* Q)
\rightarrow (T^* Q)_{\mathcal{O}_\mu}$ is also symplectic along
$\textmd{Im}(\gamma)$, and $X_{h_{\mathcal{O}_\mu}}\cdot
\bar{\varepsilon} = T\bar{\varepsilon} \cdot X_{h_{\mathcal{O}_\mu} \cdot
\bar{\varepsilon}}, $ along $\textmd{Im}(\gamma)\cap\textmd{Im}(\varepsilon). $
From the above arguments, we can obtain that
\begin{align*}
& \omega_{\mathcal{O}_\mu}(T\bar{\gamma} \cdot X_H^\varepsilon, \;
T\pi_{\mathcal{O}_\mu} \cdot w)-
\omega_{\mathcal{O}_\mu}(X_{h_{\mathcal{O}_\mu}} \cdot
\bar{\varepsilon}, \; T\pi_{\mathcal{O}_\mu} \cdot w)\\
& = \omega_{\mathcal{O}_\mu}(T\bar{\lambda}\cdot X_H\cdot \varepsilon, \; T\bar{\lambda}\cdot w)
-\omega_{\mathcal{O}_\mu}(T\bar{\varepsilon} \cdot X_{h_{\mathcal{O}_\mu} \cdot \bar{\varepsilon}}, \; T\bar{\lambda} \cdot w)\\
& = \omega_{\mathcal{O}_\mu}(T\bar{\lambda}\cdot X_H\cdot \varepsilon
- T\bar{\varepsilon} \cdot X_{h_{\mathcal{O}_\mu} \cdot \bar{\varepsilon}}, \; T\bar{\lambda}\cdot w).
\end{align*}
Because the reduced symplectic form $\omega_{\mathcal{O}_\mu}$ is non-degenerate, it follows that
$T\bar{\gamma}\cdot X_H^\varepsilon=
X_{h_{\mathcal{O}_\mu}}\cdot \bar{\varepsilon}$ is equivalent to
$T\bar{\varepsilon}\cdot X_{h_{\mathcal{O}_\mu}\cdot\bar{\varepsilon}}= T\bar{\lambda} \cdot X_H\cdot\varepsilon. $
Thus, we know that the $\varepsilon$ and $\bar{\varepsilon}$ satisfy the equation
$T\bar{\varepsilon}\cdot X_{h_{\mathcal{O}_\mu}\cdot\bar{\varepsilon}}= T\bar{\lambda} \cdot X_H\cdot\varepsilon, $
if and only if they satisfy the Type II of Hamilton-Jacobi equation $T\bar{\gamma}\cdot X_H^\varepsilon=
X_{h_{\mathcal{O}_\mu}}\cdot \bar{\varepsilon}. $
\hskip 0.3cm $\blacksquare$\\

Moreover, for the regular orbit reducible Hamiltonian system
$(T^*Q,G,\omega,H)$, we know that the dynamical vector fields
$X_{H}$ and $X_{h_{\mathcal{O}_\mu}}$ are
$\pi_{\mathcal{O}_\mu}$-related, that is, $X_{
h_{\mathcal{O}_\mu}}\cdot \pi_{\mathcal{O}_\mu} =
T\pi_{\mathcal{O}_\mu}\cdot X_{H}\cdot i_{\mathcal{O}_\mu}. $ Then
we can also prove the following Theorem 4.4, which states the
relationship between the solution of Type II of Hamilton-Jacobi equation and the
regular orbit reduction.

\begin{theorem}
For the regular orbit reducible Hamiltonian system
$(T^*Q,G,\omega,H)$, assume that $\gamma: Q \rightarrow T^*Q$ is an
one-form on $Q$, and $\lambda=\gamma \cdot \pi_{Q}: T^* Q
\rightarrow T^* Q $, and $\varepsilon: T^* Q \rightarrow T^* Q $ is
a $G$-invariant symplectic map. Moreover, assume that $\mu \in
\mathfrak{g}^\ast $ is a regular value of the momentum map
$\mathbf{J}$, and $\mathcal{O}_\mu, \; (\mu \in \mathfrak{g}^\ast),
$ is the regular reducible orbit of the Hamiltonian system, and
$\textmd{Im}(\gamma)\subset \mathbf{J}^{-1}(\mu), $ and it is
$G$-invariant, and $\varepsilon(\mathbf{J}^{-1}(\mathcal{O}_\mu))
\subset \mathbf{J}^{-1}(\mathcal{O}_\mu). $ Denote by
$\bar{\gamma}=\pi_{\mathcal{O}_\mu}(\gamma): Q \rightarrow (T^*
Q)_{\mathcal{O}_\mu} $,
$\bar{\lambda}=\pi_{\mathcal{O}_\mu}(\lambda):
\mathbf{J}^{-1}(\mathcal{O}_\mu) (\subset T^*Q) \rightarrow (T^*
Q)_{\mathcal{O}_\mu} $, and
$\bar{\varepsilon}=\pi_{\mathcal{O}_\mu}(\varepsilon):
\mathbf{J}^{-1}(\mathcal{O}_\mu) (\subset T^*Q) \rightarrow (T^*
Q)_{\mathcal{O}_\mu} $. Then $\varepsilon$ is a solution of the Type
II of Hamilton-Jacobi equation $T\gamma\cdot X_H^\varepsilon=
X_H\cdot \varepsilon, $ for the regular orbit reducible Hamiltonian
system $(T^*Q,G,\omega,H), $ if and only if $\varepsilon$ and
$\bar{\varepsilon}$ satisfy the Type II of Hamilton-Jacobi equation
$T\bar{\gamma}\cdot X_H^\varepsilon=
X_{h_{\mathcal{O}_\mu}\cdot\bar{\varepsilon}}, $ for the regular
orbit reduced Hamiltonian system $((T^\ast Q)_{\mathcal{O}_\mu},
\omega_{\mathcal{O}_\mu},h_{\mathcal{O}_\mu})$.
\end{theorem}

\noindent{\bf Proof: } Note that $\textmd{Im}(\gamma)\subset \mathbf{J}^{-1}(\mu), $ and it
is $G$-invariant, in this case, $\pi_{\mathcal{O}_\mu}^*\omega_{\mathcal{O}_\mu}=
i_{\mathcal{O}_\mu}^*\omega= \omega, $ along $\textmd{Im}(\gamma)$.
Since the dynamical vector fields
$X_{H}$ and $X_{h_{\mathcal{O}_\mu}}$ are $\pi_{\mathcal{O}_\mu}$-related, that is,
$X_{h_{\mathcal{O}_\mu}}\cdot \pi_{\mathcal{O}_\mu}= T\pi_{\mathcal{O}_\mu}\cdot X_{H}\cdot i_{\mathcal{O}_\mu}, $
and by using the reduced symplectic form $\omega_{\mathcal{O}_\mu}$, we have that
\begin{align*}
& \omega_{\mathcal{O}_\mu}(T\bar{\gamma}\cdot X_H^\varepsilon-X_{h_{\mathcal{O}_\mu}}\cdot
\bar{\varepsilon}, \; T\pi_{\mathcal{O}_\mu} \cdot w)\\
& = \omega_{\mathcal{O}_\mu}(T\bar{\gamma} \cdot X_H^\varepsilon, \;
T\pi_{\mathcal{O}_\mu} \cdot w)-
\omega_{\mathcal{O}_\mu}(X_{h_{\mathcal{O}_\mu}} \cdot
\bar{\varepsilon}, \; T\pi_{\mathcal{O}_\mu} \cdot w)\\
& = \omega_{\mathcal{O}_\mu}(T\pi_{\mathcal{O}_\mu} \cdot T\gamma \cdot X_H^\varepsilon, \; T\pi_{\mathcal{O}_\mu} \cdot w)
-\omega_{\mathcal{O}_\mu}(X_{h_{\mathcal{O}_\mu}} \cdot T\pi_{\mathcal{O}_\mu} \cdot \varepsilon, \; T\pi_{\mathcal{O}_\mu} \cdot w)\\
& = \pi_{\mathcal{O}_\mu}^*\omega_{\mathcal{O}_\mu}(T\gamma\cdot X_H^\varepsilon, \; w)
- \omega_{\mathcal{O}_\mu}(T\pi_{\mathcal{O}_\mu}\cdot X_{H}\cdot \varepsilon, \; T\pi_{\mathcal{O}_\mu} \cdot w)\\
& = \pi_{\mathcal{O}_\mu}^*\omega_{\mathcal{O}_\mu}(T\gamma\cdot X_H^\varepsilon, \; w)
-\pi_{\mathcal{O}_\mu}^*\omega_{\mathcal{O}_\mu}(X_{H}\cdot \varepsilon, \; w)\\
& = \omega(T\gamma\cdot X_H^\varepsilon, \; w)-\omega(X_{H}\cdot \varepsilon, \; w)\\
& = \omega(T\gamma\cdot X_H^\varepsilon- X_{H}\cdot \varepsilon, \; w).
\end{align*}
Because the symplectic form $\omega$ and the reduced symplectic form $\omega_{\mathcal{O}_\mu}$ are non-degenerate,
it follows that the equation
$T\bar{\gamma}\cdot X_H^\varepsilon= X_{h_{\mathcal{O}_\mu}}\cdot \bar{\varepsilon}, $
is equivalent to the equation
$T\gamma\cdot X_H^\varepsilon= X_H\cdot \varepsilon. $ Thus,
$\varepsilon$ is a solution of the Type II of Hamilton-Jacobi equation
$T\gamma\cdot X_H^\varepsilon= X_H\cdot \varepsilon, $ for the
regular orbit reducible Hamiltonian system $(T^*Q,G,\omega,H), $ if and only if
$\varepsilon$ and $\bar{\varepsilon} $ satisfy the Type II of Hamilton-Jacobi equation
$T\bar{\gamma}\cdot X_H^\varepsilon= X_{h_{\mathcal{O}_\mu}}\cdot \bar{\varepsilon}, $ for the
regular orbit reduced Hamiltonian system $((T^\ast Q)_{\mathcal{O}_\mu},
\omega_{\mathcal{O}_\mu}, h_{\mathcal{O}_\mu}). $
\hskip 0.3cm $\blacksquare$\\

\begin{remark}
If $(T^\ast Q, \omega)$ is a connected symplectic manifold, and
$\mathbf{J}:T^\ast Q\rightarrow \mathfrak{g}^\ast$ is a
non-equivariant momentum map with a non-equivariance group
one-cocycle $\sigma: G\rightarrow \mathfrak{g}^\ast$, which is
defined by $\sigma(g):=\mathbf{J}(g\cdot
z)-\operatorname{Ad}^\ast_{g^{-1}}\mathbf{J}(z)$, where $g\in G$ and
$z\in T^\ast Q$. Then we know that $\sigma$ produces a new affine
action $\Theta: G\times \mathfrak{g}^\ast \rightarrow
\mathfrak{g}^\ast $ defined by
$\Theta(g,\mu):=\operatorname{Ad}^\ast_{g^{-1}}\mu + \sigma(g)$,
where $\mu \in \mathfrak{g}^\ast$, with respect to which the given
momentum map $\mathbf{J}$ is equivariant. Assume that $G$ acts
freely and properly on $T^\ast Q$, and $\tilde{G}_\mu$ denotes the
isotropy subgroup of $\mu \in \mathfrak{g}^\ast$ relative to this
affine action $\Theta$, and $\mathcal{O}_\mu= G\cdot \mu
\subset \mathfrak{g}^\ast$ denotes the G-orbit of the point $\mu$ with respect to the action $\Theta$,
and $\mu$ is a regular value of $\mathbf{J}$.
Then the quotient space
$(T^\ast Q)_{\mathcal{O}_\mu}=\mathbf{J}^{-1}(\mathcal{O}_\mu)/ G $
is a symplectic manifold with the symplectic form
$\omega_{\mathcal{O}_\mu}$ uniquely characterized by $(4.1)$,
see Ortega and Ratiu \cite{orra04}.
Moreover, in this case,
for the given regular orbit reducible Hamiltonian system
$(T^*Q,G,\omega,H)$, we can also prove the Type I and Type II of
Hamilton-Jacobi theorem for the regular orbit reduced Hamiltonian
system $((T^\ast Q)_{\mathcal{O}_\mu}, \omega_{\mathcal{O}_\mu},
h_{\mathcal{O}_\mu} ), $ by using the above similar ways.
\end{remark}

\begin{remark}
It is worthy of note that, the one-form $\gamma: Q \rightarrow T^*Q$
may not be given by a generating function of a symplectic map, and
hence the formulations of Type I and Type II of Hamilton-Jacobi
equations for the regular orbit reducible Hamiltonian systems,
given by Theorem $4.2$ and Theorem $4.3$, have more extensive sense. On the
other hand, if $\gamma$ is a solution of the classical
Hamilton-Jacobi equation, that is, $X_H\cdot \gamma=0, $ then
$X_H^\gamma= T\pi_Q\cdot X_H\cdot \gamma=0,$ and hence from the Type
I of Hamilton-Jacobi equation, we have that $X_{h_{\mathcal{O}_\mu}}\cdot \bar{\gamma}= T\bar{\gamma}\cdot
X_H^\gamma=0. $ The equation $X_{h_{\mathcal{O}_\mu}}\cdot
\bar{\gamma}=0$ is called the classical Hamilton-Jacobi equation for
the regular orbit reduced Hamiltonian system $((T^\ast
Q)_{\mathcal{O}_\mu}, \omega_{\mathcal{O}_\mu},h_{\mathcal{O}_\mu}),
$ which shows that the dynamical vector field of the
system $((T^\ast Q)_{\mathcal{O}_\mu},
\omega_{\mathcal{O}_\mu},h_{\mathcal{O}_\mu})$ is degenerate along
$\bar{\gamma}$. In addition, for a symplectic map $\varepsilon: T^* Q \rightarrow
T^* Q $, if $X_H\cdot \varepsilon=0, $ then from the Type II of
Hamilton-Jacobi equation, we also know that $X_{h_{\mathcal{O}_\mu}}\cdot
\bar{\varepsilon}=0 $ is equivalent to $X_{h_{\mathcal{O}_\mu} \cdot
\bar{\varepsilon}}=0.$
\end{remark}

\section{Applications}

In this section, as the applications of the above theoretical
results, we consider a regular point reducible Hamiltonian system on
a Lie group, and give the Hamilton-Jacobi theorems and two types of
Lie-Poisson Hamilton-Jacobi equation for the regular point reduced
system. In particular, we show the Type I and Type II of Lie-Poisson
Hamilton-Jacobi equations for the regular point reduced rigid body
and heavy top systems, respectively. We shall follow the notations
and conventions introduced in Marsden et al. \cite{mamora90},
Marsden and Ratiu \cite{mara99}, Ortega and Ratiu \cite{orra04}, and
Marsden et al. \cite{mawazh10}.

\subsection{Lie-Poisson Hamilton-Jacobi Equation}

Let $G$ be a Lie group with Lie algebra $\mathfrak{g}$ and $T^\ast
G$ its cotangent bundle with the canonical symplectic form $\omega$.
A Hamiltonian system on $G$ is a 3-tuple $(T^\ast G,\omega, H )$,
where the function $H: T^\ast G \rightarrow \mathbb{R}$ is a
Hamiltonian, and has the associated Hamiltonian vector field $X_H$.
At first, for the Lie group $G$, the left and right translation on
$G$, defined by the map $L_g: G \rightarrow G, \; h \mapsto gh $ and
$R_g: G \rightarrow G, \; h \mapsto hg $, for someone $g \in G$,
induce the left and the right action of $G$ on itself. Let $I_g: G \to
G$; $I_g(h)=ghg^{-1}=L_g\cdot R_{g^{-1}}(h)$, for $g,h\in G$, be the
inner automorphism on $G$. The adjoint representation of the Lie group
$G$ is defined by $\operatorname{Ad}_g=T_eI_g = T_{g^{-1}}L_g \cdot
T_e R_{g^{-1}}:\mathfrak{g}\to \mathfrak{g}$. The coadjoint
representation is given by
$\operatorname{Ad}_{g^{-1}}^\ast:\mathfrak{g}^\ast\to
\mathfrak{g}^\ast$, where $\operatorname{Ad}_{g^{-1}}^\ast$ is the
dual of the linear map $\operatorname{Ad}_{g^{-1}}$, defined by
$\langle
\operatorname{Ad}_{g^{-1}}^\ast(\mu),\xi\rangle=\langle\mu,\operatorname{Ad}_{g^{-1}}(\xi)\rangle$,
where $\mu\in\mathfrak{g}^\ast$, $\xi\in\mathfrak{g}$ and
$\langle,\rangle$ denotes the pairing between $\mathfrak{g}^\ast$
and $\mathfrak{g}$. Since the coadjoint representation
$\operatorname{Ad}_{g^{-1}}^\ast:\mathfrak{g}^\ast\to
\mathfrak{g}^\ast$ can induce a left coadjoint action of $G$ on
$\mathfrak{g}^\ast$, the coadjoint orbit $\mathcal{O}_\mu$ of this
action through $\mu\in \mathfrak{g}^\ast$ is a subset of
$\mathfrak{g}^\ast$ defined by
$\mathcal{O}_\mu:=\{\operatorname{Ad}_{g^{-1}}^\ast(\mu)\in
\mathfrak{g}^\ast|g\in G\}$, and $\mathcal{O}_\mu$ is an immersed
submanifold of $\mathfrak{g}^\ast$. We know that
$\mathfrak{g}^\ast$ is a Poisson manifold with respect to the
$(\pm)$-Lie-Poisson bracket $\{\cdot,\cdot\}_\pm$ defined by
\begin{equation}
\{f,g\}_\pm(\mu):=\pm<\mu,[\frac{\delta f}{\delta \mu},\frac{\delta
g}{\delta\mu}]>,\;\; \forall f,g\in C^\infty(\mathfrak{g}^\ast),\;\;
\mu\in \mathfrak{g}^\ast,\label{5.1}
\end{equation}
where the element $\frac{\delta f}{\delta \mu}\in\mathfrak{g}$ is
defined by the equality $<v,\frac{\delta f}{\delta
\mu}>:=Df(\mu)\cdot v$, for any $v\in \mathfrak{g}^\ast$, see
Marsden and Ratiu \cite{mara99}. Thus, for the coadjoint orbit
$\mathcal{O}_\mu, \; \mu\in\mathfrak{g}^\ast$, the orbit symplectic
form can be defined by
\begin{equation}\omega_{\mathcal{O}_\mu}^\pm(\nu)(\operatorname{ad}_\xi^\ast(\nu),
\operatorname{ad}_\eta^\ast(\nu))=\pm \langle\nu,[\xi,\eta]\rangle,
\qquad \forall\; \xi,\eta\in\mathfrak{g}, \;\;
\nu\in\mathcal{O}_\mu\subset\mathfrak{g}^\ast,
\label{5.2}
\end{equation}
which are coincide with the restriction of the Lie-Poisson brackets
on $\mathfrak{g}^\ast$ to the coadjoint orbit $\mathcal{O}_\mu$.
From the symplectic stratification theorem we know that a finite
dimensional Poisson manifold is the disjoint union of its symplectic
leaves, and its each symplectic leaf is an injective immersed
Poisson submanifold whose induced Poisson structure is symplectic.
In consequence, when $\mathfrak{g}^\ast$ is endowed one of the Lie
Poisson structures $\{\cdot,\cdot\}_\pm$, the symplectic leaves of
the Poisson manifolds $(\mathfrak{g}^\ast,\{\cdot,\cdot\}_\pm)$
coincide with the connected components of the orbits of the elements
in $\mathfrak{g}^\ast$ under the coadjoint action. From Abraham and
Marsden \cite{abma78}, we know that the coadjoint orbit
$(\mathcal{O}_\mu, \omega_{\mathcal{O}_\mu}^{-}), \; \mu\in
\mathfrak{g}^\ast,$ is symplectically diffeomorphic to a regular
point reduced space $((T^\ast G)_\mu,\omega_\mu)$ of $T^*G$.\\

We now identify $T^\ast G$ and $G\times\mathfrak{g}^\ast$ locally, by using
the left translation. In fact, the map $\tilde{\lambda}: T^\ast G
\rightarrow G \times \mathfrak{g}^\ast, \;
\tilde{\lambda}(\alpha_g):=(g,(T_eL_g)^\ast \alpha_g)$, for any $\alpha_g
\in T^\ast_g G$, which defines a vector bundle isomorphism usually
referred to as the local left trivialization of $T^\ast G$. In the same
way, we can also identify tangent bundle $TG$ and
$G\times\mathfrak{g}$ locally, by using the left translation. In consequence,
we can consider the Lagrangian $L(g,\xi):TG \cong G\times
\mathfrak{g}\to \mathbb{R}$, which is usual the kinetic minus the
potential energy of the system, where $(g,\xi)\in
G\times\mathfrak{g}$, and $\xi \in \mathfrak{g}$, regarded as the
velocity of system. If we introduce the conjugate momentum
$p_i=\frac{\partial L}{\partial \xi^i}$, $i=1,\cdots,n,\; n=dim G$,
and by the Legendre transformation $FL: TG \cong
G\times\mathfrak{g}\to T^\ast G \cong G\times\mathfrak{g}^\ast$,
$(g^i,\xi^i)\to (g^i,p_i)$, we have the Hamiltonian $H(g,p):T^\ast G
\cong G\times\mathfrak{g}^\ast \to \mathbb{R}$ given by
\begin{equation}H(g^i,p_i)=\sum_{i=1}^{n}p_i\xi^i-L(g^i,\xi^i).\label{5.3}
\end{equation}
If the Hamiltonian $H(g,p):T^\ast G\cong G\times\mathfrak{g} \to
\mathbb{R}$ is left cotangent lifted $G$-action invariant, for
$\mu\in\mathfrak{g}^\ast$ we have the associated reduced Hamiltonian
$h_\mu: (T^\ast G)_\mu \cong \mathcal{O}_\mu\to \mathbb{R}$, defined
by $h_\mu\cdot \pi_\mu=H\cdot i_\mu$. By the $(\pm)$-Lie-Poisson
brackets on $\mathfrak{g}^\ast$ and the symplectic structure on the
coadjoint orbit $\mathcal{O}_\mu$, we have the reduced Hamiltonian
vector field $X_{h_\mu}$ given by
\begin{equation}   X_{h_\mu}(\nu)=\mp
\operatorname{ad}^\ast_{\delta h_\mu/\delta \nu} \nu,\quad \forall
\nu\in \mathcal{O}_\mu.\label{5.4}
\end{equation}
See Marsden and Ratiu \cite{mara99}. Thus, if the 4-tuple $(T^\ast
G,G,\omega,H )$ is a regular point reducible Hamiltonian system on
the Lie group $G$, where the Hamiltonian $H: T^\ast G\to \mathbb{R}$
is invariant of the left cotangent lifted $G$-action, for a point
$\mu\in \mathfrak{g}^\ast$, the regular value of the momentum map
$\mathbf{J}_G: T^\ast G \rightarrow \mathfrak{g}^\ast$, then the
Marsden-Weinstein reduced Hamiltonian system is $(\mathcal{O}_\mu,
\omega_{\mathcal{O}_\mu}^{-}, h_\mu ). $ Moreover, assume that
$\gamma: G \rightarrow T^*G$ is an one-form on $G$, and
$\lambda=\gamma \cdot \pi_{G}: T^* G \rightarrow T^* G $, and
$\varepsilon: T^* G \rightarrow T^* G $ is a $G_\mu$-invariant
symplectic map, such that
$\varepsilon(\mathbf{J}_G^{-1}(\mu))\subset \mathbf{J}_G^{-1}(\mu),
$ and $\textmd{Im}(\gamma)\subset \mathbf{J}_G^{-1}(\mu), $ and it
is $G_\mu$-invariant, where $G_\mu= \{g\in
G|\operatorname{Ad}_g^\ast \mu=\mu \}$ is the isotropy subgroup of
the coadjoint $G$-action at the point $\mu\in\mathfrak{g}^\ast$.
Denote by $\bar{\gamma}=\pi_\mu(\gamma): G \rightarrow
\mathcal{O}_\mu, $ and $\bar{\lambda}=\pi_\mu(\lambda):
\mathbf{J}_G^{-1}(\mu) \rightarrow \mathcal{O}_\mu, $ and
$\bar{\varepsilon}=\pi_\mu(\varepsilon):
\mathbf{J}_G^{-1}(\mu)\rightarrow \mathcal{O}_\mu, $ where
$\pi_\mu:\mathbf{J}_G^{-1}(\mu)\rightarrow (T^\ast G)_\mu$ is the
projection. By using the similar ways in the proofs of two types of
Hamilton-Jacobi theorem for the Marsden-Weinstein reduced
Hamiltonian system, see Theorem 3.3 and Theorem 3.4, we can prove
the following theorem.

\begin{theorem}
For the regular point reducible Hamiltonian system
$(T^*G,G,\omega,H)$ on the Lie group $G$, assume that $\gamma: G
\rightarrow T^*G$ is an one-form on $G$, and $\lambda=\gamma
\cdot \pi_{G}: T^* G \rightarrow T^* G $, and $\varepsilon: T^* G \rightarrow T^* G $
is a symplectic map. Denote by
$X_H^\gamma = T\pi_{G}\cdot X_H \cdot \gamma$, and
$X_H^\varepsilon = T\pi_{G}\cdot X_H \cdot \varepsilon$, where $X_{H}$ is
the dynamical vector field of $(T^*G,G,\omega,H)$. Moreover, assume that
$\mu \in \mathfrak{g}^\ast $ is a regular value of the momentum
map $\mathbf{J}_G$, and $\textmd{Im}(\gamma)\subset
\mathbf{J}_G^{-1}(\mu), $ and it is $G_\mu$-invariant,
and $\varepsilon$ is $G_\mu$-invariant and
$\varepsilon(\mathbf{J}_G^{-1}(\mu))\subset \mathbf{J}_G^{-1}(\mu). $
Denote by $\bar{\gamma}=\pi_\mu(\gamma): G \rightarrow \mathcal{O}_\mu, $ and
$\bar{\lambda}=\pi_\mu(\lambda): \mathbf{J}_G^{-1}(\mu) \rightarrow \mathcal{O}_\mu, $
and
$\bar{\varepsilon}=\pi_\mu(\varepsilon): \mathbf{J}_G^{-1}(\mu)\rightarrow \mathcal{O}_\mu. $
Then the following two assertions hold:\\
\noindent $(\mathrm{i})$
If the one-form $\gamma: G \rightarrow T^*G $ is closed with respect to
$T\pi_G: TT^* G \rightarrow TG, $
then $\bar{\gamma}$ is a solution of the Type I of Hamilton-Jacobi equation
$T\bar{\gamma}\cdot X_H^\gamma= X_{h_\mu}\cdot \bar{\gamma}; $\\
\noindent $(\mathrm{ii})$
The $\varepsilon$ and $\bar{\varepsilon} $ satisfy the Type II of Hamilton-Jacobi equation
$T\bar{\gamma}\cdot X_H^\varepsilon= X_{h_\mu}\cdot \bar{\varepsilon}, $
if and only if they satisfy
the equation $T\bar{\varepsilon}\cdot(X_{h_\mu \cdot \bar{\varepsilon}})
= T\bar{\lambda}\cdot X_H \cdot\varepsilon. $\\
Here $X_{h_\mu}$ is the dynamical vector field of the Marsden-Weinstein
reduced Hamiltonian system $(\mathcal{O}_\mu,
\omega_{\mathcal{O}_\mu}^{-}, h_\mu ), $ and $X_{h_\mu \cdot \bar{\varepsilon}}$ is the Hamiltonian
vector field of the function $h_\mu \cdot \bar{\varepsilon}: T^* Q
\rightarrow \mathbb{R}. $ \hskip 0.3cm $\blacksquare$
\end{theorem}

Note that the symplectic form on the coadjoint orbit
$\mathcal{O}_\mu$ is induced by the (-)-Lie-Poisson brackets on
$\mathfrak{g}^\ast$, then the Type I and Type II of Hamilton-Jacobi
equation, $T\bar{\gamma}\cdot X_H^\gamma= X_{h_\mu}\cdot
\bar{\gamma}, $ and $T\bar{\gamma}\cdot X_H^\varepsilon=
X_{h_\mu}\cdot \bar{\varepsilon}, $ for the Marsden-Weinstein
reduced Hamiltonian system $(\mathcal{O}_\mu,
\omega_{\mathcal{O}_\mu}^{-}, h_\mu ), $ are also called the Type I
and Type II of Lie-Poisson Hamilton-Jacobi equation, respectively.
See Marsden and Ratiu \cite{mara99}, and Ge and Marsden
\cite{gema88}.

\subsection{Hamilton-Jacobi Equations of Rigid Body }

In the following we regard the rigid body as a regular point
reducible Hamiltonian system on the rotation group $\textmd{SO}(3)$,
and give its two types of Lie-Poisson Hamilton-Jacobi equation. Note that our
description of the motion and the equations of rigid body in this
subsection follows some of the notations and conventions in
Marsden and Ratiu \cite{mara99}, Marsden \cite{ma92}.\\

It is well known that, usually, the configuration space for a
$3$-dimensional rigid body moving freely in space is
$\textmd{SE}(3)$, the six dimensional group of Euclidean (rigid)
transformations of three dimensional space $\mathbb{R}^3$, that is,
all possible rotations and translations. If translations are ignored
and only rotations are considered, then the configuration space $Q$
is $\textmd{SO}(3)$, consists of all orthogonal linear
transformations of Euclidean 3-space to itself, which have
determinant one. Its Lie algebra, denoted $\mathfrak{so}(3)$,
consists of all $3\times 3$ skew matrices. By using the isomorphism
$\hat{}:\mathbb{R}^3\to \mathfrak{so}(3)$ defined by
$$(\Omega_1,\Omega_2,\Omega_3)=\Omega\to \hat{\Omega}=\begin{bmatrix}
  0&-\Omega_3&\Omega_2\\
  \Omega_3&0&-\Omega_1\\
  -\Omega_2&\Omega_1&0
\end{bmatrix},$$
we can identify the Lie algebra $(\mathfrak{so}(3), [,])$ with
$(\mathbb{R}^3, \times )$ and the Lie algebra bracket $[,]$ on
$\mathfrak{so}(3)$ with the cross product $\times$ of vectors in
$\mathbb{R}^3$. Denote by $\mathfrak{so}^\ast(3)$ the dual of the
Lie algebra $\mathfrak{so}(3)$, and we also identity
$\mathfrak{so}^\ast(3)$ with $\mathbb{R}^3$ by pairing the Euclidean
inner product. Since the functional derivative of a function defined
on $\mathbb{R}^3$ is equal to the usual gradient of the function,
from $(5.1)$ we know that the Lie-Poisson bracket on
$\mathfrak{so}^\ast(3)$ takes the form
\begin{equation}   \{f,g\}_\pm(\Pi)=\pm\Pi\cdot (\nabla_\Pi
f\times \nabla_\Pi g ), \;\; \forall f,g\in
C^\infty(\mathfrak{so}^\ast(3)),\;\; \Pi \in
\mathfrak{so}^\ast(3).\label{5.5}
\end{equation}

The phase space of a rigid body is the cotangent bundle $T^\ast G
=T^\ast \textmd{SO}(3)\cong \textmd{SO}(3)\times
\mathfrak{so}^\ast(3)$ (locally), with the canonical symplectic form,
by the local left trivialization. We consider Lie group $G=\textmd{SO}(3), $
which acts freely and properly by the
left translation on $\textmd{SO}(3)$ itself, then the action of
$\textmd{SO}(3)$ on the phase space $T^\ast \textmd{SO}(3)$ is given by the
cotangent lift of the left translation at the identity, that is, $\Phi:
\textmd{SO}(3)\times T^\ast \textmd{SO}(3)\cong \textmd{SO}(3)\times
\textmd{SO}(3)\times \mathfrak{so}^\ast(3)\to \textmd{SO}(3)\times
\mathfrak{so}^\ast(3),$ given by $\Phi(B,(A,\Pi))=(BA,\Pi)$, for any
$A,B\in \textmd{SO}(3), \; \Pi \in \mathfrak{so}^\ast(3)$, which is
also free, proper and symplectic. Assume that the left
$\textmd{SO}(3)$ action admits an associated
$\operatorname{Ad}^\ast$-equivariant momentum map $\mathbf{J}:T^\ast
\textmd{SO}(3)\to \mathfrak{so}^\ast(3)$, and if $\Pi \in \mathfrak{so}^\ast(3)$ is a
regular value of $\mathbf{J}$, then the regular point reduced space
$(T^\ast
\textmd{SO}(3))_\Pi=\mathbf{J}^{-1}(\Pi)/\textmd{SO}(3)_\Pi$ is
symplectically diffeomorphic to the coadjoint orbit $\mathcal{O}_\Pi
\subset \mathfrak{so}^\ast(3)$.\\

Let $I$ be the moment of inertia tensor computed with respect to a
body fixed frame, which is a principal body frame, we may represent
it by the diagonal matrix diag $(I_1,I_2,I_3)$. Let
$\Omega=(\Omega_1,\Omega_2,\Omega_3)$ be the vector of angular
velocities computed with respect to the axes fixed in the body and
$(\Omega_1,\Omega_2,\Omega_3)\in \mathfrak{so}(3)$. Consider the
Lagrangian $L(A,\Omega):\textmd{TSO}(3)\cong
\textmd{SO}(3)\times\mathfrak{so}(3)\to \mathbb{R}$, which is the
total kinetic energy of the rigid body, given by
$$L(A,\Omega)=\dfrac{1}{2}\langle\Omega,\Omega\rangle
=\dfrac{1}{2}(I_1\Omega_1^2+I_2\Omega_2^2+I_3\Omega_3^2),$$ where
$A\in \textmd{SO}(3)$, $(\Omega_1,\Omega_2,\Omega_3)\in
\mathfrak{so}(3)$. If we introduce the conjugate angular momentum
$\Pi_i=\dfrac{\partial L}{\partial \Omega_i}=I_i\Omega_i$,
$i=1,2,3$, which is also computed with respect to a body fixed
frame, and by the Legendre transformation $FL:\textmd{TSO}(3)\cong
\textmd{SO}(3)\times\mathfrak{so}(3)\to T^\ast \textmd{SO}(3)\cong
\textmd{SO}(3)\times \mathfrak{so}^\ast(3), \;
(A,\Omega)\to(A,\Pi)$, where $\Pi=(\Pi_1,\Pi_2,\Pi_3)\in
\mathfrak{so}^\ast(3)$, we have the Hamiltonian $H(A,\Pi):T^\ast
\textmd{SO}(3)\cong \textmd{SO}(3)\times \mathfrak{so}^\ast(3)\to
\mathbb{R}$ given by
\begin{align*} H(A,\Pi)&=\Omega \cdot \Pi-L(A,\Omega)\\
&=I_1\Omega_1^2+I_2\Omega_2^2+I_3\Omega_3^2-\frac{1}{2}(I_1\Omega_1^2+I_2\Omega_2^2+I_3\Omega_3^2)
=\frac{1}{2}(\frac{\Pi_1^2}{I_1}+\frac{\Pi_2^2}{I_2}+\frac{\Pi_3^2}{I_3}).
\end{align*}
From the Lie-Poisson bracket of rigid body on $\mathfrak{so}^\ast(3)$, that is, for $F,K:
\mathfrak{so}^\ast(3)\to \mathbb{R}, $ we have that
$\{F,K\}_{-}(\Pi)=-\Pi\cdot(\nabla_\Pi F\times \nabla_\Pi K), $ and hence
the Hamiltonian vector field of rigid body system is given by
\begin{align*}
X_{H}(\Pi)& =\{\Pi,\; H\}_{-}=
-\Pi\cdot(\nabla_\Pi\Pi\times\nabla_\Pi
H)= -\nabla_\Pi\Pi\cdot(\nabla_\Pi H\times \Pi)\\
& =(\Pi_1,\Pi_2,\Pi_3)\times (\frac{\Pi_1}{ I_1}, \frac{\Pi_2}{
I_2}, \frac{\Pi_3}{I_3})\\ &= (
\frac{(I_2-I_3)\Pi_2\Pi_3}{I_2I_3}, \;\;
\frac{(I_3-I_1)\Pi_3\Pi_1}{I_3I_1}, \;\;
\frac{(I_1-I_2)\Pi_1\Pi_2}{I_1I_2} ),
\end{align*}
since $\nabla_\Pi\Pi=1,$ and $\nabla_{\Pi_j} H= \Pi_j/I_j , \; j= 1,2,3 $.\\

From the above expression of the Hamiltonian, we know that
$H(A,\Pi)$ is invariant under the cotangent lift of the left
$\textmd{SO}(3)$-action, $\Phi: \textmd{SO}(3)\times T^\ast
\textmd{SO}(3) \to T^\ast \textmd{SO}(3)$. For the case $\Pi_0 =\mu
\in \mathfrak{so}^\ast(3)$ is a regular value of $\mathbf{J}$, we
have the reduced Hamiltonian
$h_{\mathcal{O}_\mu}(\Pi):\mathcal{O}_\mu\subset
\mathfrak{so}^\ast(3)\to \mathbb{R}$ given by
$h_{\mathcal{O}_\mu}(\Pi)\cdot
\pi_{\mathcal{O}_\mu}=H(A,\Pi)|_{\mathcal{O}_\mu}$. Moreover, note
that for $F_{\mathcal{O}_\mu},K_{\mathcal{O}_\mu}: \mathcal{O}_\mu
\to \mathbb{R}$, we have that
$\omega_{\mathcal{O}_\mu}^{-}(X_{F_{\mathcal{O}_\mu}},
X_{K_{\mathcal{O}_\mu}})=
\{F_{\mathcal{O}_\mu},K_{\mathcal{O}_\mu}\}_{-}|_{\mathcal{O}_\mu}$.
Thus, for the reduced Hamiltonian $h_{\mathcal{O}_\mu}(\Pi):
\mathcal{O}_\mu \to \mathbb{R}$, we have the Hamiltonian vector
field $X_{h_{\mathcal{O}_\mu}}(K_{\mathcal{O}_\mu})
=\{K_{\mathcal{O}_\mu},h_{\mathcal{O}_\mu}\}_{-}|_{\mathcal{O}_\mu
}. $ \\

In the following we shall derive the Type I and Type II of
Lie-Poisson Hamilton-Jacobi equation for the
regular point reduced rigid body system
$(\mathcal{O}_\mu,\omega_{\mathcal{O}_\mu}^{-},h_{\mathcal{O}_\mu}).$
Assume that $\gamma: \textmd{SO}(3) \rightarrow T^*
\textmd{SO}(3)$ is an one-form on $\textmd{SO}(3)$,
$\gamma(A)=(\gamma_1,\gamma_2,\gamma_3,\gamma_4,\gamma_5,\gamma_6),$
and
$\lambda=\gamma \cdot \pi_{\textmd{SO}(3)}: T^* \textmd{SO}(3)
\rightarrow T^* \textmd{SO}(3), $
$\lambda(A,\Pi)=(\lambda_1,\lambda_2,\lambda_3,\lambda_4,\lambda_5,\lambda_6),$ and
$\lambda_i(A,\Pi)=\gamma_i(A)\cdot \pi_{\textmd{SO}(3)}, \; i=1, \cdots, 6.$
For the regular value of $\mathbf{J}$, $\mu \in \mathfrak{so}^\ast(3),$
$\textmd{Im}(\gamma)\subset \mathbf{J}^{-1}(\mu), $ and it is
$\textmd{SO}(3)_\mu$-invariant, and $\bar{\gamma}=\pi_\mu(\gamma):
\textmd{SO}(3) \rightarrow \mathcal{O}_\mu, $
$\bar{\gamma}(A)=
(\bar{\gamma}_1,\bar{\gamma}_2,\bar{\gamma}_3) \in
\mathcal{O}_\mu (\subset \mathfrak{so}^\ast(3)), $ where
$\pi_\mu: \mathbf{J}^{-1}(\mu)\rightarrow \mathcal{O}_\mu. $ We choose that $\Pi=
(\Pi_1,\Pi_2,\Pi_3)=
(\bar{\gamma}_1,\bar{\gamma}_2,\bar{\gamma}_3) $,
then
$h_{\mathcal{O}_\mu} \cdot \bar{\gamma}: \textmd{SO}(3)
\rightarrow \mathbb{R} $ is given by
$$ h_{\mathcal{O}_\mu} \cdot \bar{\gamma}(A)=
H(A,\Pi)|_{\mathcal{O}_\mu} \cdot \bar{\gamma}(A)
=\frac{1}{2}(\frac{\bar{\gamma}_1^2}{I_1}+\frac{\bar{\gamma}_2^2}{I_2}+\frac{\bar{\gamma}_3^2}{I_3})
, $$ and the vector field
\begin{align*}
& X_{h_{\mathcal{O}_\mu}}(\Pi) \cdot \bar{\gamma}(A)
=\{\Pi,h_{\mathcal{O}_\mu} \}_{-}|_{\mathcal{O}_\mu} \cdot
\bar{\gamma}(A)\\
& = -\Pi\cdot(\nabla_\Pi\Pi\times\nabla_\Pi
(h_{\mathcal{O}_\mu}))\cdot \bar{\gamma} =
-\nabla_\Pi\Pi\cdot(\nabla_\Pi (h_{\mathcal{O}_\mu})\times
\Pi)\cdot \bar{\gamma}\\
& =(\Pi_1,\Pi_2,\Pi_3)\times (\frac{\Pi_1}{ I_1}, \frac{\Pi_2}{
I_2}, \frac{\Pi_3}{I_3})\cdot \bar{\gamma}\\ &= (
\frac{(I_2-I_3)\bar{\gamma}_2\bar{\gamma}_3}{I_2I_3}, \;\;
\frac{(I_3-I_1)\bar{\gamma}_3\bar{\gamma}_1}{I_3I_1}, \;\;
\frac{(I_1-I_2)\bar{\gamma}_1\bar{\gamma}_2}{I_1I_2} ),
\end{align*}
On the other hand, from the expression of
Hamiltonian vector field $X_H$, we have that
\begin{align*}
X_H^\gamma & =T\pi_{\textmd{SO}(3)}\cdot X_H\cdot\gamma= X_H\cdot\gamma \\
& =(
\frac{(I_2-I_3)\gamma_5\gamma_6}{I_2I_3}, \;\;
\frac{(I_3-I_1)\gamma_6\gamma_4}{I_3I_1}, \;\;
\frac{(I_1-I_2)\gamma_4\gamma_5}{I_1I_2} ).
\end{align*}
If $\gamma$ is closed with respect to
$T\pi_{\textmd{SO}(3)}: TT^* \textmd{SO}(3) \rightarrow T\textmd{SO}(3), $
then $\pi_{\textmd{SO}(3)}^*(\mathbf{d}\gamma)=0.$ We choose that
$(\gamma_4,\gamma_5,\gamma_6)=\Pi=(\Pi_1,\Pi_2,\Pi_3)=
(\bar{\gamma}_1,\bar{\gamma}_2,\bar{\gamma}_3), $ and hence
\begin{align*}
T\bar{\gamma}\cdot X_H^\gamma=(
\frac{(I_2-I_3)\bar{\gamma}_2\bar{\gamma}_3}{I_2I_3}, \;\;
\frac{(I_3-I_1)\bar{\gamma}_3\bar{\gamma}_1}{I_3I_1}, \;\;
\frac{(I_1-I_2)\bar{\gamma}_1\bar{\gamma}_2}{I_1I_2} )
= X_{h_{\mathcal{O}_\mu}} \cdot \bar{\gamma}.
\end{align*}
Thus, the Type I of Lie-Poisson Hamilton-Jacobi equation for the
regular point reduced rigid body system
$(\mathcal{O}_\mu,\omega_{\mathcal{O}_\mu}^{-},h_{\mathcal{O}_\mu})$ holds.\\

Next, for any $\textmd{SO}(3)_\mu$-invariant symplectic map $\varepsilon: T^* \textmd{SO}(3)
\rightarrow T^* \textmd{SO}(3),$ assume that $\varepsilon(A,\Pi)
=(\varepsilon_1,\varepsilon_2,\varepsilon_3,\varepsilon_4,\varepsilon_5,\varepsilon_6),$
and $\varepsilon(\mathbf{J}^{-1}(\mu))\subset \mathbf{J}^{-1}(\mu). $
Denote by $\bar{\varepsilon}=\pi_\mu(\varepsilon): \mathbf{J}^{-1}(\mu)\rightarrow \mathcal{O}_\mu, $
$\bar{\varepsilon}(A,\Pi)=(\bar{\varepsilon}_1,\bar{\varepsilon}_2,\bar{\varepsilon}_3) \in
\mathcal{O}_\mu, $ and
$\bar{\lambda}=\pi_\mu(\lambda): \mathbf{J}^{-1}(\mu) \rightarrow
\mathcal{O}_\mu, $ and
$\bar{\lambda}(A,\Pi)=
(\bar{\lambda}_1,\bar{\lambda}_2,\bar{\lambda}_3) \in
\mathcal{O}_\mu. $ We choose that $\Pi=
(\Pi_1,\Pi_2,\Pi_3)=
(\bar{\varepsilon}_1,\bar{\varepsilon}_2,\bar{\varepsilon}_3) $,
then
$h_{\mathcal{O}_\mu} \cdot \bar{\varepsilon}: T^*\textmd{SO}(3)
\rightarrow \mathbb{R} $ is given by
$$ h_{\mathcal{O}_\mu} \cdot \bar{\varepsilon}(A,\Pi)=
H(A,\Pi)|_{\mathcal{O}_\mu} \cdot \bar{\varepsilon}(A,\Pi)
=\frac{1}{2}(\frac{\bar{\varepsilon}_1^2}{I_1}+\frac{\bar{\varepsilon}_2^2}{I_2}+\frac{\bar{\varepsilon}_3^2}{I_3})
, $$ and the vector field
\begin{align*}
& X_{h_{\mathcal{O}_\mu}}(\Pi) \cdot \bar{\varepsilon}
=\{\Pi,h_{\mathcal{O}_\mu} \}_{-}|_{\mathcal{O}_\mu} \cdot
\bar{\varepsilon}(A,\Pi)
= -\Pi\cdot(\nabla_\Pi\Pi\times\nabla_\Pi
(h_{\mathcal{O}_\mu}))\cdot \bar{\varepsilon}\\
&= (
\frac{(I_2-I_3)\bar{\varepsilon}_2\bar{\varepsilon}_3}{I_2I_3}, \;\;
\frac{(I_3-I_1)\bar{\varepsilon}_3\bar{\varepsilon}_1}{I_3I_1}, \;\;
\frac{(I_1-I_2)\bar{\varepsilon}_1\bar{\varepsilon}_2}{I_1I_2} ),
\end{align*}
On the other hand, note that
\begin{align*}
X_H^\varepsilon & =T\pi_{\textmd{SO}(3)}\cdot X_H\cdot\varepsilon= X_H\cdot\varepsilon \\
& =(
\frac{(I_2-I_3)\varepsilon_5\varepsilon_6}{I_2I_3}, \;\;
\frac{(I_3-I_1)\varepsilon_6\varepsilon_4}{I_3I_1}, \;\;
\frac{(I_1-I_2)\varepsilon_4\varepsilon_5}{I_1I_2} ),
\end{align*}
and
\begin{align*}
T\bar{\gamma}\cdot X_H^\varepsilon=(
\frac{(I_2-I_3)\bar{\gamma}_2\bar{\gamma}_3}{I_2I_3}, \;\;
\frac{(I_3-I_1)\bar{\gamma}_3\bar{\gamma}_1}{I_3I_1}, \;\;
\frac{(I_1-I_2)\bar{\gamma}_1\bar{\gamma}_2}{I_1I_2} ),
\end{align*}
and
\begin{align*}
T\bar{\lambda}\cdot X_H \cdot \varepsilon=(
\frac{(I_2-I_3)\bar{\lambda}_2\bar{\lambda}_3}{I_2I_3}, \;\;
\frac{(I_3-I_1)\bar{\lambda}_3\bar{\lambda}_1}{I_3I_1}, \;\;
\frac{(I_1-I_2)\bar{\lambda}_1\bar{\lambda}_2}{I_1I_2} ).
\end{align*}
Thus, when we choose that $\Pi=(\Pi_1,\Pi_2,\Pi_3)\in
\mathcal{O}_\mu, $ and
$\Pi=(\bar{\gamma}_1,\bar{\gamma}_2,\bar{\gamma}_3)=
(\bar{\varepsilon}_1,\bar{\varepsilon}_2,\bar{\varepsilon}_3)=
(\bar{\lambda}_1,\bar{\lambda}_2,\bar{\lambda}_3), $ we must have
that
\begin{align*}
T\bar{\gamma}\cdot X_H^\varepsilon=X_{h_{\mathcal{O}_\mu}} \cdot \bar{\varepsilon}
=T\bar{\lambda}\cdot X_H \cdot \varepsilon.
\end{align*}
Since the map $\varepsilon: T^* \textmd{SO}(3)
\rightarrow T^* \textmd{SO}(3)$ is symplectic, then
$T\bar{\varepsilon}\cdot X_{h_{\mathcal{O}_\mu} \cdot \bar{\varepsilon}}
=X_{h_{\mathcal{O}_\mu}} \cdot \bar{\varepsilon}. $
Thus, in this case, we must have that
$\varepsilon$ and $\bar{\varepsilon} $ are the solution of the Type II of Lie-Poisson
Hamilton-Jacobi equation
$T\bar{\gamma}\cdot X_H^\varepsilon= X_{h_{\mathcal{O}_\mu}}\cdot \bar{\varepsilon}, $
for the regular point reduced rigid body system
$(\mathcal{O}_\mu,\omega_{\mathcal{O}_\mu}^{-},h_{\mathcal{O}_\mu})$, if and only if they satisfy
the equation $T\bar{\varepsilon}\cdot(X_{h_{\mathcal{O}_\mu} \cdot \bar{\varepsilon}})
= T\bar{\lambda}\cdot X_H \cdot\varepsilon. $\\

\begin{remark}
It is worthy of note that, if the one-form $\gamma: \textmd{SO}(3) \rightarrow T^*
\textmd{SO}(3)$ is determined by a generating function of a symplectic map, then it is a
solution of classical Hamilton-Jacobi equation $X_H\cdot \gamma=0.$
From Remark 3.7 we know that
the classical Lie-Poisson Hamilton-Jacobi equation for the
regular point reduced rigid body system
$(\mathcal{O}_\mu,\omega_{\mathcal{O}_\mu}^{-},h_{\mathcal{O}_\mu})$ is given by
$X_{h_{\mathcal{O}_\mu}}\cdot \bar{\gamma}=0, $ that is,
\begin{equation*}
 \left\{\begin{aligned}
 & \frac{(I_2-I_3)\bar{\gamma}_2\bar{\gamma}_3}{I_2I_3}=0, \\
 & \frac{(I_3-I_1)\bar{\gamma}_3\bar{\gamma}_1}{I_3I_1}=0, \\
 & \frac{(I_1-I_2)\bar{\gamma}_1\bar{\gamma}_2}{I_1I_2}=0,
\end{aligned} \right. \label{5.6}
\end{equation*}
where $\bar{\gamma}=\pi_\mu(\gamma):
\textmd{SO}(3) \rightarrow \mathcal{O}_\mu. $
\end{remark}

To sum up the above discussion, we have the following proposition.
\begin{proposition}
If the 4-tuple $(T^\ast\textmd{SO}(3), \textmd{SO}(3),\omega, H)$ is a
regular point reducible Hamiltonian system, then for a point $\mu \in
\mathfrak{so}^\ast(3)$, the regular value of the momentum map
$\mathbf{J}:T^\ast \textmd{SO}(3)\to \mathfrak{so}^\ast(3)$, the
Marsden-Weinstein reduced system is the 3-tuple
$(\mathcal{O}_\mu,\omega_{\mathcal{O}_\mu}^{-},h_{\mathcal{O}_\mu}
),$ where $\mathcal{O}_\mu \subset \mathfrak{so}^\ast(3)$ is the
coadjoint orbit, $\omega_{\mathcal{O}_\mu}^{-}$ is the orbit symplectic
form on $\mathcal{O}_\mu $, which is induced by the rigid body
Poisson bracket on $\mathfrak{so}^\ast(3)$,
$h_{\mathcal{O}_\mu}(\Pi)\cdot\pi_{\mathcal{O}_\mu}=H(A,\Pi)
=H(A,\Pi)|_{\mathcal{O}_\mu}$. Assume that $\gamma: \textmd{SO}(3)
\rightarrow T^* \textmd{SO}(3)$ is an one-form on
$\textmd{SO}(3)$,
and $\lambda=\gamma \cdot \pi_{\textmd{SO}(3)}:
T^* \textmd{SO}(3) \rightarrow T^* \textmd{SO}(3), $ and $\varepsilon:
T^* \textmd{SO}(3) \rightarrow T^* \textmd{SO}(3) $ is a
$\textmd{SO}(3)_\mu$-invariant symplectic map.
Denote by
$X_H^\gamma = T\pi_{\textmd{SO}(3)}\cdot X_H \cdot \gamma$, and
$X_H^\varepsilon = T\pi_{\textmd{SO}(3)}\cdot X_H \cdot \varepsilon$.
Moreover, assume that $\textmd{Im}(\gamma)\subset \mathbf{J}^{-1}(\mu), $ and it is
$\textmd{SO}(3)_\mu$-invariant,
and $\varepsilon(\mathbf{J}^{-1}(\mu))\subset \mathbf{J}^{-1}(\mu). $
Denote by $\bar{\gamma}=\pi_\mu(\gamma):
\textmd{SO}(3) \rightarrow \mathcal{O}_\mu, $ and
$\bar{\lambda}=\pi_\mu(\lambda): \mathbf{J}^{-1}(\mu) \rightarrow
\mathcal{O}_\mu, $ and
$\bar{\varepsilon}=\pi_\mu(\varepsilon): \mathbf{J}^{-1}(\mu)\rightarrow \mathcal{O}_\mu. $
Then the following two assertions hold:\\
\noindent $(\mathrm{i})$
If the one-form $\gamma: \textmd{SO}(3) \rightarrow T^*\textmd{SO}(3) $ is closed with respect to
$T\pi_{\textmd{SO}(3)}: TT^* \textmd{SO}(3) \rightarrow T\textmd{SO}(3), $
then $\bar{\gamma}$ is a solution of the Type I of Lie-Poisson Hamilton-Jacobi equation
$T\bar{\gamma}\cdot X_H^\gamma= X_{h_{\mathcal{O}_\mu}}\cdot \bar{\gamma}; $\\
\noindent $(\mathrm{ii})$
The $\varepsilon$ and $\bar{\varepsilon} $ satisfy the Type II of Lie-Poisson Hamilton-Jacobi equation
$T\bar{\gamma}\cdot X_H^\varepsilon= X_{h_{\mathcal{O}_\mu}}\cdot \bar{\varepsilon}, $
if and only if they satisfy
the equation $T\bar{\varepsilon}\cdot(X_{h_{\mathcal{O}_\mu} \cdot \bar{\varepsilon}})
= T\bar{\lambda}\cdot X_H \cdot\varepsilon. $ \hskip 0.3cm $\blacksquare$
\end{proposition}

\subsection{Hamilton-Jacobi Equations of Heavy Top }

In the following we regard the heavy top as a regular point
reducible Hamiltonian system on the Euclidean group
$\textmd{SE}(3)$, and give its two types of Lie-Poisson Hamilton-Jacobi equation.
Note that our description of the motion and the equations of heavy
top follows some of the notations and conventions in
Marsden and Ratiu \cite{mara99}, Marsden \cite{ma92}.\\

We know that a heavy top is by definition a rigid body with a fixed
point in $\mathbb{R}^3$ and moving in gravitational field. Usually,
exception of the singular point, its physical phase space is $T^\ast
\textmd{SO}(3)$ and the symmetry group is $S^1$, regarded as
rotations about the z-axis, the axis of gravity, this is because
gravity breaks the symmetry and the system is no longer
$\textmd{SO}(3)$ invariant. By the semidirect product reduction
theorem, see Marsden et al. \cite{mamiorpera07} and \cite{mawazh10},
we show that the reduction of $T^\ast \textmd{SO}(3)$ by $S^1$ gives
a space which is symplectically diffeomorphic to the reduced space
obtained by the reduction of $T^\ast \textmd{SE}(3)$ by the left action
of $\textmd{SE}(3)$, that is, the coadjoint orbit
$\mathcal{O}_{(\mu,a)} \subset \mathfrak{se}^\ast(3)\cong T^\ast
\textmd{SE}(3)/\textmd{SE}(3)$. In fact, in this case, we can
identify the phase space $T^\ast \textmd{SO}(3)$ with the reduction
of the cotangent bundle of the special Euclidean group
$\textmd{SE}(3)=\textmd{SO}(3)\circledS \mathbb{R}^3$ by the
Euclidean translation subgroup $\mathbb{R}^3$ and identifies the
symmetry group $S^1$ with isotropy group $G_a=\{ A\in
\textmd{SO}(3)\mid Aa=a \}=S^1$, which is Abelian and
$(G_a)_{\mu_a}= G_a =S^1,\; \forall \mu_a \in \mathfrak{g}^\ast_a$,
where $a$ is a vector aligned with the direction of gravity and
where $\textmd{SO}(3)$ acts on $\mathbb{R}^3$ in the standard way.\\

Now we consider the cotangent bundle $T^\ast G =T^\ast
\textmd{SE}(3)\cong \textmd{SE}(3)\times \mathfrak{se}^\ast(3)$
(locally), with the canonical symplectic form, by the local left
trivialization. We consider Lie group $G=\textmd{SE}(3)=
\textmd{SO}(3)\circledS \mathbb{R}^3 $ with Lie algebra
$\mathfrak{se}(3)=\mathfrak{so}(3) \circledS \mathbb{R}^3$, which is
a semidirect product Lie group and acts freely and properly by the
left translation on $\textmd{SE}(3)$ itself, then the action of
$\textmd{SE}(3)$ on the phase space $T^\ast \textmd{SE}(3)$ is given
by cotangent lift of the left translation at the identity, that is,
$\Phi: \textmd{SE}(3)\times T^\ast \textmd{SE}(3) \cong
\textmd{SE}(3)\times \textmd{SE}(3)\times \mathfrak{se}^\ast(3)\to
\textmd{SE}(3)\times \mathfrak{se}^\ast(3),$ is given by
$\Phi((B,u),(A,v,\Pi,w))=(BA,v,\Pi,w)$, for any $A,B\in
\textmd{SO}(3), \; \Pi \in \mathfrak{so}^\ast(3), \; u,v,w \in
\mathbb{R}^3$, which is also free, proper and symplectic. Assume
that the left $\textmd{SE}(3)$ action admits an associated
$\operatorname{Ad}^\ast$-equivariant momentum map $\mathbf{J}:T^\ast
\textmd{SE}(3)\to \mathfrak{se}^\ast(3)$, and if $(\Pi,w) \in
\mathfrak{se}^\ast(3)$ is a regular value of $\mathbf{J}$, then the
regular point reduced space $(T^\ast
\textmd{SE}(3))_{(\Pi,w)}=\mathbf{J}^{-1}(\Pi,w)/\textmd{SE}(3)_{(\Pi,w)}$
is symplectically diffeomorphic to the coadjoint orbit
$\mathcal{O}_{(\Pi,w)} \subset \mathfrak{se}^\ast(3)$.\\

Let $I=diag(I_1,I_2,I_3)$ be the moment of inertia of the heavy top
in the body-fixed frame, that is, it is in principal body frame. Let
$\Omega=(\Omega_1,\Omega_2,\Omega_3)$ be the vector of heavy top
angular velocities computed with respect to the axes fixed in the
body and $(\Omega_1,\Omega_2,\Omega_3)\in \mathfrak{so}(3)$. Let
$\Gamma$ be the unit vector viewed by an observer moving with the
body, $m$ be that total mass of the system, $g$ be the magnitude of
the gravitational acceleration, $\chi$ be the unit vector on the
line connecting the origin $O$ to the center of mass of the system,
and $h$ be the length of this segment.\\

Note that Lie algebra $\mathfrak{se}(3)=\mathfrak{so}(3) \circledS
\mathbb{R}^3$ and its dual
$\mathfrak{se}^\ast(3)=\mathfrak{so}^\ast(3) \circledS
\mathbb{R}^3$, we consider the Lagrangian
$L(A,v,\Omega,\Gamma):\textmd{TSE}(3)\cong
\textmd{SE}(3)\times\mathfrak{se}(3)\to \mathbb{R}$ , which is the
total kinetic minus potential energy of the heavy top, given by
$$L(A,v,\Omega,\Gamma)=\dfrac{1}{2}\langle\Omega,\Omega\rangle-mgh\Gamma\cdot\chi
=\dfrac{1}{2}(I_1\Omega_1^2+I_2\Omega_2^2+I_3\Omega_3^2)-mgh\Gamma\cdot\chi,$$
where $(A,v)\in \textmd{SE}(3)$,
$\Omega=(\Omega_1,\Omega_2,\Omega_3)\in \mathfrak{so}(3)$,
$\Gamma\in\mathbb{R}^3$, and the variable $\Gamma$ is regarded as a
parameter with respect to potential energy of the heavy top. If we
introduce the conjugate angular momentum $\Pi_i=\dfrac{\partial
L}{\partial \Omega_i}=I_i\Omega_i, \; i=1,2,3,$ and by the Legendre
transformation with the parameter $\Gamma$, $FL:\textmd{TSE}(3)\cong
\textmd{SE}(3)\times\mathfrak{so}(3) \circledS \mathbb{R}^3\to
T^\ast \textmd{SE}(3)\cong \textmd{SE}(3)\times
\mathfrak{so}^\ast(3) \circledS \mathbb{R}^3,\quad
(A,v,\Omega,\Gamma)\to(A,v,\Pi,\Gamma)$, where
$\Pi=(\Pi_1,\Pi_2,\Pi_3)\in \mathfrak{so}^\ast(3)$, we have the
Hamiltonian $H(A,v,\Pi,\Gamma): T^\ast \textmd{SE}(3)\cong
\textmd{SE}(3)\times \mathfrak{so}^\ast(3) \circledS \mathbb{R}^3\to
\mathbb{R}$ given by
\begin{align*} H(A,v,\Pi,\Gamma)&=\Omega\cdot
\Pi-L(A,v,\Omega,,\Gamma)\\
&=I_1\Omega_1^2+I_2\Omega_2^2+I_3\Omega_3^2-
\frac{1}{2}(I_1\Omega_1^2+I_2\Omega_2^2+I_3\Omega_3^2)+mgh\Gamma\cdot\chi\\
&=\frac{1}{2}(\frac{\Pi_1^2}{I_1}+\frac{\Pi_2^2}{I_2}+\frac{\Pi_3^2}{I_3})+mgh\Gamma\cdot\chi.
\end{align*}
From the semidirect product Poisson bracket, see Marsden et al.
\cite{mamiorpera07}, we can get the heavy top Poisson bracket on
$\mathfrak{se}^\ast(3)$, that is, for $F,K: \mathfrak{se}^\ast(3)\to
\mathbb{R}, $ we have that
\begin{equation}
\{F,K\}_{-}(\Pi,\Gamma)=-\Pi\cdot(\nabla_\Pi F\times\nabla_\Pi
K)-\Gamma\cdot(\nabla_\Pi F\times \nabla_\Gamma K-\nabla_\Pi K\times
\nabla_\Gamma F). \label{4.6}
\end{equation}
Thus, the Hamiltonian vector fields of heavy top system are given by
\begin{align*}
 X_{H}(\Pi)& = \{\Pi,\; H \}_{-}= -\Pi\cdot(\nabla_\Pi\Pi\times\nabla_\Pi
H) -\Gamma\cdot(\nabla_\Pi\Pi\times\nabla_\Gamma
H-\nabla_\Pi H \times\nabla_\Gamma\Pi)\\
& =(\Pi_1,\Pi_2,\Pi_3)\times (\frac{\Pi_1}{ I_1}, \frac{\Pi_2}{
I_2}, \frac{\Pi_3}{I_3})
+mgh(\Gamma_1,\Gamma_2,\Gamma_3)\times (\chi_1,\chi_2,\chi_3)\\
& = ( \frac{(I_2-I_3)\Pi_2\Pi_3}{I_2I_3}+
mgh(\Gamma_2\chi_3-\Gamma_3\chi_2),\;\;
\frac{(I_3-I_1)\Pi_3\Pi_1}{I_3I_1} +
mgh(\Gamma_3\chi_1-\Gamma_1\chi_3),\\
& \;\;\;\;\;\;
\frac{(I_1-I_2)\Pi_1\Pi_2}{I_1I_2}
 + mgh(\Gamma_1\chi_2-\Gamma_2\chi_1) ),
\end{align*}
since $\nabla_\Pi\Pi=1, \; \nabla_\Gamma\Pi =0, \;
\chi=(\chi_1,\chi_2,\chi_3), $ and $\nabla_{\Pi_j}
H= \Pi_j/I_j , \; \nabla_{\Gamma_j}
H= mgh\chi_j, \; j= 1,2,3. $

\begin{align*}
 X_{H}(\Gamma)&
= \{\Gamma,\; H\}_{-} =-\Pi\cdot(\nabla_\Pi\Gamma\times\nabla_\Pi
H)-\Gamma\cdot(\nabla_\Pi\Gamma\times\nabla_\Gamma
H-\nabla_\Pi H \times\nabla_\Gamma\Gamma) \\
& =\nabla_\Gamma\Gamma\cdot(\Gamma\times\nabla_\Pi
H)=(\Gamma_1,\Gamma_2,\Gamma_3)\times (\frac{\Pi_1}{ I_1},
\frac{\Pi_2}{ I_2}, \frac{\Pi_3}{I_3})\\
& = ( \frac{I_2\Gamma_2\Pi_3-
I_3\Gamma_3\Pi_2}{I_2I_3}, \;\;
\frac{I_3\Gamma_3\Pi_1-
I_1\Gamma_1\Pi_3}{I_3I_1}, \;\;
\frac{I_1\Gamma_1\Pi_2-
I_2\Gamma_2\Pi_1}{I_1I_2} ),
\end{align*}
since $\nabla_\Gamma \Gamma =1, \; \nabla_\Pi\Gamma =0, $ and
$\nabla_{\Pi_j} H= \Pi_j/I_j , \; j= 1,2,3.$\\

From the above expression of the Hamiltonian, we know that
$H(A,v,\Pi,\Gamma)$ is invariant under the cotangent lift of the
left $\textmd{SE}(3)$-action, $\Phi:\textmd{SE}(3)\times T^\ast
\textmd{SE}(3)\to T^\ast \textmd{SE}(3)$. For the case
$(\Pi_0,\Gamma_0)=(\mu,a)\in \mathfrak{se}^\ast(3)$ is a regular
value of $\mathbf{J}$, we have the reduced Hamiltonian
$h_{\mathcal{O}_{(\mu,a)}}(\Pi,,\Gamma):\mathcal{O}_{(\mu,a)}(\subset
\mathfrak{se}^\ast (3))\to \mathbb{R}$ given by
$h_{\mathcal{O}_{(\mu,a)}}(\Pi,\Gamma)\cdot
\pi_{\mathcal{O}_{(\mu,a)}}
=H(A,v,\Pi,\Gamma)|_{\mathcal{O}_{(\mu,a)}}$. Moreover, note that
for $F_{\mathcal{O}_{(\mu,a)}},K_{\mathcal{O}_{(\mu,a)}}:
\mathcal{O}_{(\mu,a)} \to \mathbb{R}$, we have that
$$\omega_{\mathcal{O}_{(\mu,a)}}^{-}(X_{F_{\mathcal{O}_{(\mu,a)}}},
X_{K_{\mathcal{O}_{(\mu,a)}}})=
\{F_{\mathcal{O}_{(\mu,a)}},K_{\mathcal{O}_{(\mu,a)}}\}_{-}|_{\mathcal{O}_{(\mu,a)}}.
$$ Thus, for the reduced Hamiltonian $h_{\mathcal{O}_{(\mu,a)}}(\Pi,\Gamma):
\mathcal{O}_{(\mu,a)} \to \mathbb{R}$, we have the Hamiltonian
vector field
$X_{h_{\mathcal{O}_{(\mu,a)}}}(K_{\mathcal{O}_{(\mu,a)}})
=\{K_{\mathcal{O}_{(\mu,a)}},h_{\mathcal{O}_{(\mu,a)}}\}_{-}|_{\mathcal{O}_{(\mu,a)}
}.$ \\

In the following we shall derive the Type I and Type II of
Lie-Poisson Hamilton-Jacobi equation for the
regular point reduced heavy top system
$(\mathcal{O}_{(\mu,a)},\omega^{-}_{\mathcal{O}_{(\mu,a)}},h_{\mathcal{O}_{(\mu,a)}}).$
Assume that $\gamma: \textmd{SE}(3) \rightarrow T^*
\textmd{SE}(3)$ is an one-form on $\textmd{SE}(3)$,
$\gamma(A,v)=(\gamma_1,\cdots,\gamma_{12}),$
and $\lambda=\gamma \cdot \pi_{\textmd{SE}(3)}: T^* \textmd{SE}(3)
\rightarrow T^* \textmd{SE}(3), $
$\lambda(A,v,\Pi,\Gamma)=(\lambda_1,\cdots,\lambda_{12}),$ and
$\lambda_i(A,v,\Pi,\Gamma)=\gamma_i(A,v)\cdot \pi_{\textmd{SE}(3)}, \; i=1, \cdots, 12.$
For the regular value of $\mathbf{J}$, $(\mu,a) \in \mathfrak{se}^\ast(3),$
$\textmd{Im}(\gamma)\subset \mathbf{J}^{-1}((\mu,a)), $ and it is
$\textmd{SE}(3)_{(\mu,a)}$-invariant, and
$\bar{\gamma}=\pi_{(\mu,a)}(\gamma): \textmd{SE}(3) \rightarrow
\mathcal{O}_{(\mu,a)}, $
$\bar{\gamma}(A,v)=
(\bar{\gamma}_1,\bar{\gamma}_2,\bar{\gamma}_3,\bar{\gamma}_4,\bar{\gamma}_5,\bar{\gamma}_6) \in
\mathcal{O}_{(\mu,a)} (\subset \mathfrak{se}^\ast(3)), $ where
$\pi_{(\mu,a)}: \mathbf{J}^{-1}((\mu,a))\rightarrow \mathcal{O}_{(\mu,a)}. $ We choose that $\Pi=
(\Pi_1,\Pi_2,\Pi_3)=
(\bar{\gamma}_1,\bar{\gamma}_2,\bar{\gamma}_3), \; \Gamma= (\Gamma_1,\Gamma_2,\Gamma_3)=
(\bar{\gamma}_4,\bar{\gamma}_5,\bar{\gamma}_6), $
then
$h_{\mathcal{O}_{(\mu,a)}} \cdot \bar{\gamma}: \textmd{SE}(3)
\rightarrow \mathbb{R} $ is given by
\begin{align*}
& h_{\mathcal{O}_{(\mu,a)}} \cdot \bar{\gamma}(A,v,\Pi,\Gamma)=
H(A,v,\Pi,\Gamma)|_{\mathcal{O}_{(\mu,a)}} \cdot \bar{\gamma}(A,v,\Pi,\Gamma)\\
& =\frac{1}{2}(\frac{\bar{\gamma}_1^2}{I_1}
+\frac{\bar{\gamma}_2^2}{I_2}+\frac{\bar{\gamma}_3^2}{I_3}) +
mgh(\bar{\gamma}_4\cdot\chi_1+ \bar{\gamma}_5\cdot\chi_2+
\bar{\gamma}_6\cdot\chi_3),
\end{align*} and the vector field
\begin{align*}
& X_{h_{\mathcal{O}_{(\mu,a)}}}(\Pi) \cdot
\bar{\gamma}
= \{\Pi,h_{\mathcal{O}_{(\mu,a)}}\}_{-}|_{\mathcal{O}_{(\mu,a)}}\cdot \bar{\gamma}(A,v,\Pi,\Gamma)\\
& = -\Pi\cdot(\nabla_\Pi\Pi\times\nabla_\Pi
(h_{\mathcal{O}_{(\mu,a)}})) \cdot
\bar{\gamma}-\Gamma\cdot(\nabla_\Pi\Pi\times\nabla_\Gamma
(h_{\mathcal{O}_{(\mu,a)}})-\nabla_\Pi
(h_{\mathcal{O}_{(\mu,a)}}) \times\nabla_\Gamma\Pi)\cdot \bar{\gamma}\\
& =(\Pi_1,\Pi_2,\Pi_3)\times (\frac{\Pi_1}{ I_1}, \frac{\Pi_2}{
I_2}, \frac{\Pi_3}{I_3}) \cdot \bar{\gamma}
+mgh(\Gamma_1,\Gamma_2,\Gamma_3)\times (\chi_1,\chi_2,\chi_3)\cdot \bar{\gamma}\\
& = ( \frac{(I_2-I_3)\bar{\gamma}_2\bar{\gamma}_3}{I_2I_3}+
mgh(\bar{\gamma}_5\chi_3-\bar{\gamma}_6\chi_2),\;\;
\frac{(I_3-I_1)\bar{\gamma}_3\bar{\gamma}_1}{I_3I_1} +
mgh(\bar{\gamma}_6\chi_1-\bar{\gamma}_4\chi_3),\\
& \;\;\;\;\;\;
\frac{(I_1-I_2)\bar{\gamma}_1\bar{\gamma}_2}{I_1I_2}
 + mgh(\bar{\gamma}_4\chi_2-\bar{\gamma}_5\chi_1) ),
\end{align*}
and
\begin{align*}
& X_{h_{\mathcal{O}_{(\mu,a)}}}(\Gamma) \cdot
\bar{\gamma}
= \{\Gamma,h_{\mathcal{O}_{(\mu,a)}}\}_{-}|_{\mathcal{O}_{(\mu,a)}}\cdot \bar{\gamma}(A,v,\Pi,\Gamma)\\
& =-\Pi\cdot(\nabla_\Pi\Gamma\times\nabla_\Pi
(h_{\mathcal{O}_{(\mu,a)}})) \cdot \bar{\gamma}-
\Gamma\cdot(\nabla_\Pi\Gamma\times\nabla_\Gamma
(h_{\mathcal{O}_{(\mu,a)}})-\nabla_\Pi
(h_{\mathcal{O}_{(\mu,a)}}) \times\nabla_\Gamma\Gamma) \cdot \bar{\gamma}\\
& =\nabla_\Gamma\Gamma\cdot(\Gamma\times\nabla_\Pi
(h_{\mathcal{O}_{(\mu,a)}}))\cdot \bar{\gamma}
=(\Gamma_1,\Gamma_2,\Gamma_3)\times (\frac{\Pi_1}{ I_1},
\frac{\Pi_2}{ I_2}, \frac{\Pi_3}{I_3})\cdot
\bar{\gamma}\\
& = ( \frac{I_2\bar{\gamma}_5\bar{\gamma}_3-
I_3\bar{\gamma}_6\bar{\gamma}_2}{I_2I_3}, \;\;
\frac{I_3\bar{\gamma}_6\bar{\gamma}_1-
I_1\bar{\gamma}_4\bar{\gamma}_3}{I_3I_1}, \;\;
\frac{I_1\bar{\gamma}_4\bar{\gamma}_2-
I_2\bar{\gamma}_5\bar{\gamma}_1}{I_1I_2} ).
\end{align*}
On the other hand, from the expressions of
Hamiltonian vector fields $X_H(\Pi)$ and $X_H(\Gamma)$, we have that
\begin{align*}
X_H(\Pi)^\gamma & =T\pi_{\textmd{SE}(3)}\cdot X_H(\Pi)\cdot\gamma= X_H(\Pi)\cdot\gamma \\
& = ( \frac{(I_2-I_3)\gamma_8\gamma_9}{I_2I_3}+
mgh(\gamma_{11}\chi_3-\gamma_{12}\chi_2),\;\;
\frac{(I_3-I_1)\gamma_9\gamma_7}{I_3I_1} +
mgh(\gamma_{12}\chi_1-\gamma_{10}\chi_3),\\
& \;\;\;\;\;\;
\frac{(I_1-I_2)\gamma_7\gamma_8}{I_1I_2}
 + mgh(\gamma_{10}\chi_2-\gamma_{11}\chi_1) ),
\end{align*}
and
\begin{align*}
X_H(\Gamma)^\gamma & =T\pi_{\textmd{SE}(3)}\cdot X_H(\Gamma)\cdot\gamma= X_H(\Gamma)\cdot\gamma \\
& = ( \frac{I_2\gamma_{11}\gamma_9-
I_3\gamma_{12}\gamma_8}{I_2I_3}, \;\;
\frac{I_3\gamma_{12}\gamma_7-
I_1\gamma_{10}\gamma_9}{I_3I_1}, \;\;
\frac{I_1\gamma_{10}\gamma_8-
I_2\gamma_{11}\gamma_7}{I_1I_2} ).
\end{align*}
If $\gamma$ is closed with respect to
$T\pi_{\textmd{SE}(3)}: TT^* \textmd{SE}(3) \rightarrow T\textmd{SE}(3), $
then $\pi_{\textmd{SE}(3)}^*(\mathbf{d}\gamma)=0.$ We choose that
$(\gamma_7,\gamma_8,\gamma_9)=\Pi=(\Pi_1,\Pi_2,\Pi_3)=
(\bar{\gamma}_1,\bar{\gamma}_2,\bar{\gamma}_3), \;
(\gamma_{10},\gamma_{11},\gamma_{12})=\Gamma=(\Gamma_1,\Gamma_2,\Gamma_3)=
(\bar{\gamma}_4,\bar{\gamma}_5,\bar{\gamma}_6), $ and hence
\begin{align*}
& T\bar{\gamma}\cdot X_H(\Pi)^\gamma \\
& =( \frac{(I_2-I_3)\bar{\gamma}_2\bar{\gamma}_3}{I_2I_3}+
mgh(\bar{\gamma}_5\chi_3-\bar{\gamma}_6\chi_2),\;\;
\frac{(I_3-I_1)\bar{\gamma}_3\bar{\gamma}_1}{I_3I_1} +
mgh(\bar{\gamma}_6\chi_1-\bar{\gamma}_4\chi_3),\\
& \;\;\;\;\;\;
\frac{(I_1-I_2)\bar{\gamma}_1\bar{\gamma}_2}{I_1I_2}
 + mgh(\bar{\gamma}_4\chi_2-\bar{\gamma}_5\chi_1)) \\
& = X_{h_{\mathcal{O}_{(\mu,a)}}}(\Pi) \cdot \bar{\gamma},
\end{align*}
and
\begin{align*}
& T\bar{\gamma}\cdot X_H(\Gamma)^\gamma \\
& = ( \frac{I_2\bar{\gamma}_5\bar{\gamma}_3-
I_3\bar{\gamma}_6\bar{\gamma}_2}{I_2I_3}, \;\;
\frac{I_3\bar{\gamma}_6\bar{\gamma}_1-
I_1\bar{\gamma}_4\bar{\gamma}_3}{I_3I_1}, \;\;
\frac{I_1\bar{\gamma}_4\bar{\gamma}_2-
I_2\bar{\gamma}_5\bar{\gamma}_1}{I_1I_2} )\\
& = X_{h_{\mathcal{O}_{(\mu,a)}}}(\Gamma) \cdot \bar{\gamma}.
\end{align*}
Thus, $T\bar{\gamma}\cdot X_H^\gamma= X_{h_{\mathcal{O}_{(\mu,a)}}}\cdot \bar{\gamma}, $ that is,
the Type I of Lie-Poisson Hamilton-Jacobi equation for the
regular point reduced heavy top system
$(\mathcal{O}_{(\mu,a)},\omega^{-}_{\mathcal{O}_{(\mu,a)}},h_{\mathcal{O}_{(\mu,a)}})$ holds.\\

Next, for any $\textmd{SE}(3)_{(\mu,a)}$-invariant symplectic map
$\varepsilon: T^* \textmd{SE}(3) \rightarrow T^* \textmd{SE}(3),$
assume that $\varepsilon(A,v,\Pi,\Gamma) =(\varepsilon_1,\cdots,
\varepsilon_{12}),$ and
$\varepsilon(\mathbf{J}^{-1}((\mu,a)))\subset
\mathbf{J}^{-1}((\mu,a)). $ Denote by
$\bar{\varepsilon}=\pi_{(\mu,a)}(\varepsilon):
\mathbf{J}^{-1}((\mu,a)) \rightarrow \mathcal{O}_{(\mu,a)}, $
$\bar{\varepsilon}(A,v,\Pi,\Gamma)=(\bar{\varepsilon}_1,\cdots,\bar{\varepsilon}_6)
\in \mathcal{O}_{(\mu,a)}(\subset \mathfrak{se}^\ast(3)), $ and
$\bar{\lambda}=\pi_{(\mu,a)}(\lambda): T^* \textmd{SE}(3)
\rightarrow \mathcal{O}_{(\mu,a)}, $ and
$\bar{\lambda}(A,v,\Pi,\Gamma)=
(\bar{\lambda}_1,\cdots,\bar{\lambda}_6) \in \mathcal{O}_{(\mu,a)}
(\subset \mathfrak{se}^\ast(3)). $ We choose that
$\Pi=(\Pi_1,\Pi_2,\Pi_3)=(\bar{\varepsilon}_1,\bar{\varepsilon}_2,\bar{\varepsilon}_3),
\; \Gamma= (\Gamma_1,\Gamma_2,\Gamma_3)=
(\bar{\varepsilon}_4,\bar{\varepsilon}_5,\bar{\varepsilon}_6), $
then $h_{\mathcal{O}_{(\mu,a)}} \cdot \bar{\varepsilon}:
T^*\textmd{SE}(3) \rightarrow \mathbb{R} $ is given by
\begin{align*}
& h_{\mathcal{O}_{(\mu,a)}} \cdot \bar{\varepsilon}(A,v,\Pi,\Gamma)=
H(A,v,\Pi,\Gamma)|_{\mathcal{O}_{(\mu,a)}} \cdot \bar{\varepsilon}(A,v,\Pi,\Gamma)\\
& =\frac{1}{2}(\frac{\bar{\varepsilon}_1^2}{I_1}
+\frac{\bar{\varepsilon}_2^2}{I_2}+\frac{\bar{\varepsilon}_3^2}{I_3}) +
mgh(\bar{\varepsilon}_4\cdot\chi_1+ \bar{\varepsilon}_5\cdot\chi_2+
\bar{\varepsilon}_6\cdot\chi_3),
\end{align*} and the vector field
\begin{align*}
& X_{h_{\mathcal{O}_{(\mu,a)}}}(\Pi) \cdot
\bar{\varepsilon}
= \{\Pi,h_{\mathcal{O}_{(\mu,a)}}\}_{-}|_{\mathcal{O}_{(\mu,a)}}\cdot \bar{\varepsilon}(A,v,\Pi,\Gamma)\\
& = ( \frac{(I_2-I_3)\bar{\varepsilon}_2\bar{\varepsilon}_3}{I_2I_3}+
mgh(\bar{\varepsilon}_5\chi_3-\bar{\varepsilon}_6\chi_2),\;\;
\frac{(I_3-I_1)\bar{\varepsilon}_3\bar{\varepsilon}_1}{I_3I_1} +
mgh(\bar{\varepsilon}_6\chi_1-\bar{\varepsilon}_4\chi_3),\\
& \;\;\;\;\;\;
\frac{(I_1-I_2)\bar{\varepsilon}_1\bar{\varepsilon}_2}{I_1I_2}
 + mgh(\bar{\varepsilon}_4\chi_2-\bar{\varepsilon}_5\chi_1) ),
\end{align*}
and
\begin{align*}
& X_{h_{\mathcal{O}_{(\mu,a)}}}(\Gamma) \cdot
\bar{\varepsilon}
= \{\Gamma,h_{\mathcal{O}_{(\mu,a)}}\}_{-}|_{\mathcal{O}_{(\mu,a)}}\cdot \bar{\varepsilon}(A,v,\Pi,\Gamma)\\
& = ( \frac{I_2\bar{\varepsilon}_5\bar{\varepsilon}_3-
I_3\bar{\varepsilon}_6\bar{\varepsilon}_2}{I_2I_3}, \;\;
\frac{I_3\bar{\varepsilon}_6\bar{\varepsilon}_1-
I_1\bar{\varepsilon}_4\bar{\varepsilon}_3}{I_3I_1}, \;\;
\frac{I_1\bar{\varepsilon}_4\bar{\varepsilon}_2-
I_2\bar{\varepsilon}_5\bar{\varepsilon}_1}{I_1I_2} ).
\end{align*}
On the other hand, note that
\begin{align*}
X_H(\Pi)^\varepsilon & =T\pi_{\textmd{SO}(3)}\cdot X_H(\Pi)\cdot\varepsilon= X_H(\Pi)\cdot\varepsilon \\
& = ( \frac{(I_2-I_3)\varepsilon_8\varepsilon_9}{I_2I_3}+
mgh(\varepsilon_{11}\chi_3-\varepsilon_{12}\chi_2),\;\;
\frac{(I_3-I_1)\varepsilon_9\varepsilon_7}{I_3I_1} +
mgh(\varepsilon_{12}\chi_1-\varepsilon_{10}\chi_3),\\
& \;\;\;\;\;\;
\frac{(I_1-I_2)\varepsilon_7\varepsilon_8}{I_1I_2}
 + mgh(\varepsilon_{10}\chi_2-\varepsilon_{11}\chi_1) ),
\end{align*}
and
\begin{align*}
X_H(\Gamma)^\varepsilon & =T\pi_{\textmd{SO}(3)}\cdot X_H(\Gamma)\cdot\varepsilon= X_H(\Gamma)\cdot\varepsilon \\
& = ( \frac{I_2\varepsilon_{11}\varepsilon_9-
I_3\varepsilon_{12}\varepsilon_8}{I_2I_3}, \;\;
\frac{I_3\varepsilon_{12}\varepsilon_7-
I_1\varepsilon_{10}\varepsilon_9}{I_3I_1}, \;\;
\frac{I_1\varepsilon_{10}\varepsilon_8-
I_2\varepsilon_{11}\varepsilon_7}{I_1I_2} ),
\end{align*}
then we have that
\begin{align*}
T\bar{\gamma}\cdot X_H(\Pi)^\varepsilon
& =( \frac{(I_2-I_3)\bar{\gamma}_2\bar{\gamma}_3}{I_2I_3}+
mgh(\bar{\gamma}_5\chi_3-\bar{\gamma}_6\chi_2),\;\;
\frac{(I_3-I_1)\bar{\gamma}_3\bar{\gamma}_1}{I_3I_1} +
mgh(\bar{\gamma}_6\chi_1-\bar{\gamma}_4\chi_3),\\
& \;\;\;\;\;\;
\frac{(I_1-I_2)\bar{\gamma}_1\bar{\gamma}_2}{I_1I_2}
 + mgh(\bar{\gamma}_4\chi_2-\bar{\gamma}_5\chi_1)),
\end{align*}
and
\begin{align*}
T\bar{\gamma}\cdot X_H(\Gamma)^\varepsilon
& = ( \frac{I_2\bar{\gamma}_5\bar{\gamma}_3-
I_3\bar{\gamma}_6\bar{\gamma}_2}{I_2I_3}, \;\;
\frac{I_3\bar{\gamma}_6\bar{\gamma}_1-
I_1\bar{\gamma}_4\bar{\gamma}_3}{I_3I_1}, \;\;
\frac{I_1\bar{\gamma}_4\bar{\gamma}_2-
I_2\bar{\gamma}_5\bar{\gamma}_1}{I_1I_2} ),
\end{align*}
as well as
\begin{align*}
T\bar{\lambda}\cdot X_H(\Pi) \cdot \varepsilon
& =( \frac{(I_2-I_3)\bar{\lambda}_2\bar{\lambda}_3}{I_2I_3}+
mgh(\bar{\lambda}_5\chi_3-\bar{\lambda}_6\chi_2),\;\;
\frac{(I_3-I_1)\bar{\lambda}_3\bar{\lambda}_1}{I_3I_1} +
mgh(\bar{\lambda}_6\chi_1-\bar{\lambda}_4\chi_3),\\
& \;\;\;\;\;\;
\frac{(I_1-I_2)\bar{\lambda}_1\bar{\lambda}_2}{I_1I_2}
 + mgh(\bar{\lambda}_4\chi_2-\bar{\lambda}_5\chi_1)),
\end{align*}
and
\begin{align*}
T\bar{\lambda}\cdot X_H(\Gamma) \cdot \varepsilon
& = ( \frac{I_2\bar{\lambda}_5\bar{\lambda}_3-
I_3\bar{\lambda}_6\bar{\lambda}_2}{I_2I_3}, \;\;
\frac{I_3\bar{\lambda}_6\bar{\lambda}_1-
I_1\bar{\lambda}_4\bar{\lambda}_3}{I_3I_1}, \;\;
\frac{I_1\bar{\lambda}_4\bar{\lambda}_2-
I_2\bar{\lambda}_5\bar{\lambda}_1}{I_1I_2} ).
\end{align*}
Thus, when we choose that $(\Pi,\Gamma)=(\Pi_1,\Pi_2,\Pi_3,\Gamma_1,\Gamma_2,\Gamma_3)\in \mathcal{O}_{(\mu,a)}, $
and $\Pi=(\bar{\gamma}_1,\bar{\gamma}_2,\bar{\gamma}_3)=
(\bar{\varepsilon}_1,\bar{\varepsilon}_2,\bar{\varepsilon}_3)=
(\bar{\lambda}_1,\bar{\lambda}_2,\bar{\lambda}_3), $ and
$\Gamma=(\bar{\gamma}_4,\bar{\gamma}_5,\bar{\gamma}_6)=
(\bar{\varepsilon}_4,\bar{\varepsilon}_5,\bar{\varepsilon}_6)=
(\bar{\lambda}_4,\bar{\lambda}_5,\bar{\lambda}_6), $ we must have that
\begin{align*}
& T\bar{\gamma}\cdot X_H(\Pi)^\varepsilon =X_{h_{\mathcal{O}_{(\mu,a)}}}(\Pi) \cdot \bar{\varepsilon}
=T\bar{\lambda}\cdot X_H(\Pi) \cdot \varepsilon, \\
& T\bar{\gamma}\cdot X_H(\Gamma)^\varepsilon=X_{h_{\mathcal{O}_{(\mu,a)}}}(\Gamma) \cdot \bar{\varepsilon}
=T\bar{\lambda}\cdot X_H(\Gamma) \cdot \varepsilon.
\end{align*}
Since the map $\varepsilon: T^* \textmd{SE}(3)
\rightarrow T^* \textmd{SE}(3)$ is symplectic, then
$T\bar{\varepsilon}\cdot X_{h_{\mathcal{O}_{(\mu,a)}} \cdot \bar{\varepsilon}}
=X_{h_{\mathcal{O}_{(\mu,a)}}} \cdot \bar{\varepsilon}. $
Thus, in this case, we must have that
$\varepsilon$ and $\bar{\varepsilon} $ are the solution of the Type II of Lie-Poisson
Hamilton-Jacobi equation
$T\bar{\gamma}\cdot X_H^\varepsilon= X_{h_{\mathcal{O}_{(\mu,a)}}}\cdot \bar{\varepsilon}, $
for the regular point reduced heavy top system
$(\mathcal{O}_{(\mu,a)},\omega^{-}_{\mathcal{O}_{(\mu,a)}},h_{\mathcal{O}_{(\mu,a)}})$, if and only if they satisfy
the equation $T\bar{\varepsilon}\cdot(X_{h_{\mathcal{O}_{(\mu,a)}} \cdot \bar{\varepsilon}})
= T\bar{\lambda}\cdot X_H \cdot\varepsilon. $\\

\begin{remark}
It is worthy of note that, if the one-form $\gamma: \textmd{SE}(3) \rightarrow T^*
\textmd{SE}(3)$ is determined by a generating function of a symplectic map, then it is a
solution of classical Hamilton-Jacobi equation $X_H\cdot \gamma=0.$
From Remark 3.7 we know that
the classical Lie-Poisson Hamilton-Jacobi equation for the regular point reduced
heavy top system  $(\mathcal{O}_{(\mu,a)} ,
\omega_{\mathcal{O}_{(\mu,a)}}^{-},h_{\mathcal{O}_{(\mu,a)}} )$ is given by
$X_{h_{\mathcal{O}_{(\mu,a)}}}\cdot \bar{\gamma}=0, $ that is,
\begin{equation*}
 \left\{\begin{aligned}
 & \frac{(I_2-I_3)\bar{\gamma}_2\bar{\gamma}_3}{I_2I_3}+
mgh(\bar{\gamma}_5\chi_3-\bar{\gamma}_6\chi_2)=0, \\
 & \frac{(I_3-I_1)\bar{\gamma}_3\bar{\gamma}_1}{I_3I_1} +
mgh(\bar{\gamma}_6\chi_1-\bar{\gamma}_4\chi_3)=0, \\
 & \frac{(I_1-I_2)\bar{\gamma}_1\bar{\gamma}_2}{I_1I_2}
 + mgh(\bar{\gamma}_4\chi_2-\bar{\gamma}_5\chi_1)=0,\\
 & \frac{I_2\bar{\gamma}_5\bar{\gamma}_3-
I_3\bar{\gamma}_6\bar{\gamma}_2}{I_2I_3}=0, \\
 & \frac{I_3\bar{\gamma}_6\bar{\gamma}_1-
I_1\bar{\gamma}_4\bar{\gamma}_3}{I_3I_1}=0, \\
 & \frac{I_1\bar{\gamma}_4\bar{\gamma}_2-
I_2\bar{\gamma}_5\bar{\gamma}_1}{I_1I_2}=0,
\end{aligned} \right. \label{5.8}
\end{equation*}
where $\bar{\gamma}=\pi_{(\mu,a)}(\gamma): \textmd{SE}(3) \rightarrow
\mathcal{O}_{(\mu,a)}. $
\end{remark}

To sum up the above discussion, we have the following proposition.
\begin{proposition}
If the 4-tuple $(T^\ast \textmd{SE}(3),\textmd{SE}(3),\omega, H)$ is a
regular point reducible Hamiltonian system, then for a point $(\mu,a)\in
\mathfrak{se}^\ast(3)$, the regular value of the momentum map
$\mathbf{J}: T^\ast \textmd{SE}(3)\to \mathfrak{se}^\ast(3)$, the
Marsden-Weinstein reduced system is 3-tuple $(\mathcal{O}_{(\mu,a)},
\omega^{-}_{\mathcal{O}_{(\mu,a)}},h_{\mathcal{O}_{(\mu,a)}})$,
where $\mathcal{O}_{(\mu,a)} \subset \mathfrak{se}^\ast(3)$ is the
coadjoint orbit, $\omega_{\mathcal{O}_{(\mu,a)}}$ is the orbit
symplectic form on $\mathcal{O}_{(\mu,a)}$, which is induced by the
heavy top Poisson bracket on $\mathfrak{se}^\ast(3)$,
$h_{\mathcal{O}_{(\mu,a)}}(\Pi,\Gamma)\cdot \pi_{\mathcal{O}_{(\mu,a)}}=H(A,v,\Pi,\Gamma)
=H(A,v,\Pi,\Gamma)|_{\mathcal{O}_{(\mu,a)}}$. Assume that $\gamma:
\textmd{SE}(3) \rightarrow T^* \textmd{SE}(3)$ is an one-form
on $\textmd{SE}(3)$,
and $\lambda=\gamma \cdot \pi_{\textmd{SE}(3)}:
T^* \textmd{SE}(3) \rightarrow T^* \textmd{SE}(3), $ and $\varepsilon:
T^* \textmd{SE}(3) \rightarrow T^* \textmd{SE}(3), $ is a $\textmd{SE}(3)_\mu$-invariant symplectic map.
Denote by
$X_H^\gamma = T\pi_{\textmd{SE}(3)}\cdot X_H \cdot \gamma$, and
$X_H^\varepsilon = T\pi_{\textmd{SE}(3)}\cdot X_H \cdot \varepsilon$.
Moreover, assume that
$\textmd{Im}(\gamma)\subset \mathbf{J}^{-1}(\mu,a), $ and it is
$\textmd{SE}(3)_{(\mu,a)}$-invariant, where
$\textmd{SE}(3)_{(\mu,a)}$ is the isotropy subgroup of coadjoint
$\textmd{SE}(3)$-action at the point
$(\mu,a)\in\mathfrak{se}^\ast(3)$, and $\varepsilon(\mathbf{J}^{-1}(\mu,a))\subset \mathbf{J}^{-1}(\mu,a). $
Denote by
$\bar{\gamma}=\pi_{(\mu,a)}(\gamma): \textmd{SE}(3) \rightarrow
\mathcal{O}_{(\mu,a)}, $ and $\bar{\lambda}=\pi_{(\mu,a)}(\lambda):
\mathbf{J}^{-1}(\mu,a) \rightarrow \mathcal{O}_{(\mu,a)}, $ and
$\bar{\varepsilon}=\pi_{(\mu,a)}(\varepsilon): \mathbf{J}^{-1}(\mu,a)\rightarrow \mathcal{O}_{(\mu,a)}. $
Then the following two assertions hold:\\
\noindent $(\mathrm{i})$
If the one-form $\gamma: \textmd{SE}(3) \rightarrow T^*\textmd{SE}(3) $ is closed with respect to
$T\pi_{\textmd{SE}(3)}: TT^* \textmd{SE}(3) \rightarrow T\textmd{SE}(3), $
then $\bar{\gamma}$ is a solution of the Type I of Lie-Poisson Hamilton-Jacobi equation
$T\bar{\gamma}\cdot X_H^\gamma= X_{h_{\mathcal{O}_{(\mu,a)}}}\cdot \bar{\gamma}; $\\
\noindent $(\mathrm{ii})$
The $\varepsilon$ and $\bar{\varepsilon} $ satisfy the Type II of Lie-Poisson Hamilton-Jacobi equation
$T\bar{\gamma}\cdot X_H^\varepsilon= X_{h_{\mathcal{O}_{(\mu,a)}}}\cdot \bar{\varepsilon}, $
if and only if they satisfy
the equation $T\bar{\varepsilon}\cdot(X_{h_{\mathcal{O}_{(\mu,a)}} \cdot \bar{\varepsilon}})
= T\bar{\lambda}\cdot X_H \cdot\varepsilon. $ \hskip 0.3cm $\blacksquare$
\end{proposition}

In the following we shall introduce briefly some topics in future
research. At first, we note that the theory of controlled mechanical
systems has formed an important subject in recent years. In
particular, in Marsden et al.\cite{mawazh10}, the authors set up the
regular reduction theory of RCH systems on a symplectic fiber
bundle, by using momentum map and the associated reduced symplectic
forms and from the viewpoint of completeness of regular symplectic
reduction, and some generalizations for Poisson structure are given
in Wang and Zhang \cite{wazh12} and Ratiu and Wang \cite{rawa12}.
Since the Hamilton-Jacobi theory is developed based on the
Hamiltonian picture of dynamics, it is natural idea to extend the
Hamilton-Jacobi theory to the (regular) controlled Hamiltonian
system and its a variety of the reduced systems, and it is also
possible to describe the relationship between the RCH-equivalence
for the controlled Hamiltonian systems and the solutions of the
corresponding Hamilton-Jacobi equations. Next, if the Hamiltonian
system we considered has nonholonomic constraints, in general, the
dynamical vector field of nonholonomic Hamiltonian system is not
Hamiltonian, however, it can be described by the dynamical vector
field of a distributional Hamiltonian system. Thus, it is possible
to set up the Hamilton-Jacobi theory for the nonholonomic
Hamiltonian system and the nonholonomic reducible Hamiltonian system
on a cotangent bundle by using the distributional Hamiltonian system
and the reduced distributional Hamiltonian system, see de Le\'{o}n
and Wang \cite{lewa15}. Finally, we also note that there have been a
lot of beautiful results of reduction theory of Hamiltonian systems
in celestial mechanics, hydrodynamics and plasma physics. So, it is
an important topic to study the application of reduction and
Hamilton-Jacobi theory in celestial mechanics, hydrodynamics and
plasma physics. These are our goals in future research.\\

\noindent{\bf Acknowledgments:} Especially grateful to
Professor Tudor S. Ratiu, Professor Manuel de Le\'{o}n,
Professor Arjan Van der Schaft and Professor Juan-Pablo Ortega for their
help and guiding in the study of geometric mechanics.
H. Wang's research was partially supported by Nankai University, 985 Project
and the Key Laboratory of Pure Mathematics and Combinatorics, Ministry of
Education, China.


\end{document}